\documentclass[preprint,5p,twocolumn]{elsarticle}
\pdfoutput=1 
\usepackage{epsf}
\usepackage{graphicx}
\usepackage{psfig}
\usepackage{dcolumn}
\usepackage{bm}
\usepackage{subfigure}
\usepackage{color}
\usepackage{amsthm}
\usepackage{amssymb}
\usepackage{amsmath}
\usepackage{pb-diagram}
\usepackage{algorithm}
\usepackage{algorithmic}
\usepackage{hyperref}

\newtheorem{thm}{Theorem}[section]

\newtheorem{defn}[thm]{Definition}

\newtheorem{rem}[thm]{Remark}

\newcommand{\setof}[1]{\left\{ {#1}\right\}}

\newcommand{\R}{{\mathbb{R}}}
\newcommand{\T}{{\mathbb{T}}}
\newcommand{\Z}{{\mathbb{Z}}}

\newcommand{\cR}{{\mathcal R}}

\newcommand{\sC}{{\mathsf C}}

\newcommand{\sR}{{\mathsf R}}

\def\setof#1{\left\{{#1}\right\}}
\def\ang#1{\langle {#1} \rangle}

\definecolor{gray85}{gray}{0.85} 
\definecolor{gray8}{gray}{0.8} 
\definecolor{gray7}{gray}{0.7} 
\definecolor{gray6}{gray}{0.6} 
\definecolor{gray5}{gray}{0.5} 
\definecolor{gray4}{gray}{0.4} 
\definecolor{gray35}{gray}{0.35} 

\newcommand{\pd}{{\mathsf{PD}}}
\newcommand{\per}{{\mathsf{Per}}}

\begin{document}

\begin{frontmatter}

\title{Analysis of Kolmogorov Flow and Rayleigh-B\'enard Convection using Persistent Homology}

\author[ru]{Miroslav Kram\' ar}
\ead{miroslav@math.rutgers.edu}

\author[ru]{Rachel Levanger}
\ead{rachel@math.rutgers.edu}

\author[gt]{Jeffrey Tithof}
\ead{jtithof3@gatech.edu}

\author[gt]{Balachandra Suri}
\ead{balachandra@gatech.edu}

\author[vt]{Mu Xu}
\ead{xumu8621@vt.edu}

\author[vt]{Mark Paul}
\ead{mrp@vt.edu}

\author[gt]{Michael F. Schatz}
\ead{michael.schatz@physics.gatech.edu}

\author[ru]{Konstantin Mischaikow}
\ead{mischaik@math.rutgers.edu}

\address[ru]{Department of Mathematics,
Hill Center-Busch Campus,
Rutgers University,
110 Frelingheusen Rd,
Piscataway, NJ  08854-8019, USA}
\address[gt]{Center for Nonlinear Science and School of Physics,
Georgia Institute of Technology,
Atlanta, GA 30332-0430, USA}
\address[vt]{Department of Mechanical Engineering,
Virginia Tech, Blacksburg, VA 24061, USA}

\begin{abstract}
We use persistent homology to build a quantitative understanding of large complex 
systems that are driven far-from-equilibrium; in particular, we analyze image time 
series of flow field patterns from numerical simulations of two important problems in 
fluid dynamics: Kolmogorov flow and Rayleigh-B\'enard convection.  
For each image we compute a persistence diagram to yield a reduced description of 
the flow field; by applying different metrics to the space of persistence diagrams, we relate characteristic 
features in persistence diagrams to the geometry of the corresponding flow patterns.  
We also examine the dynamics of the flow patterns by a second application of persistent 
homology to the time series of persistence diagrams.  We demonstrate that persistent
homology provides an effective method both for quotienting  out symmetries in families of solutions and for identifying multiscale recurrent dynamics.  Our approach is quite general and it is anticipated to be applicable to a broad range of open problems exhibiting complex spatio-temporal behavior.
\end{abstract}

\end{frontmatter}

\section{Introduction}

We introduce new mathematical techniques for analyzing complex  spatiotemporal nonlinear dynamics and demonstrate their efficacy in problems from two different paradigms in hydrodynamics.
Our approach employs methods from algebraic topology; earlier efforts have shown that computing the homology of topological spaces associated to scalar or vector fields generated by complex systems can provide new insights into dynamics \cite{ gameiro:kalies:mischaikow, gameiro:mischaikow:wanner, kramar:goullet:kondic:mischaikow:2014,krishnan:kurtuldu:schatz:gameiro:mischaikow:madruga, kurtuldu:mischaikow:schatz:11b, kurtuldu:mischaikow:schatz:11a}. We extend prior work by using a relatively new tool called persistent homology \cite{edelsbrunner:harer, carlsson, weinberger}.  

Complex spatiotemporal systems often exhibit complicated pattern evolution. The patterns are given by scalar or vector fields representing the state of the system under study. Persistent homology can be viewed as a map $\pd$ that assigns to every field a collection of points in $\R^2$, called a {\it persistence diagram}.  For a given scalar field $f : D \to R$, the points in the persistence diagram $\pd(f)$ encode geometric features of the sub-level sets $\sC(f,\theta) = \setof{x \in D \mid f(x) < \theta}$ for all values of $\theta$. A feature encoded by the point $(\theta_b, \theta_d) \in \pd(f)$  appears in  $\sC(f,\theta_b)$ for the first time and disappears in $\sC(f,\theta_d))$. Therefore, $\theta_b$ and $\theta_d$ are called birth and death coordinates of this feature. The  lifespan $\theta_d - \theta_b > 0$ indicates the prominence of the feature. In particular, features with long lifespans are considered important and features with short lifespans are often associated with noise. Thus, the persistence diagram is a highly simplified representation
of the field generating the pattern.

The space of all persistence diagrams, $\per$, can be endowed with a  variety of metrics under which $\pd$ is a continuous function.  This has several important implications that we exploit in this paper. First, continuity implies that small changes in the field pattern, e.g.\ bounded errors associated with measurements or numerical approximations, lead to small changes in the persistence diagrams.  Second, by using different metrics, we can vary our focus of interest between larger and smaller changes in the persistence diagrams.  Moreover, by comparing different metrics, we can infer if the changes in a pattern  affect  geometric features with longer or shorter life spans. Finally,  since, applying the map $\pd$ to a time series of patterns produces a time series in $\per$, the distance between the consecutive data points  in $\per$ can be used to quantify the average rate at which the geometry of the patterns is changing. 

As mentioned above, the dynamics of spatiotemporal systems are characterized by the time-evolution of the patterns corresponding to the fields generated by the system.  However,  capturing these vector fields, either experimentally or numerically, results in  multi-scale high dimensional data sets. In order to efficiently analyze these data sets, a dimension reduction must be performed. We use persistent homology to  perform nonlinear dimension reduction from a time series of patterns to a time series of persistence diagrams.  We show  that this reduction  can cope with redundancies introduced by symmetries (both discrete and continuous) present in the system. In particular, this approach directly quotients out symmetries and, thereby, permits easy identification of solutions that lie on a group orbit.

Separately, we also apply persistent homology to extract information about dynamical structures in the reduced data.  Characterizing dynamics  in the space of persistence diagrams cannot be done using conventional methods (e.g., time delay embeddings), since choosing a coordinate system in $\per$ is currently an open problem.  However, since $\per$ is a metric space, the geometry of the point cloud $X$, generated by the time series of the reduced data, is encoded  by a scalar field which assigns to each point in $\per$ its distance to $X$. We show how persistent homology may be applied to describe dynamics by characterizing the geometry of $X$. 

An outline of the paper is as follows. In Section~\ref{sec:systems} we present a brief overview of the two fluid flows examined in this paper: (1) Kolmogorov flow and (2) Rayleigh-B\'{e}nard convection.  We note here, for emphasis, that while persistent homology can be applied to vector fields,  it will be sufficient for this paper to focus on  scalar fields drawn from these systems (specifically, one component of the vorticity field for Kolmogorov flow, and the  temperature field for Rayleigh-B\'{e}nard convection).   

In Section~\ref{sec:PH*} we discuss  key issues related to the application of persistent homology. By now, the mathematical theory of persistent homology is well developed. Therefore, our main emphasis is on the computational aspect of passing from the data to the persistence diagrams. Section~\ref{sec:PH}  describes  the correspondence between the geometric features of a scalar field and the points in its corresponding persistence diagram. Section~\ref{sec:SpaceOfPD} discusses the structure of the space $\per$ and the properties of the associated metrics.

In Sections~\ref{sec:Distances} and \ref{sec:AnalyzingPointCloud} we discuss how these metrics can be used to analyze dynamics.  First, we interpret  distance between the persistence diagrams representing the consecutive data points in the time series as a rate at which  geometry of the corresponding scalar fields is changing. Second, we motivate and explain the procedure for extracting the geometric structure of the point cloud in $\per$.

We close  the paper by applying the developed techniques  to the following problems. In Section~\ref{sec:FixedPoints}, we identify distinct classes of symmetry-related equilibria for Kolomogorov flow.  In Section~\ref{sec:PeriodicOrbit}, we show that a relative periodic orbit for Kolmogorov flow  collapses to a closed loop in $\per$. Finally, in Section~\ref{sec:RBC}, we  deal with identifying recurrent dynamics that occur on different time scales in our study of Rayleigh-B\'{e}nard convection flow.

\section{The Systems to be Studied}
\label{sec:systems}
\subsection{Kolmogorov Flow} 
\label{Intro_Kolmogorov}

For the study of turbulence in two dimensions, Kolmogorov proposed a model flow where the two-dimensional (2D) velocity field  ${\bf u}(x,y,t)$ is given by 

\begin{eqnarray}\label{eq:q2dns}
\frac{\partial\bf{u}}{\partial t} + \beta {\bf u}\cdot{\bf{\nabla}}{\bf u} &=& -\frac{1}{\rho}{\bf{\nabla}}p + {\nu} {\bf\nabla}^2 {\bf u} -\alpha {\bf{u}} + {\bf f}\\
{\bf \nabla} \cdot {\bf u} &=& 0 \nonumber
\end{eqnarray}
(with $\beta = 1$ and $\alpha = 0$), where $p(x,y)$ is the pressure field, $\nu$ is the kinematic viscosity, $\rho$ is fluid density, and 
${\bf f} = \chi\sin(\kappa y) \hat{{\bf x}}$ is the forcing that drives the flow \citep{arnold_1960}.  Laboratory experiments in electromagnetically-driven shallow layers of electrolyte can exhibit flow dynamics that  are well-described by Equations (\ref{eq:q2dns}) with appropriate choices of $\beta$ and $\alpha$ to capture  three-dimensional effects, which are commonly present in experiments \cite{suri_2014}.   In this paper, we refer to  all models described by Equations (\ref{eq:q2dns}) (including experimentally-realistic versions) as Kolmogorov flows.

It is convenient to use  the vorticity-stream function formulation \citep{panton_2006} to study Kolmogorov flow analytically and numerically. Equations (\ref{eq:q2dns}), written in terms of the z-component of the vorticity field $\omega = (\nabla \times {\bf u}) \cdot \hat{\bf{k}}$,  a scalar field, take the form

\begin{equation}\label{eq:q2dvor_nd}
\frac{\partial\omega}{\partial t} + \beta{\bf u}\cdot{\bf \nabla}\omega=
\nu \nabla^2\omega - \alpha\omega + \chi\kappa\cos(\kappa y).
\end{equation}

For the current study, we choose $\beta = 0.83$, $\nu = 3.26 \times 10^{-6}$ m$^2$/s, $\alpha = 0.063$ s$^{-1}$, $\rho = 959$ kg/m$^3$, and $\lambda = 2\pi/\kappa =  0.0254$ m.    We express the strength of the forcing in terms of a  non-dimensional parameter, the Reynolds number $Re = \sqrt{\frac{\lambda^3\chi}{8 \nu^2}}$.

Equation (\ref{eq:q2dvor_nd}) is solved numerically by using a pseudo-spectral method \citep{mitchell_2013}, assuming periodic boundary conditions in both $x$ and $y$ directions, i.e., $\omega(x,y) = {\omega}(x+L_x,y) = {\omega}(x,y+L_y)$, where  $L_x = 0.085$ m and $L_y = 4\lambda = 0.1016$ m are the dimensions of the domain in the $x$ and $y$ directions, respectively.  

It is important to note that Equation (\ref{eq:q2dvor_nd}), with periodic boundary conditions, is invariant under any combination of three distinct coordinate transformations: (1) a translation along $x$: $\mathcal{T}_{\delta x}(x,y)=(x+\delta x,y)$, $\delta x \in [0,L_x]$; (2) a rotation by $\pi$: ${\mathcal R}(x,y)=(-x,-y)$; and (3) a reflection and a shift: $\mathcal{D}(x,y)=(-x,y+\lambda/2)$.  Because of these symmetries, each particular solution to Equation (\ref{eq:q2dvor_nd}) generates a set of solutions which are dynamically equivalent.  Physically, invariance under continuous translation leads to the existence of relative equilibria (REQ) and relative periodic orbit (RPO) solutions, in addition to equilibria (EQ) and periodic orbit (PO) solutions.

\begin{figure}[t!]
\centering
\subfigure[]{\includegraphics[width=2.5in]{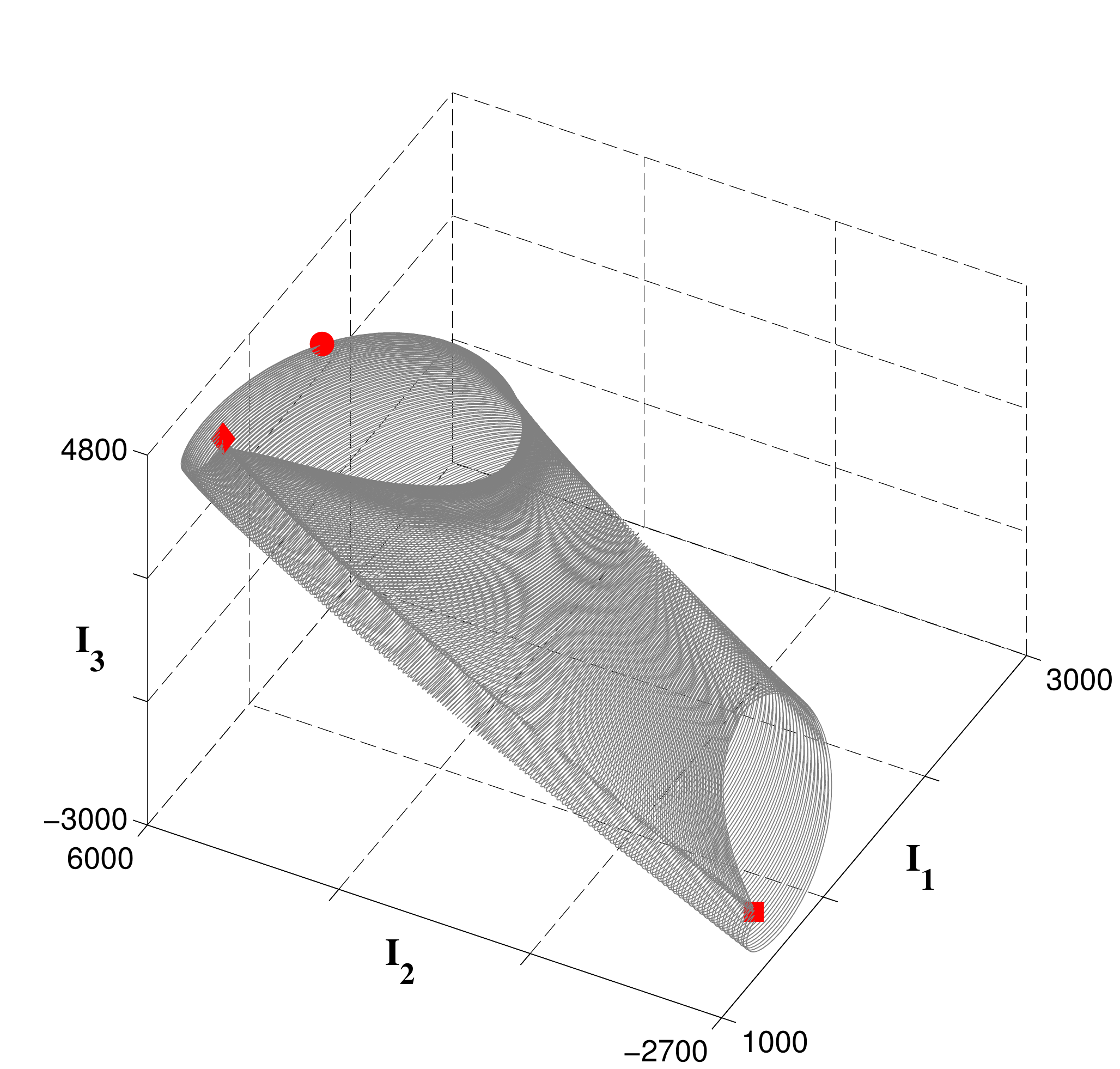}}
\subfigure[]{\includegraphics[width=2.5in]{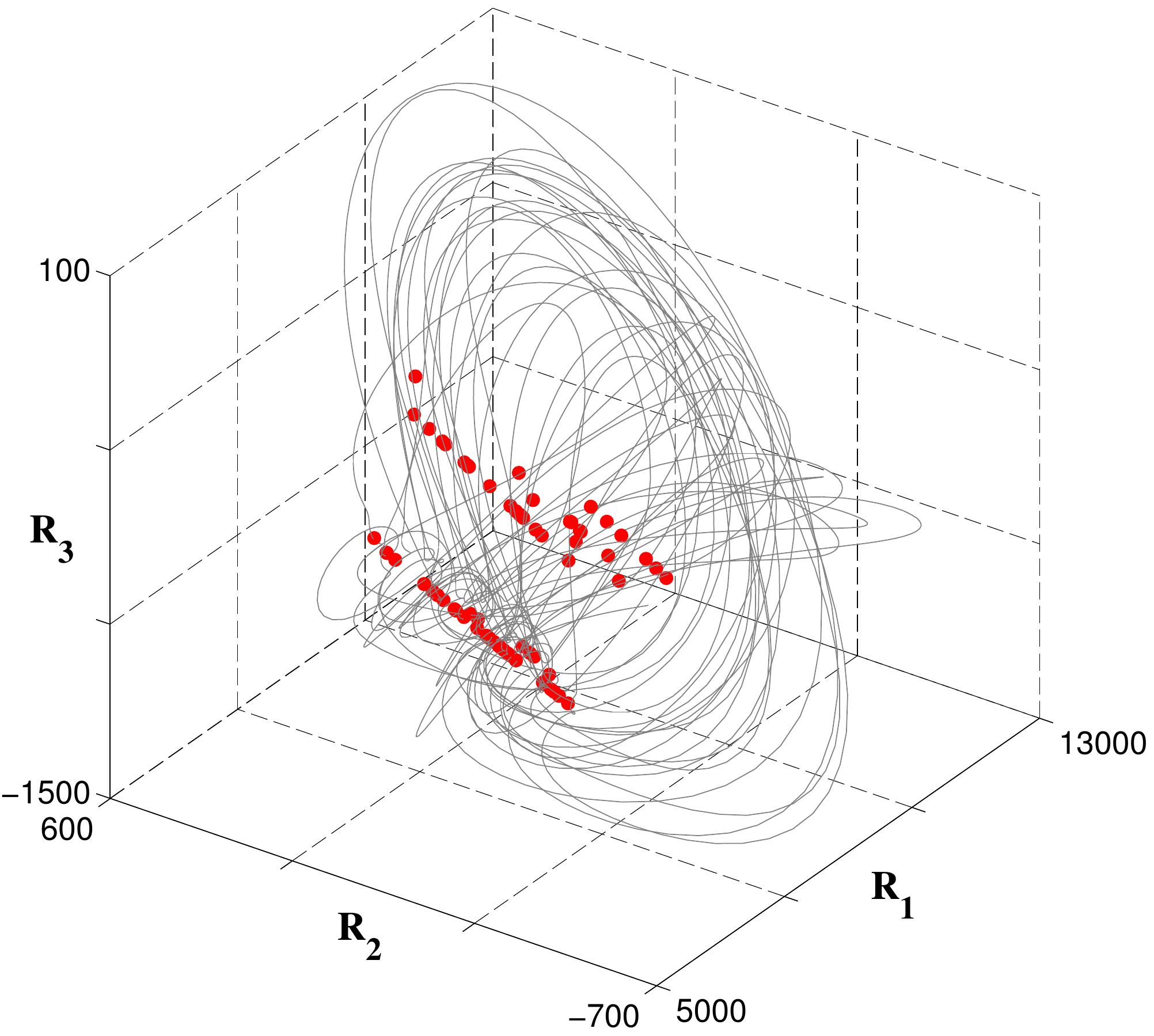}}
\caption{(a) Three-dimensional projections of a stable RPO at $Re=25.43$ from the Kolmogorov flow using the imaginary part of the three dominant Fourier modes, $I_1$, $I_2$, and $I_3$. The gray line indicates the evolution of a RPO; three snapshots sampled from that orbit are indicated by a red diamond, a red circle, and a red square, which are analyzed below. (b) Three-dimensional projections of a turbulent trajectory, at $Re=26.43$, using the real parts of the three dominant Fourier modes,  $R_1$, $R_2$, and $R_3$. The gray line indicates the chaotic evolution of the flow, which is influenced by the presence of unstable fixed points, indicated by red circles, which are also analyzed below.}
\label{fig:projections}
\end{figure}

For $Re = 25.43$, the flow is characterized by a steady RPO; Figure \ref{fig:projections}(a) shows a projection, plotted using the three dominant Fourier modes of this RPO. The RPO has a period 34.78 seconds and a drift speed $1.354 \times 10^{-6}$ m/s. The tunnel-like structure is a result of the periodic motion superposed over the slow drift along the $x$-direction.  For larger forcing ($Re =  26.43$), the flow becomes weakly turbulent, as can be seen from the Fourier projections in Figure \ref{fig:projections}(b). The turbulent dynamics in this regime are of great interest as the flow explores a region of the state space which contains ``weakly" unstable EQ, PO, REQ, and RPO solutions. Recent theoretical advances have shown that the identification of these solutions could aid the understanding of weakly turbulent dynamics \citep{chandler_2013}. For instance, if the turbulent trajectory is close to an EQ solution ($\omega_0$), which is characterized by $\partial \omega_0 / \partial t = 0$, we would expect the instantaneous rate of change of $\omega$ to be relatively small, i.e., $\partial \omega/\partial t \approx 0$. Similarly, a close pass to a PO solution would mean $\omega(t+T) \approx \omega(t)$, where $T$ is the period of the PO that is guiding the dynamics of the turbulent trajectory. The turbulent trajectory depicted in Figure \ref{fig:projections}(b) passes close to unstable EQ and REQ solutions which are indicated by the red dots.

\begin{figure}
\centering
\subfigure[]{\raisebox{0.02in}{\includegraphics[width=1.25in]{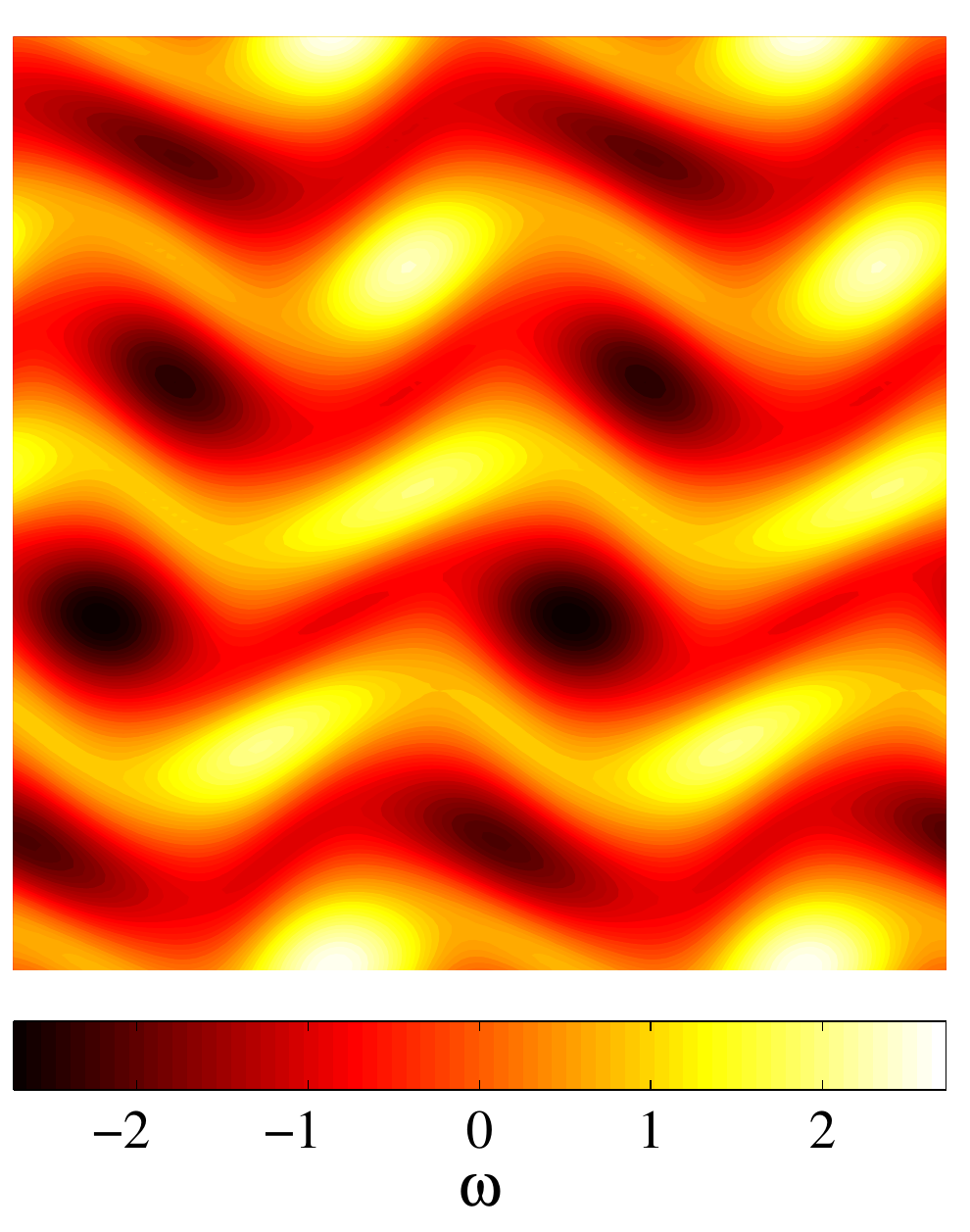}}}
\qquad
\subfigure[]{\includegraphics[width=1.75in]{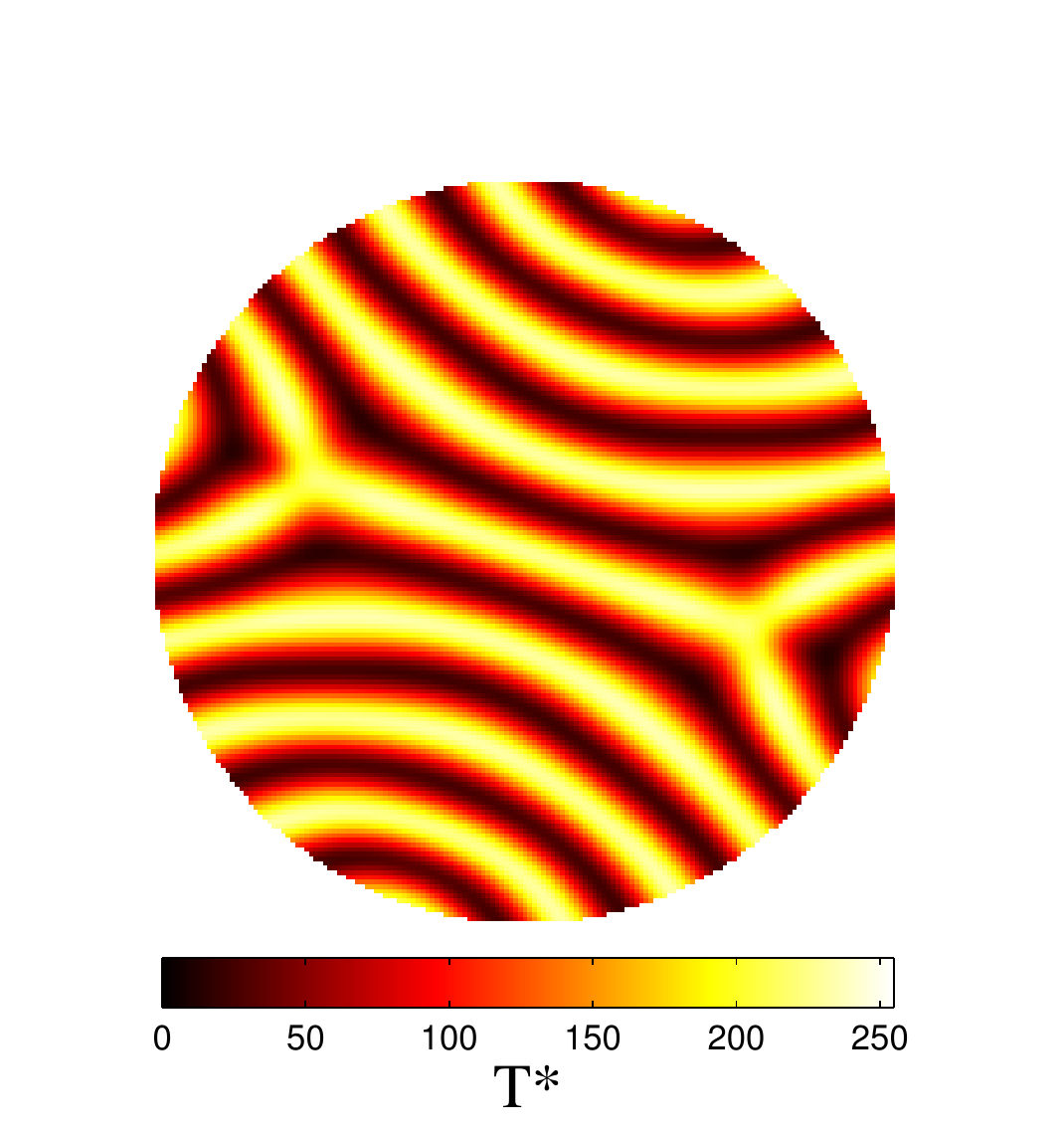}}
\caption{(a) A snapshot of the $z$-component of the vorticity field $\omega$ for Kolmogorov flow from the stable relative periodic orbit found at $Re=25.43$. (b) A snapshot of the renormalized 8-bit mid-plane temperature field $T^*$ for Rayleigh-B\'{e}nard convection from the stable almost-periodic orbit found at $Ra=3000$ and $Pr=1$.}
\label{fig:KLF}
\end{figure}

\subsection {Rayleigh-B\'{e}nard Convection}
\label{subsec:RBC}

Rayleigh-B\'enard convection is a canonical pattern forming system that has been used to gain many new fundamental insights into the spatiotemporal dynamics of systems that are driven far-from-equilibrium~\cite{cross:1993,bodenschatz:2000}.  Rayleigh-B\'enard convection is the buoyancy driven fluid flow that occurs when a shallow layer of fluid is heated uniformly from below in a gravitational field.  The dynamics are governed by the Boussinesq equations,
\begin{eqnarray}
	Pr^{-1} \left(\frac{\partial \mathbf{u}}{\partial t} + \mathbf{u}\cdot\mathbf{\nabla} \mathbf{u} \right)  &=& -\mathbf{\nabla} p + \nabla^{2}\mathbf{u}+ Ra T \hat{\mathbf{z}},  \label{eq:rbc1}\\
	\frac{\partial T}{\partial t} + \mathbf{u} \cdot \mathbf{\nabla} T &=& \nabla^{2} T,  \label{eq:rbc2}\\
	\mathbf{\nabla} \cdot \mathbf{u} &=& 0,
 \label{eq:rbc3}
\end{eqnarray}
where $\mathbf{u}(x,y,z,t)$ is a vector field of the fluid velocity, $p(x,y,z,t)$ is the pressure field, and $T(x,y,z,t)$ is the temperature field.  In our notation, the origin of the Cartesian coordinates $(x,y,z)$ at the center of the domain are at the lower heated plate where $\hat{\mathbf{z}}$ is a unit vector opposing gravity.  Equations~(\ref{eq:rbc1})-(\ref{eq:rbc3}) represent the conservation of momentum, energy, and mass, respectively.  The equations have been nondimensionalized using the vertical diffusion time of heat as the time scale, the layer depth as the length scale, and the constant temperature difference between the lower and upper plates as the temperature scale.

In our work, we consider Rayleigh-B\'enard convection in a shallow domain with a cylindrical cross-section.  The no-slip fluid boundary condition $\mathbf{u}=0$ is applied to all material surfaces.  The lower and upper plates are held at a constant temperature where $T(z=0)=1$ and $T(z=1)=0$, respectively.  The lateral sidewalls of the cylindrical container are assumed to be perfectly conducting, which yields $T(z)=1-z$.

The dynamics can be described using three non-dimensional parameters.  The Rayleigh number $Ra$ represents the ratio of buoyancy to viscous forces. At the critical value $Ra_c = 1708$, an infinite layer of fluid undergoes a bifurcation to straight and parallel convection rolls.  For increasing values of the Rayleigh number $Ra>Ra_c$, the dynamics become periodic, chaotic, and eventually turbulent.  The Prandtl number $Pr$ is the ratio of the momentum and thermal diffusivities.  For typical experiments using compressed gasses, $Pr \approx 1$.  Lastly, the aspect ratio of the cylindrical domain $\Gamma$ is the ratio of the domain's radius to its depth.  

We numerically integrate Eqs.~(\ref{eq:rbc1})-(\ref{eq:rbc3}) using a highly parallel spectral element algorithm that has been tailored for the study of convection (c.f.~\cite{paul:2003}).  Figure~\ref{fig:KLF}(b) shows a typical pattern from a numerical simulation of Rayleigh-B\'enard convection. In this simulation, $Ra=3000$, $Pr=1$, and the aspect ratio of the domain is $\Gamma=10$.  The numerical simulation is initiated from a field of small random perturbations to the temperature field and is integrated for long times.  Figure~\ref{fig:KLF}(b) illustrates the fluid temperature field at the horizontal mid-plane ($z=1/2$), where light is warm rising fluid and dark is cool falling fluid.  This image is a snap shot in time of a time-dependent pattern where the dynamics are nearly periodic in time.  The pattern shown does not include the region near the sidewall.  Specifically, a distance of one-layer depth from the lateral sidewall is not shown (this distance is approximately the width of a convection roll).  This is done to remove the complex fluid flow that occurs in the small region near the sidewalls to allow our diagnostics to focus upon the bulk patterns and dynamics (c.f.~\cite{bodenschatz:2000}).

\section{Persistent Homology}
\label{sec:PH*}

The aim of this paper is to introduce an approach for analyzing the dynamics of the pattern evolution in spatiotemporal systems. This is done in two steps. First, we perform nonlinear data reduction, and then we extract information about the dynamical structures from this reduced data. We formulate both of these tasks in terms of analyzing the structure of the sub-level sets of a scalar function $f : D \to \R$, where $D$ is a topological space. Tools from algebraic topology, homology in particular, are used to capture and quantify the geometry of the sub-level sets.

Recall that given any topological space $Z$, homology assigns to $Z$ a sequence of vector spaces $H_k(Z)$, $k=0,1,\ldots$. The dimension of $H_k(Z)$ is called the $k$-th {\em Betti number} and is denoted by $\beta_k(Z)$. Betti numbers provide geometric information about $X$:  $\beta_0(Z)$ is the number of connected components, or pieces, of $Z$; $\beta_1(Z)$ indicates the number of loops or tunnels in $Z$;  and $\beta_2(Z)$ is the number of cavities.

Our goal is to understand structure of the sub-level sets
\begin{equation}
\label{eq:sublevel}
\sC(f,\theta) = \setof{x\in D \mid f(x)\leq \theta},
\end{equation}
for all values of $\theta \in \R$.
As we vary $\theta$, the number of components, loops, and cavities in $\sC(f,\theta)$ changes, implying that $\beta_k(\sC(f,\theta))$, $k=0,1, 2$, also changes.  (See Section~\ref{sec:PH} for examples.)  What is remarkable is that, under very weak conditions, we can choose bases for the vector spaces $H_k(\sC(f,\theta))$ over all values of $\theta$ such that, given a basis element of $H_k(\sC(f,\theta))$, we can identify a unique value $\theta_b$ at which this basis element appears and a unique value $\theta_d$ at which this basis element disappears. We refer to $\theta_b$ as the {\em birth value}, $\theta_d$ as the {\em death value}, and the pair $(\theta_b,\theta_d)\in\R^2$ as a {\em persistence point} corresponding to the chosen basis element of $H_k(\sC(f,\theta))$. The difference $\theta_d-\theta_b$ is called the {\em life span} of the persistence point. Longer life spans are associated with geometric features that persist through larger variations of $\theta$, and persistence diagrams are a codification of this information.  Given a scalar field $f$, the set of associated \emph{persistence diagrams} are denoted by $\pd(f) = \setof{\pd_k(f)}$, where $\pd_k(f)$ consists of all  persistence points corresponding to the $k$-th level of homology  (keeping track of multiple copies of a single point), along with infinitely many points at each point along the diagonal $\theta_b=\theta_d$. The reason for the inclusion of the diagonal is made clear in Definition~\ref{defn:bottleneck}, when we define metrics on the space of persistence diagrams.

For the systems introduced in Section~\ref{sec:systems}, we first use persistent homology as a nonlinear data reduction method. For Kolmogorov flow we study the scalar field ${\omega\colon D\to \R}$, the $z$-component of the vorticity  field, while for Rayleigh-B\'enard convection we study the scalar field $T\colon D\to \R$, the temperature field at the mid-plane. It is important to note that the domains for these two scalar fields are different.
For Kolmogorov flow, the domain $D$ is a torus since we are using periodic boundary conditions, while for Rayleigh-B\'enard convection, $D$ is a disk. For the disk, we need only to concern ourselves with the vector spaces $H_k(\sC(\omega,\theta))$ for $k=0,1$. However, for the torus, the vector spaces $H_2(\sC(T,\theta))$ also need to be considered, since the torus encloses a three-dimensional cavity. In section~\ref{sec:PH}, we explain how the persistence diagrams $\pd(f)$ capture important information about the patterns given by the scalar fields $\omega$ and $T$.

The set of all persistence diagrams $\pd$ is a metric space, denoted by $\per$ (see Section~\ref{sec:SpaceOfPD}). Since we are studying the evolution of Kolmogorov flow and Rayleigh-B\'enard convection,  we have  time series of the vorticity $\setof{\omega_i}$ and temperature $\setof{T_i}$ fields, and, therefore, we have time series of persistence diagrams $\setof{\pd(\omega_i)}$ and $\setof{\pd(T_i)}$. We view each of these time series as a point cloud $X\subset \per$. To extract information about  dynamical structures present in the time series, we use persistent homology a second time to quantify the geometry associated with this point cloud. This is achieved by introducing a new scalar function $f \colon \per \to [0,\infty)$ that gives the distance from any point in $\per$ to the point cloud $X$ and is defined by
\begin{equation}
\label{eq:metricScalar}
f(x) := d(x,X) := \min_{x_i \in X}d(x,x_i),
\end{equation}
where $d$ is an appropriate metric  on the space of persistence diagrams.  The associated sub-level sets $\sC(f,\theta)$ are once again given by \eqref{eq:sublevel}.

To carry out the steps mentioned above requires the ability to compute the persistence diagrams $\pd(f)$. To do this, we need to calculate $H_k(\sC(f,\theta))$, which requires  a discrete representation of $\sC(f,\theta)$ called a {\em complex}. In the context of nonlinear data reduction, we make use of a {\em cubical} complex.  When  analyzing the geometry of the point cloud, we approximate  $\sC(f,\theta)$ using a {\em Vietoris-Rips} complex, which is a special form of a simplicial complex. This is a classical subject and thus there are a variety of references providing precise definitions of complexes, e.g. \cite{edelsbrunner:harer} for Vietoris-Rips complexes and \cite{kaczynski:mischaikow:mrozek} for cubical complexes, discussions of issues related to approximations \cite{kramar:goullet:kondic:mischaikow:2014}, and how one proceeds from a complex to computing persistent homology \cite{mischaikow:nanda, edelsbrunner:harer}. The homological computations in this paper were performed using the Perseus software \cite{perseus}.

The numerical data for the vorticity and the temperature fields is presented in the form of piecewise-constant functions defined on a rectangular lattice.  For Kolmogorov flow,  values of $\omega$ are reported in double precision. Recall that the vector spaces $H_k(\sC(\omega,\theta))$ can only change for $\theta \in \Theta$, where  $\Theta$ is the finite set of values that $\omega$ attains on the given lattice. Each of the sets $\sC(\omega,\theta)$ is a cubical complex, and we  use the Perseus software to compute the corresponding persistence diagrams  using only the values $\theta \in \Theta$. Numerical simulations for   Rayleigh-B\'{e}nard convection are carried out with high precision as well. However, keeping in mind our goal to compare the numerical  simulations  with  experimental data, we convert the temperature field to an $8$-bit temperature field $T^*$ (an integer-valued function with values between $0$ and $255$), which can be obtained experimentally.  Consequences of this rescaling are examined in Section~\ref{sec:Distances}.
 
\begin{figure}[t!]
\centering
\subfigure[]{\includegraphics[width=1.5in]{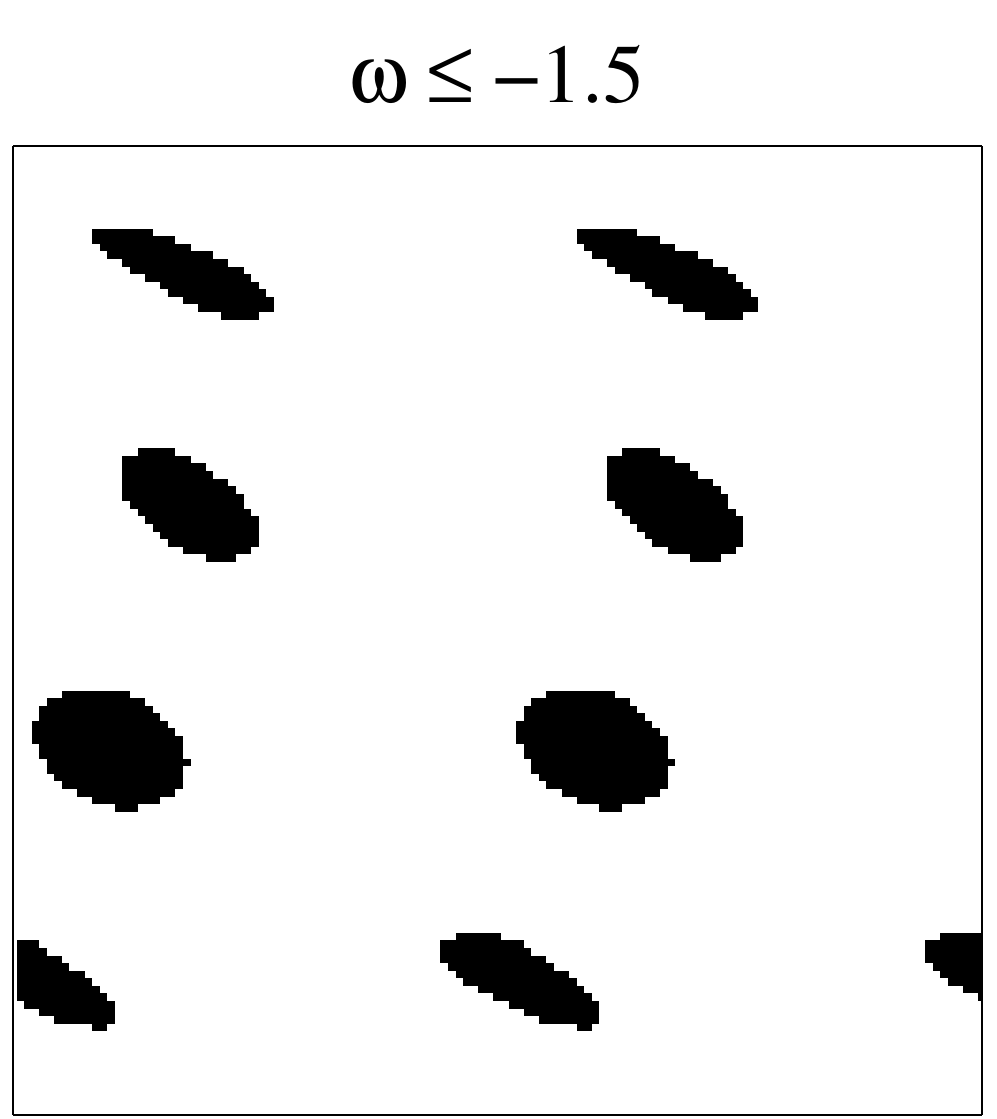}}\qquad
\subfigure[]{\includegraphics[width=1.5in]{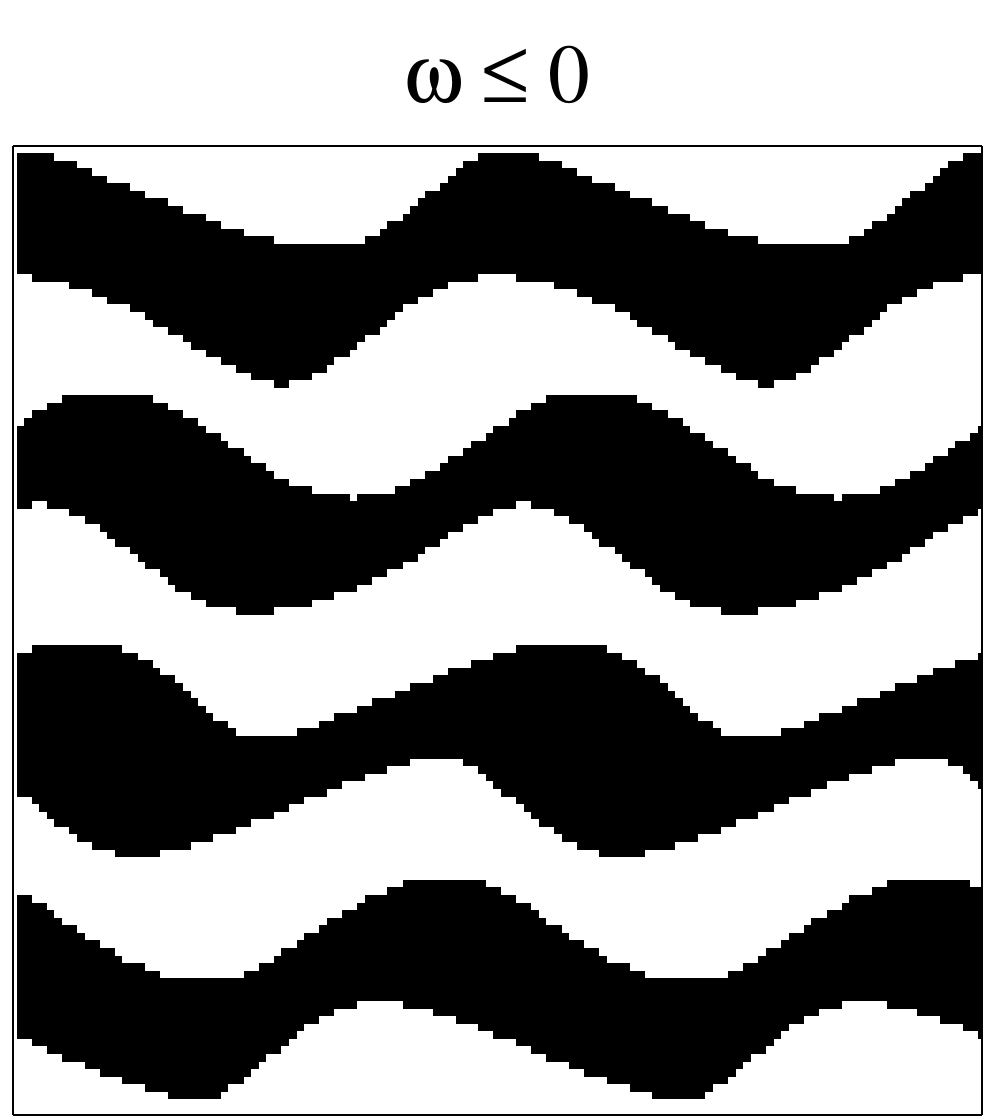}}\\
\subfigure[]{\includegraphics[width=1.5in]{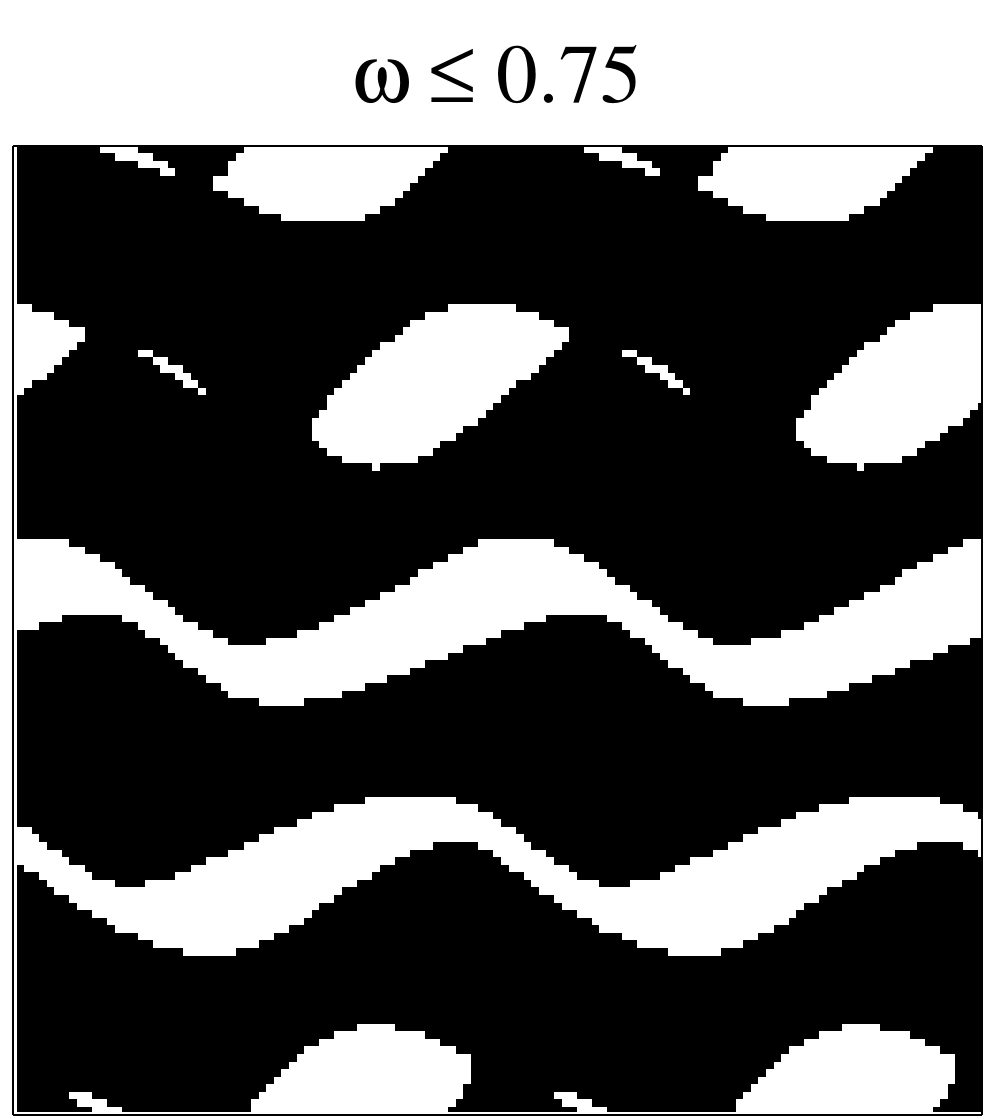}}\qquad
\subfigure[]{\includegraphics[width=1.5in]{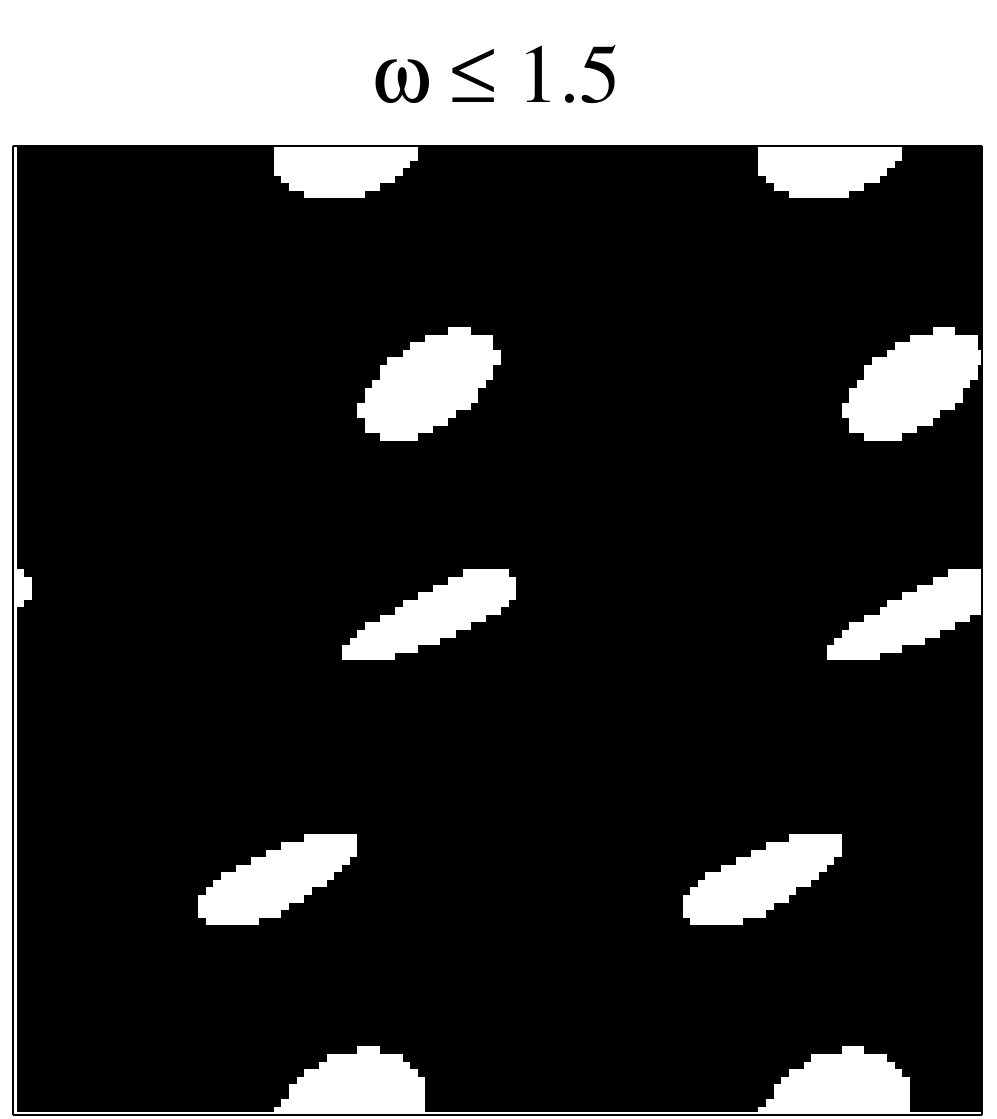}}\\
\subfigure[]{\includegraphics[width=1.6in]{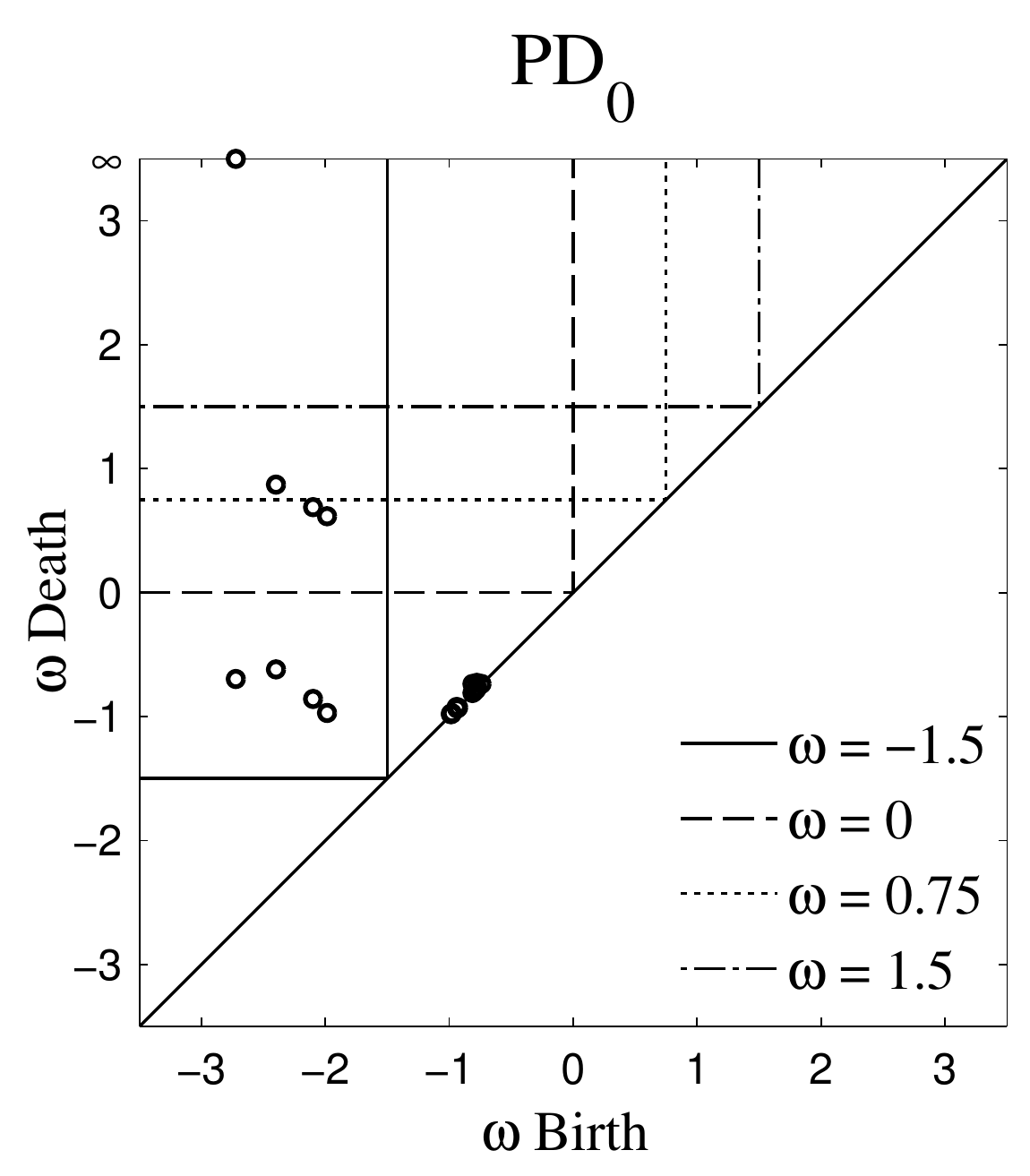}}\qquad
\subfigure[]{\includegraphics[width=1.6in]{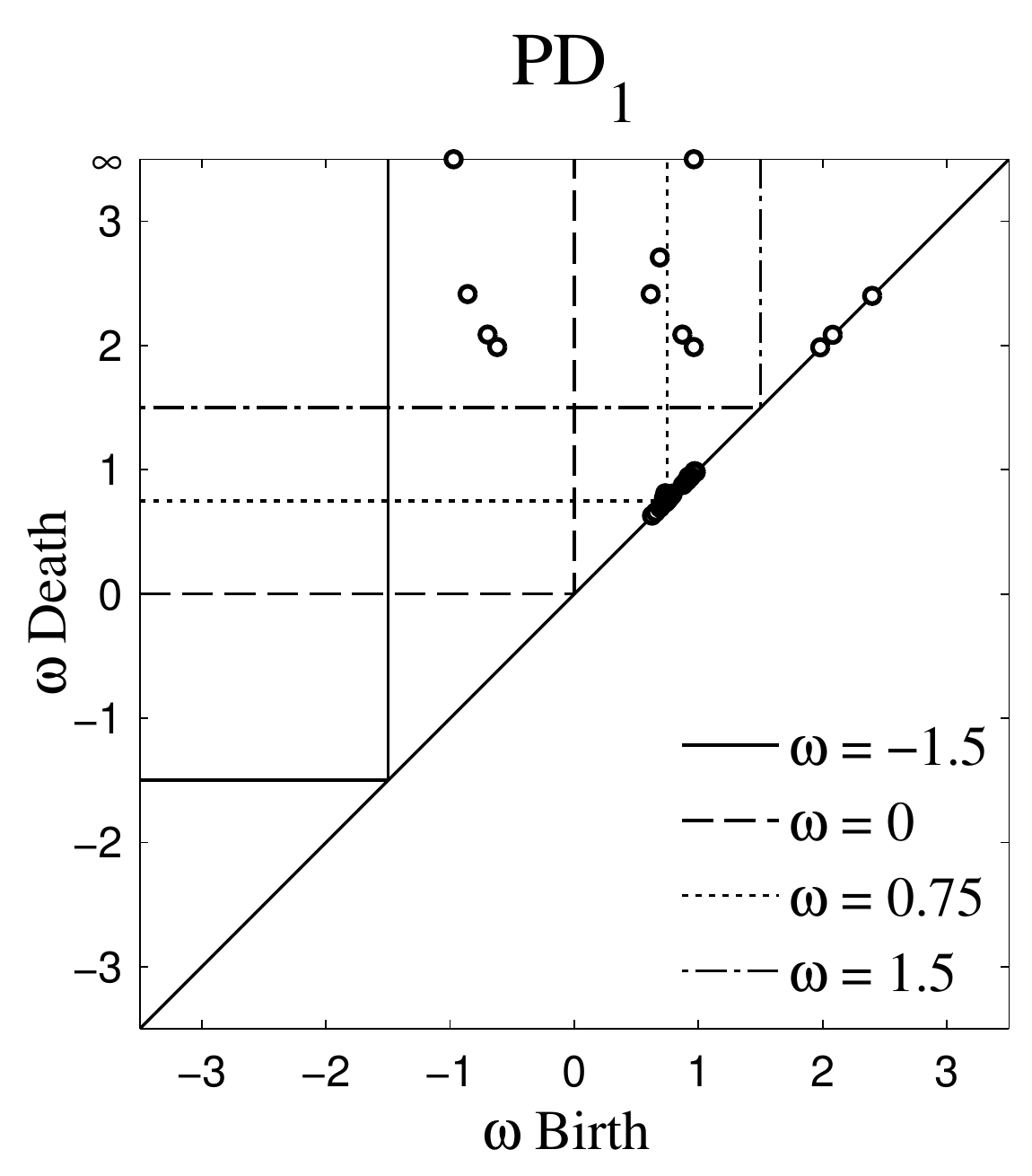}}\\
\caption{(a-d) Sub-level sets $\sC(\omega,\theta) = \setof{ x \in D : \omega(x) \leq \theta}$ of the vorticity field,  shown in Figure~\ref{fig:KLF}(a), for different  values of $\theta$, depicted in black.  (e)  $\pd_0(\omega)$ and (f) $\pd_1(\omega)$ persistence diagrams of the vorticity field indicate the values of $\theta$ at which the connected components  and loops appear and disappear (merge together). Every point $(\theta_b,\theta_d)$ in the $\pd_0(\omega)$ ($\pd_1(\omega)$) persistence diagram corresponds to a connected component (loop) that appears in $\sC(\omega, \theta_b)$ for the first time and is present in every set $\sC(\omega,\theta)$, for $b < \theta < d$, but   not  in   $\sC(\omega,\theta_d)$. A connected component disappears by merging with a previously existing component and  a loop disappears when it is filled in. Video 1 of the supplementary materials provides an animation.
\label{fig:PDs_kolmogorov}
}
\end{figure}

\section{Interpreting Persistence Diagrams}
\label{sec:PH}

The purpose of this section is to provide intuition and interpretation of the information that persistence diagrams present. As indicated in the previous section, we are interested in the diagrams $\pd_k(\omega)$, $k=0,1,2$, of the vorticity field for Kolmogorov flow, and the diagrams $\pd_k(T^*)$, $k=0,1,$ of the temperature field for Rayleigh-B\'{e}nard convection, shown in Figure~\ref{fig:KLF}.

\begin{figure}[t!]
\centering
\subfigure[]{\includegraphics[width=1.5in]{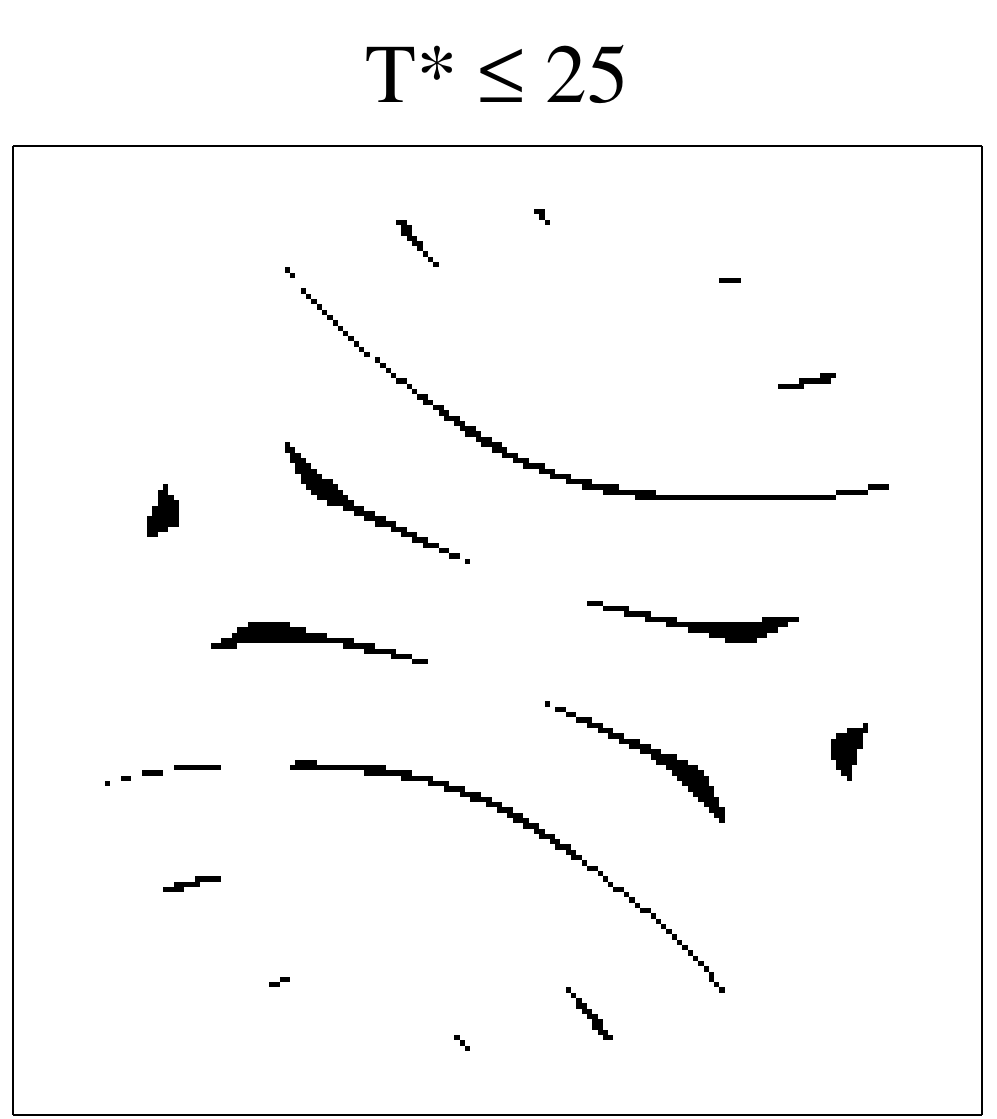}}\qquad
\subfigure[]{\includegraphics[width=1.5in]{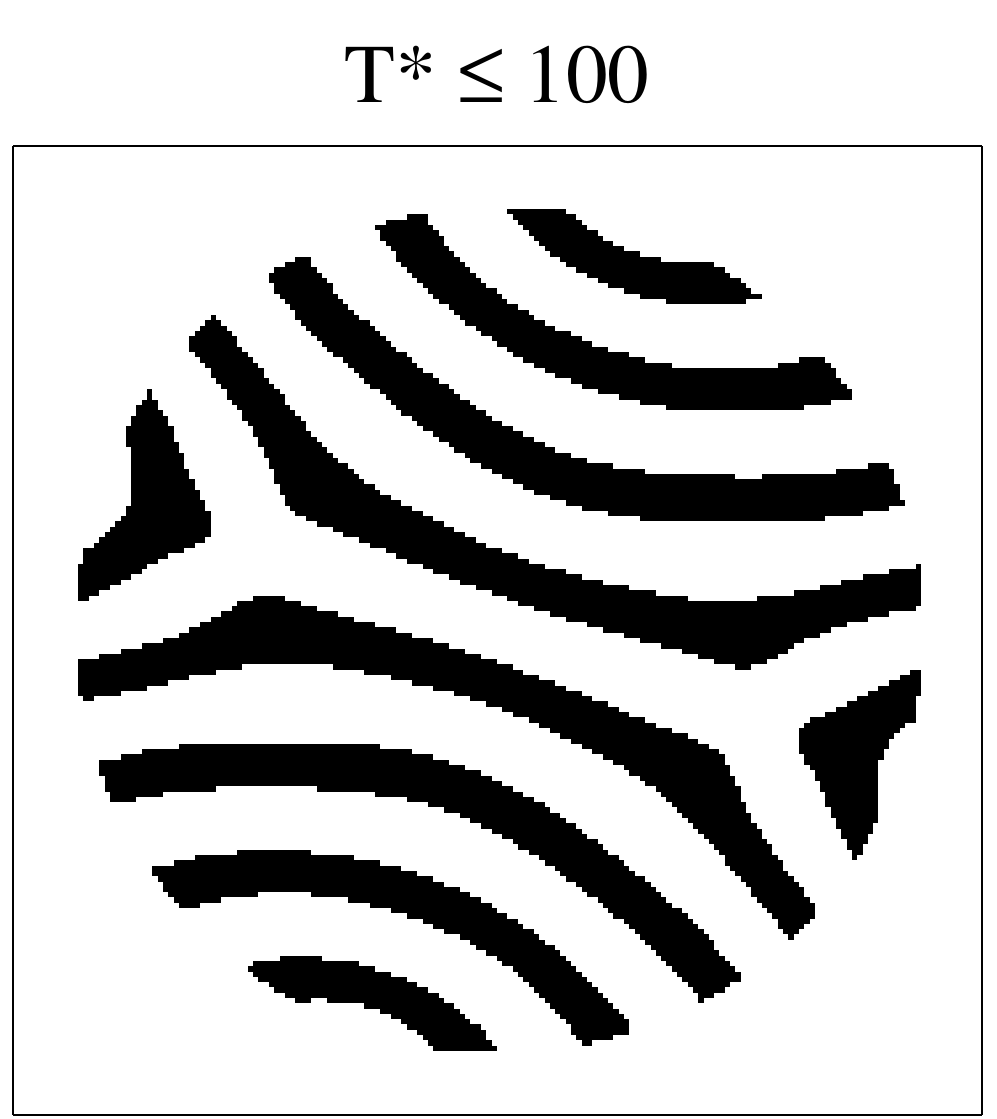}}\\
\subfigure[]{\includegraphics[width=1.5in]{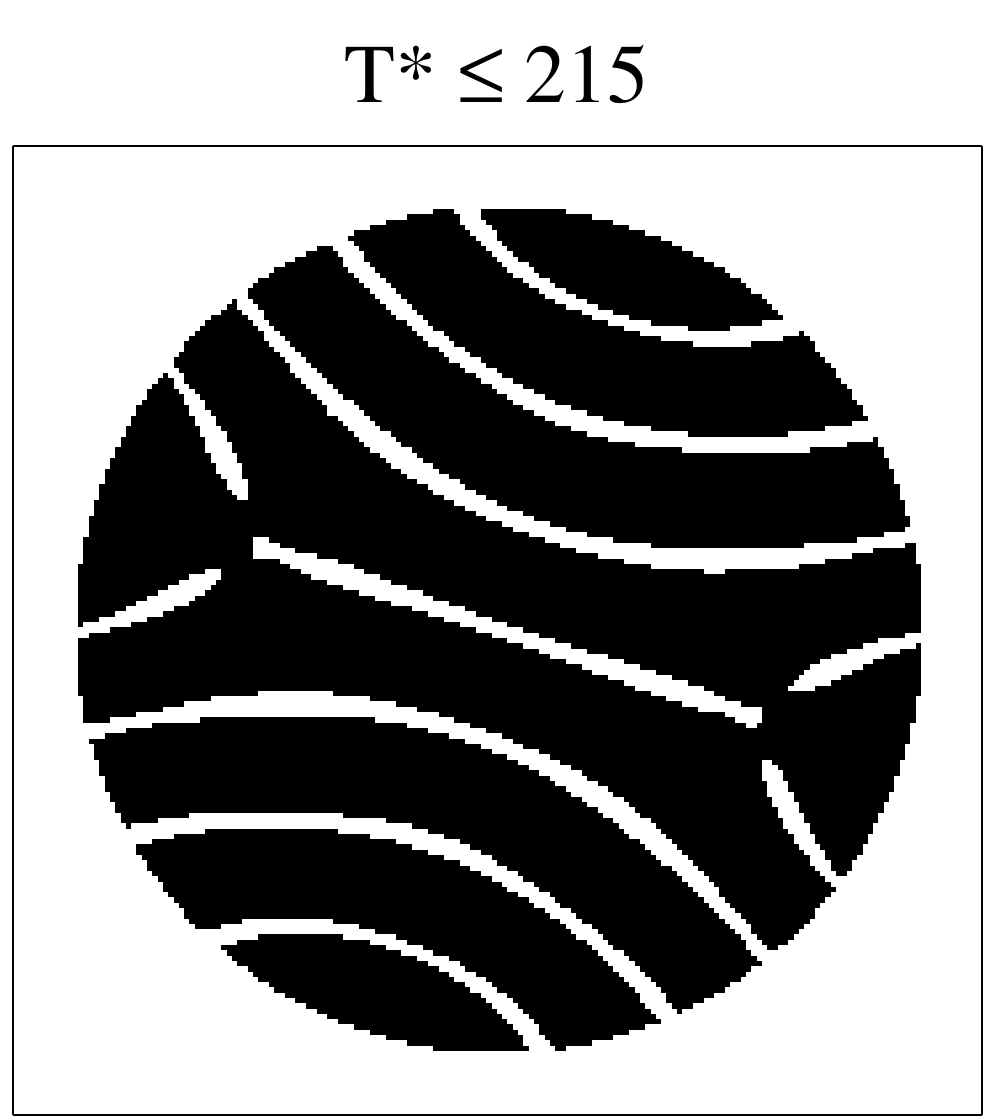}}\qquad
\subfigure[]{\includegraphics[width=1.5in]{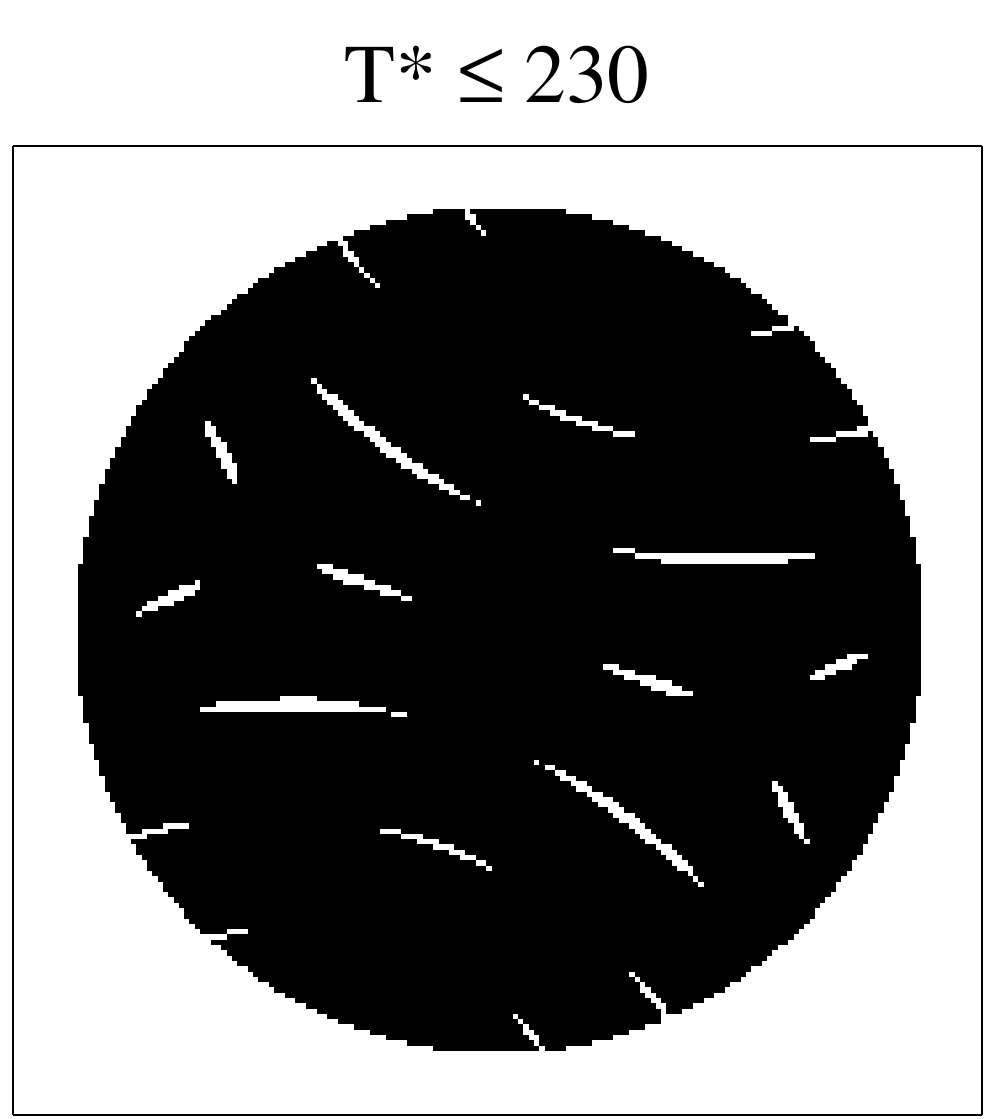}}\\
\subfigure[]{\includegraphics[width=1.6in]{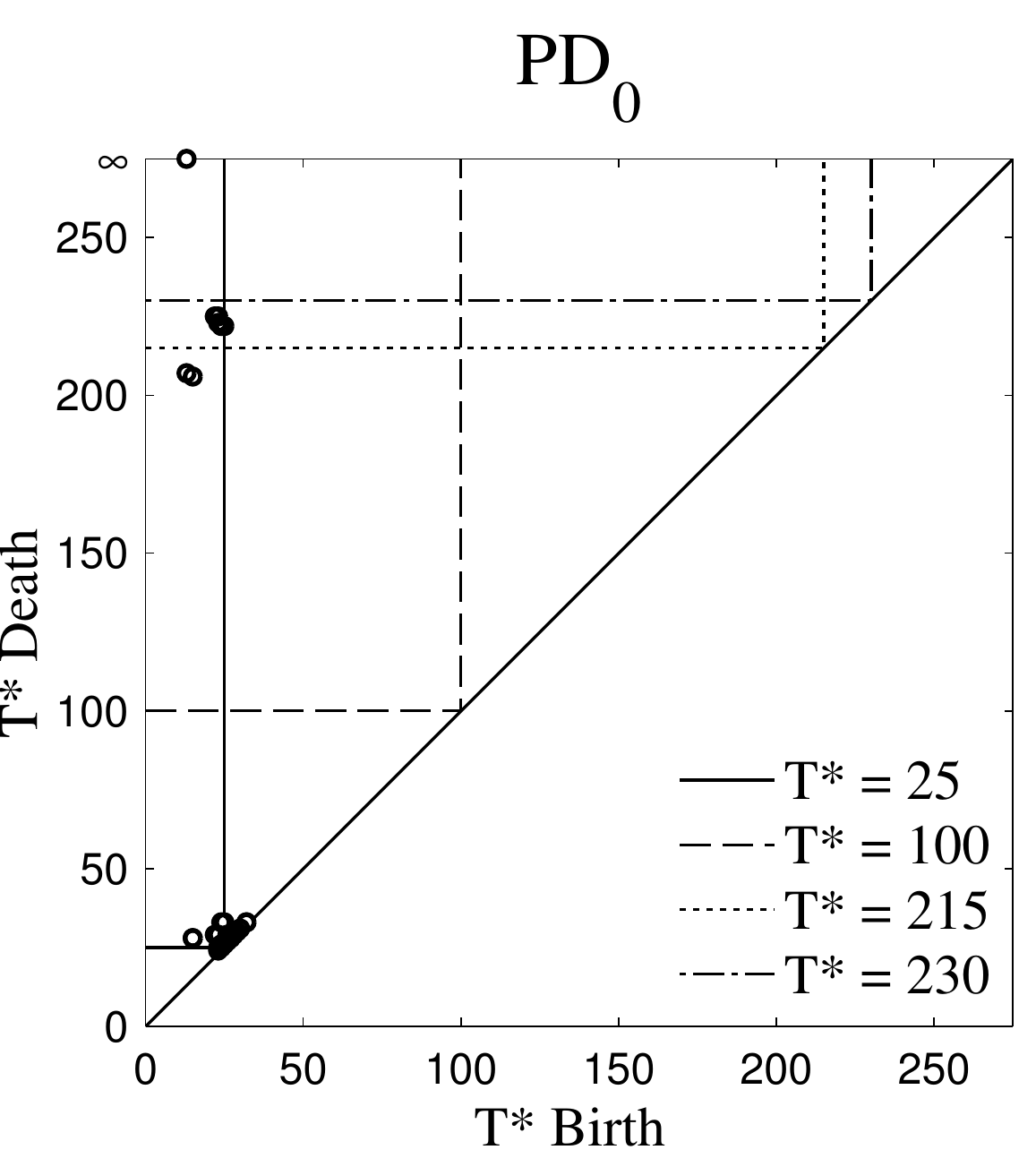}}\qquad
\subfigure[]{\includegraphics[width=1.6in]{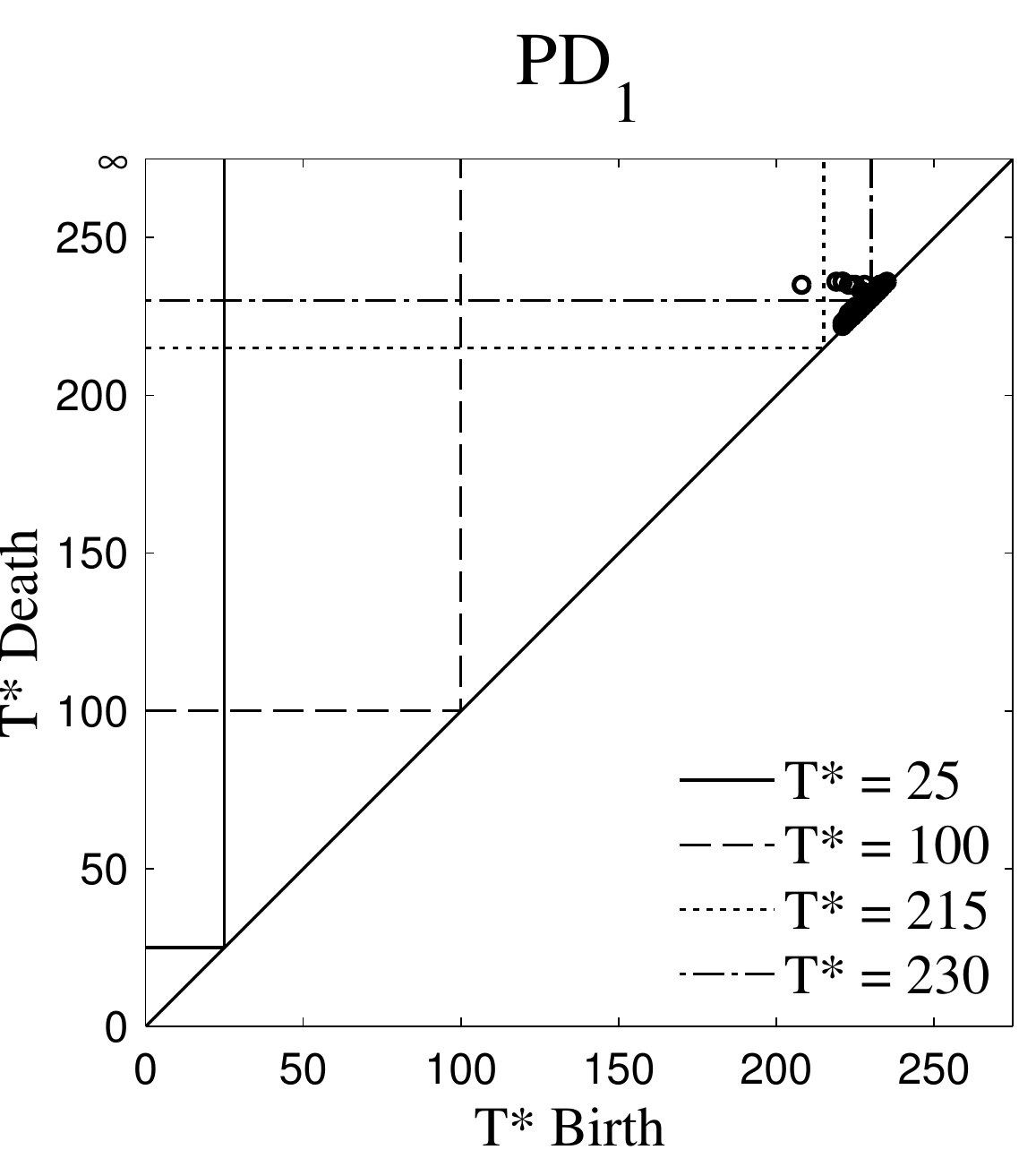}}\\
\caption{(a-d) Sub-level sets $\sC(T^*,\theta) = \setof{ x \in D : \T^*(x) \leq \theta}$ of the the renormalized 8-bit temperature field $T^*$,  shown in Figure~\ref{fig:KLF}(b), for different  values of $\theta$, depicted in black.  As before the persistence diagrams (e)  $\pd_0(\omega)$ and (f) $\pd_1(\omega)$  indicate the values of $\theta$ at which the connected components  and loops appear and disappear (merge together). Video 2 of the supplementary materials provides an animation.}
\label{fig:PDs_panam}
\end{figure}

We begin by discussing $\pd_0(\omega)$,  shown in Figure~\ref{fig:PDs_kolmogorov}(e), computed from a single time snapshot of the vorticity field for the Kolmogorov flow.  The minimum value of the vorticity field is $ -2.7206$, and therefore, $\sC(\omega,\theta)=\emptyset$ for all $\theta < -2.7206$. At $\theta = -2.7206$, two components appear, indicated by the two persistence points with birth value $\theta_b = -2.7206$. The death value of one of these two persistence points is $\theta_d = -0.697$, and so  the two components merge at this value, resulting in a single component. This explains the persistence point $(-2.7206, -0.697)$. The reason the other persistence point is denoted by $(-2.7206, \infty)$, with $\theta_d = \infty$, is because when features merge, a choice must be made about which topological feature (in this case, a connected component)  dies.  Having a consistent choice of basis over all values of $\theta$ requires that the homology generator associated with the geometric feature that has the larger birth value die first. If the birth values are the same, then it does not matter which topological feature with this birth value is chosen to be the one that persists. In particular, this implies that the generator associated with one of these two initial components can never die. 

Figure~\ref{fig:PDs_kolmogorov}(a) indicates the subset of $D$ corresponding to $\sC(\omega,-1.5)$. We remind the reader that the domain $D$ for Kolmogorov flow is a torus, since the left (top) and right (bottom) boundaries  are identified. Therefore,  $\sC(\omega,-1.5)$ consists  of eight distinct  connected components instead of nine. 

The existence of these eight connected components can also be extracted from $\pd_0(\omega)$, shown in Figure~\ref{fig:PDs_kolmogorov}(e). Observe that these connected components correspond to connected regions with birth value $\theta_b \leq -1.5$ and death value $\theta_d > -1.5$.   In Figure~\ref{fig:PDs_kolmogorov}(e), this corresponds to the eight points in the rectangular region  $R_{-1.5} := \setof{ (\theta_b,\theta_d)\in\R^2\mid \theta_b \leq -1.5 \text{ and } \theta_d > -1.5}$.

Figure~\ref{fig:PDs_kolmogorov}(b) indicates that $\sC(\omega,0)$ consists of four connected horizontal bands, which agrees with the number of persistence points in the rectangular region $R_{0}= \setof{ (\theta_b,\theta_d)\in\R^2\mid  \theta_b \leq 0 \text{ and } \theta_d > 0}$ of $\pd_0(\omega)$. Each stripe is created as two distinct components present in Figure~\ref{fig:PDs_kolmogorov}(a) grow and merge, causing one component to die each time.  The deaths of these components are captured by the points in the rectangle $R_{-1.5}$ which are not  in the rectangle $R_{0}$, since these are components that are born before $\theta = -1.5$ but die before $\theta = 0$.

Three  horizontal stripes merge together before  $ \theta = 0.75$, as indicated by two points inside the rectangle $R_{-1.5}$ that are not  in the rectangle $R_{0.75}$.  The two remaining connected components merge together soon thereafter, and for all greater threshold values, there is only one connected  component. 

To finish our analysis of $\pd_0(\omega)$, we  turn our attention to the persistence points close to the diagonal.  These have very short life spans, which suggests that these features may be numerical artifacts.  In our example, these points represent  the narrow horizontal bands formed in between two connected components before they merge into a single band (see video 1 available in the supplementary materials). These narrow bands are formed by  small oscillations of the vorticity field at the places where the field is almost constant. 

We now turn our attention to  the $\pd_1(\omega)$ persistence diagram, which characterizes loops in $\sC(\omega,\theta)$.   Appendix~\ref{sec:HomologyOfTorus} provides a detailed discussion  of independent loops  on a torus.  From $\pd_1(\omega)$, we see that the first loop appears at threshold $\theta = -0.963$. It corresponds to one  of the four horizontal bands shown in Figure~\ref{fig:PDs_kolmogorov}(b).  Each horizontal band  generates a single independent loop, corroborated by the existence of four persistence points in the rectangle $R_{0}$ of $\pd_1(\omega)$.

We note that the full torus has two loops captured by homology. This is expressed in $\pd_1(\omega)$ by the  two persistence points with $\theta_d = \infty$. Observe that the first loop that appears at $\theta = -0.963$ is equivalent to one of the toral loops, thus it cannot be killed by any other loop, and hence is captured by the persistence point $(-0.963,\infty)$. The other three loops present at $\theta = 0$  correspond to the same toral loop and thus must die. In fact, they do so by $\theta = 2.5$. Note that the birth values $\theta_b$ of these persistence points are close to the death values $\theta_d$ of the persistence points in $R_{-1.5}\setminus R_{0}$ of $\pd_0(\omega)$. This implies that shortly after the components merge, they form horizontal bands across the entire domain.

New loops are also created as the bands start merging. If two horizontal bands are connected by $n$ links, then the number of loops generated by this object (two bands plus the links)  is  $(1 + n)$. Thus, the first additional  loop appears when a second link is created (see Appendix~\ref{sec:HomologyOfTorus}). In our example, this happens near the threshold $ 0.75$. 

In Figure~\ref{fig:PDs_kolmogorov}(c), there are four distinct links between the two horizontal bands at the top of the figure.  The small punctures visible in  Figure~\ref{fig:PDs_kolmogorov}(c) are filled in almost immediately, and the four links merge into two distinct links.  The points in $\pd_1(\omega)$ that are close to the diagonal capture this behavior.  The other two links are present for a  wider range of thresholds, and the loop they generate is represented by one of the persistence points in $\pd_1(\omega)$ with birth coordinate slightly smaller  than  $0.75$.  The horizontal band at the top and the horizontal band at the bottom are linked in a similar fashion. This explains the presence of another point with  birth coordinate slightly smaller  than  $0.75$. 

At $\theta = 0.932$, a connection from the top to the bottom boundary is created.   This loop is homologically equivalent to the second of the two independent  loops of the torus, and hence is identified by the persistence point  $(0.932,\infty)$.  As  the threshold passes the value $1.988$, the punctures shown in Figure~\ref{fig:PDs_kolmogorov}(d) start disappearing and the corresponding  loops start dying.  Again, there are $1 + n$ independent loops for $n > 0$ punctures.    Since the maximum value of $\omega$ is $2.7092$, the sub-level set  is the whole torus for any threshold above this, i.e.\ $\sC(\omega,\theta)=D$ for all $\theta \geq 2.7092$. In this case, there are no more punctures, and the rectangle $R_{2.7092}$ contains only two persistence points.

Finally, we address the $\pd_2(\omega)$ persistence diagram, not shown for brevity.  This diagram contains a single persistence point at $(\theta_b,\infty)=(3.01,\infty)$. The birth coordinate, $\theta_b$, corresponds to the minimum value of $\theta$ for which $\sC(\omega,\theta) = D$, the whole torus. Since $\sC(\omega,\theta) = D$ for all $\theta \geq \theta_b$, this point never dies.

We now discuss the persistence diagrams for the temperature field $T^*$ shown in Figure~\ref{fig:PDs_panam} for Rayleigh-B\'{e}nard convection.  Again, beginning with $\pd_0(T^*)$, the points  with short life spans correspond the large number of small connected components that make up $\sC(T^*,25)$, as shown in  Figure~\ref{fig:PDs_panam}(a). Points with long life spans represent the well-defined connected components shown in  Figure~\ref{fig:PDs_panam}(b).  From the persistence diagram, we can see that these  components merge almost simultaneously at two threshold values, $\theta\approx 210$ and $\theta \approx 225$.  

Turning to $\pd_1(T^*)$, we note that the domain of the temperature field is a disk, so the independent loops correspond to punctures inside of the disk.  The diagram $\pd_1(T^*)$ indicates that there are no loops with long life spans, and the loops that do appear do so roughly at the same threshold values at which the dominant components merge. These features are due to the small fluctuations of the temperature field close to the critical values at which different rolls merge together. This  is  consistent with their short life spans.

\section{The Space of Persistence Diagrams}
\label{sec:SpaceOfPD}

As explained in the previous section, a persistence diagram codifies, in a reasonably compact form, considerable information about the geometry of a scalar function. As suggested by the examples, we use persistence diagrams to provide a reduced description of the state of the dynamical system of interest at any given point in time. Therefore, to analyze the dynamics, we need to be able to compare one collection of persistence diagrams $\pd$ (corresponding to a snapshot of the flow pattern at an instant of time) to another collection of diagrams $\pd'$ (from another flow snapshot).  There are a variety of metrics that can be imposed on persistence diagrams \cite{Chazal:2011:PCR:1998196.1998212, chazal:desilva:oudot, chazal:desilva:glisse:oudot, bubenik}.   The metrics used in this paper rely on pairing  the  points $p \in \pd_k$ in a one-to-one correspondence (bijection)  with the  points in $\pd'_k$.  According to the definition,  every persistence diagram contains an infinite number of copies of the diagonal. Hence, there are many different bijections $\gamma$ between $\pd_k$ and $\pd'_k$.   Roughly speaking, the distance between $\pd$ and $\pd'$ is defined  using   the bijections that  ``minimize the shift'' in the mapping of the points $p$ from $\pd_k$ to  $\gamma(p)$ in $\pd'_k$.  This notion is made more precise in the following definition.

\begin{defn}
\label{defn:bottleneck}
{\em
Let $\pd = \{\pd_k\}$ and $\pd' = \{\pd'_k\}$ be two collections of persistence diagrams.  The {\em bottleneck distance} between $\pd$ and $\pd'$ is defined to be 
\begin{equation}
d_B(\pd,\pd') = \max_{k}\inf_{\gamma\colon \pd_k \to \pd'_k} \sup_{p\in\pd_k} \| p -\gamma(p)\|_\infty , 
\label{eqn:dBDef}
\end{equation}
where $\|(a_0,b_0)-(a_1,b_1)\|_\infty := \max\setof{|a_0-a_1| , |b_0-b_1|}$ and $\gamma$  ranges over all bijections between persistence points. Similarly, the {\em degree-$p$ Wasserstein distance} is defined as
\begin{equation}
d_{W^p} (\pd,\pd')  =   \left[ \sum_{k} \inf_{\gamma\colon \pd_k \to \pd'_k} \sum_{p\in\pd_k} \| p -\gamma(p)\|^p_\infty \right]^{1/p}.
\label{eqn:dWDef}
\end{equation}
}
\end{defn}

Roughly speaking, a function $f\colon D\to\R$ is {\emph{tame} if, for every $\theta\in\R$,  the vector space $H_k(f^{-1}((-\infty,\theta]))$ is finite dimensional for every $k$, and there are only finitely-many thresholds at which the  vector spaces change (for a precise definition see \cite{edelsbrunner:harer}). For our purposes, it suffices to remark that if  $f$ is a piecewise-constant function on a finite complex, then $f$ is tame. In particular, the numerically-computed vorticity field $\omega$ and $8$-bit temperature field $T^*$ are tame functions.

For the remainder of this paper, we use $\per_k$ to denote the set of persistence diagrams corresponding to $H_k$ and $\per$ to denote the set of all persistence diagrams.  Let $T(D,\R)$ denote the set of tame functions $f\colon D\to \R$ equipped with the $L^\infty$ norm. A fundamental result \cite{edelsbrunner:harer} is that, using the Wasserstein or bottleneck metrics, $\pd\colon T(D,\R)\to\per$ is a Lipschitz-continuous function.   In particular, if $f,g\in T(D,\R)$, then
\begin{equation}
\label{eqn::Stability}
d_B(\pd(f),\pd(g)) \leq \sup_{x\in D} | f(x)- g(x)|.
\end{equation}

These results on Lipschitz continuity have two important implications for this work, both stemming from the fact that our analysis is based on numerical simulations. Assume for the moment that $f\colon D\to \R$ denotes the exact solution at a given time to either Kolmogorov flow or the Boussinesq equations. Ideally, we want to understand $\pd(f)$. Our computations of persistent homology are based on $\sC(\tilde{f},\theta_i)$, a cubical complex defined in terms of the numerically-reported values $\theta_i$, where $\tilde{f}$ represents the associated piecewise-constant function. If the numerical approximation $\tilde{f}$ satisfies $\sup_{x\in D} | f(x)- \tilde{f}(x)| \leq \varepsilon$, then by \eqref{eqn::Stability} we have a bound on the bottleneck distance between the actual persistence diagram $\pd(f)$ and the computed persistence diagram $\pd(\tilde{f})$, so that $d_B(\pd(f),\pd(\tilde{f})) \leq \varepsilon$.  Figure~\ref{fig:bottleneck} provides a schematic justification of this claim.

\begin{figure}[t]
\centering
\subfigure[]{\includegraphics[width=1.7in]{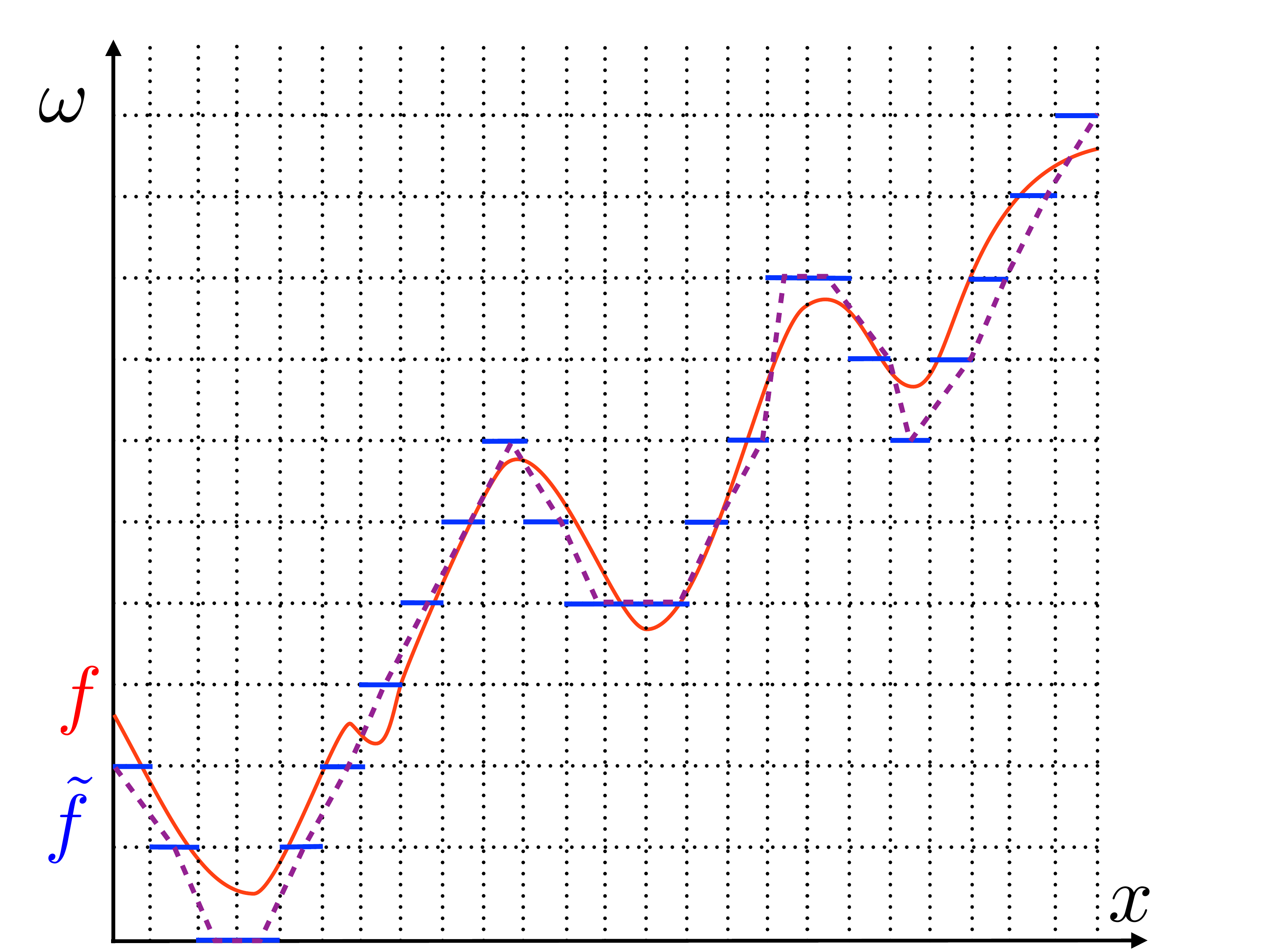}}
\subfigure[]{\includegraphics[width=1.7in]{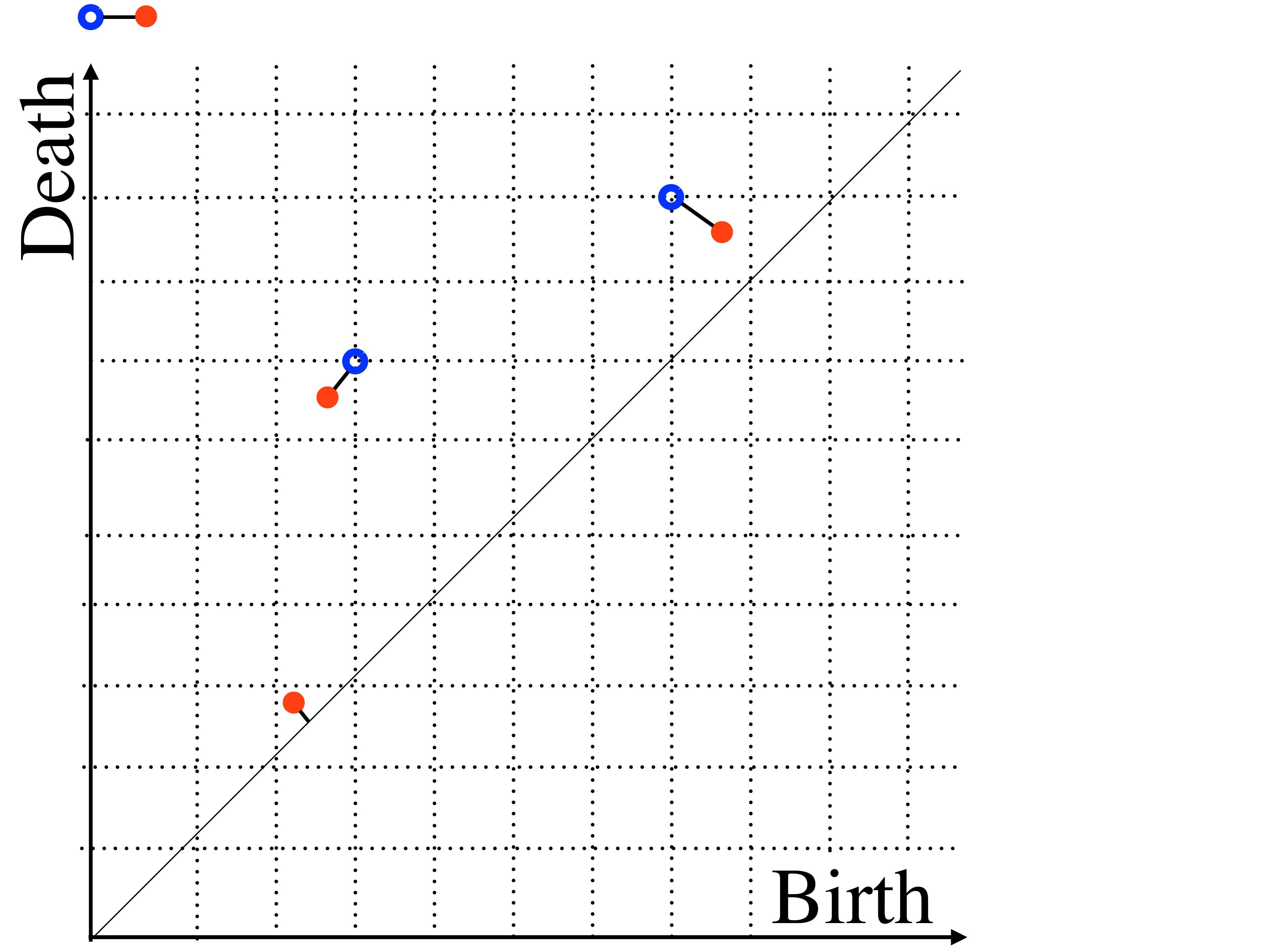}}
	
\caption{(a) A one-dimensional scalar field $f \colon D \to \R$ and its piecewise constant approximation $\tilde{f}$. The dashed line shows a tame approximation of $f$ whose persistence diagram is identical to $\pd(\tilde{f})$. (b) Points in the persistence diagram $\pd(f)$ are given by closed symbols and the points in   $\pd(\tilde{f})$ are represented by open symbols. Two points on the top have infinite death coordinate. Lines connecting the points represent matching of the persistence points  for which the  bottleneck distance $d_B(\pd(f),\pd(\tilde{f}))$ is realized. It follows from (\ref{eqn::Stability}) that $d_B(\pd(f),\pd(\tilde{f}))  \leq \sup_{x\in D} | f(x)- \tilde{f}(x)|$.  
\label{fig:bottleneck}}
\end{figure}

As indicated in the introduction, persistent homology is invariant under certain continuous deformations of the domain.  To be more precise, if $g\colon D\to D$ is a homeomorphism and $f\colon D\to \R$, then $\pd(f\circ g) = \pd(f)$.  Of particular relevance to this paper is a function $g$ which arises as a symmetric action on the domain. In this paper, we work with piecewise-constant numerical approximations of the actual functions of interest, and we cannot assume that this equality holds. However, if $f$ is given and $f' = f\circ g$, where $g$ is as above, and we have an $L^\infty$ bound $\varepsilon$ on the difference between the approximation and the true function, then by \eqref{eqn::Stability},
\begin{equation}
d_B(D(\tilde{f}),D(\tilde{f'})) \leq 2\varepsilon.
\end{equation}
In summary, under the assumption of bounded noise or errors from numerical simulations (or experimental data), we have explicit control of the errors of the distances in $\per$.

\section{Using Metrics in the Space of Persistence Diagrams}
\label{sec:Distances}

The goal of this section is twofold: one, to provide intuition about the information contained in the different metrics, and two, to suggest how viewing a time series in $\per$ can provide insight into the underlying dynamics.  

We begin by remarking that the bottleneck distance $d_B$ measures only the single largest difference between the persistence diagrams and ignores the rest. The Wasserstein distance $d_{W^p}$ includes all differences between the diagrams. Thus, it is always true that
\begin{equation}
d_B \leq d_{W^p}.
\end{equation}
The sensitivity of the Wasserstein metric to small differences (possibly due to noise) can be modulated by the choice of the value of $p$, i.e.\ if $p > q$, then one expects $d_{W^p}$  to be less sensitive to small changes than $d_{W^q}$. In this paper, we restrict ourselves to the bottleneck distance $d_B$ and the  Wasserstein distances $d_{W^p}$ for
$p=1,2$.

The most obvious use of these metrics is to identify or distinguish patterns. As an example, we consider patterns along an orbit from the Kolmogorov flow. As indicated in Section~\ref{Intro_Kolmogorov}, this particular trajectory arises from a periodic orbit  with a slow drift along an orbit of continuous symmetry. In particular, we consider the three time points indicated in Figure~\ref{fig:projections}(a): two that appear to differ by the continuous symmetry, and a third that lies on the `opposite' side of the periodic orbit. Plots of the associated vorticity fields at these points (see Figure~\ref{fig:3States}) agree with this characterization of the time points. We want to identify this information through the associated persistence diagrams $\pd^a$, $\pd^b$, and $\pd^c$, shown in Figure~\ref{fig:3StatesDiagram}.  Indeed, the plots of $\pd_k^a$ and $\pd_k^b$ are difficult to distinguish, but $\pd_k^c$ is clearly distinct. To quantify this difference, we make use of the distances between  the persistence diagrams  using $d_B$, $d_{W^2}$, and $d_{W^1}$. These values are recorded in Table~\ref{table:ThreeNetworksFullDistance}. Not surprisingly, the distances between $\pd^a$ and $\pd^b$ are much smaller than the distances between $\pd^a$ and $\pd^c$. We want to use these distances, as opposed to the detailed information in the persistence diagrams, to obtain rough information about how the pattern at Figure~\ref{fig:3States}(a) differs from the pattern at Figure~\ref{fig:3States}(c).

\begin{figure}[t]
\centering
\subfigure[]{\includegraphics[height=1.055in]{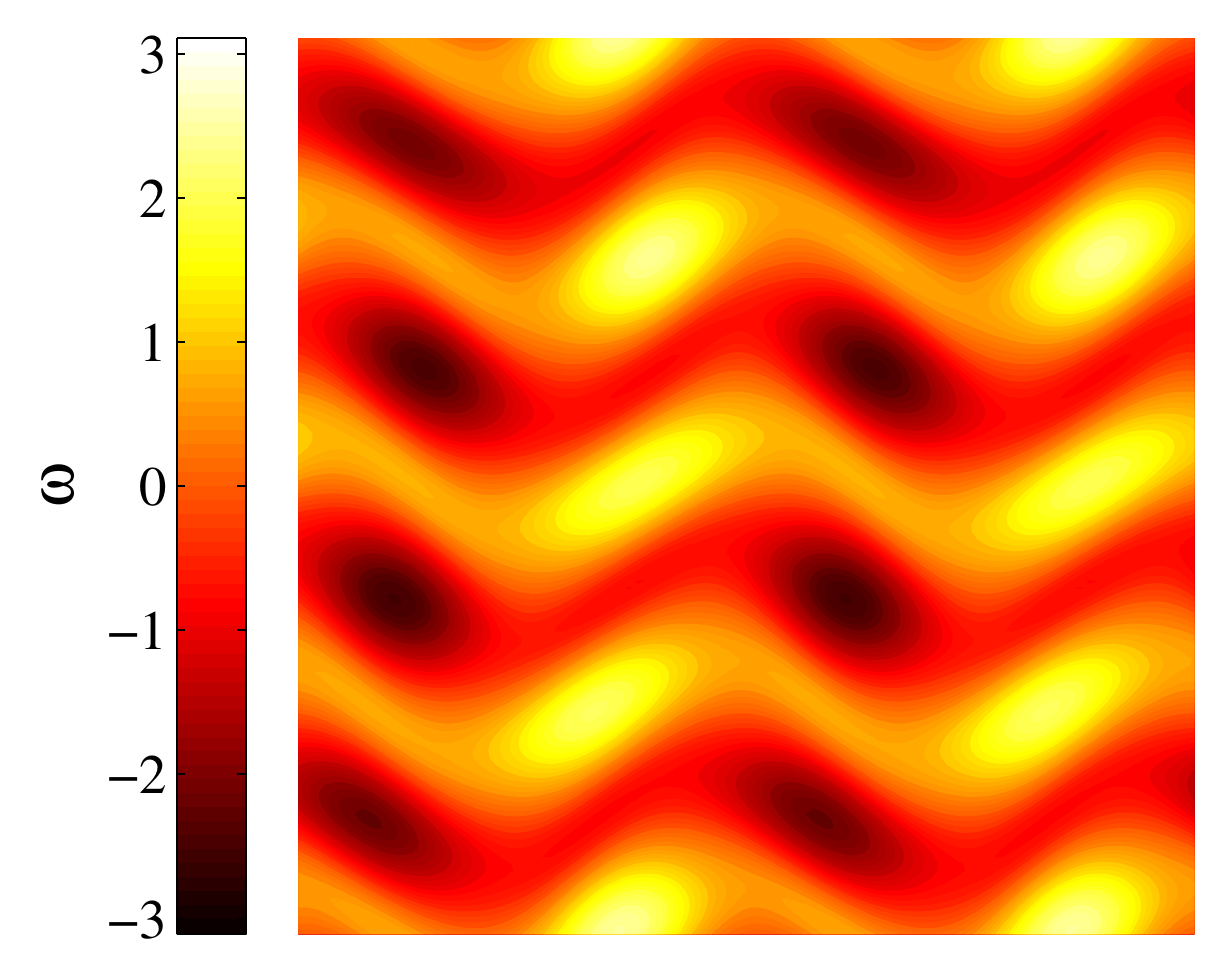}} 
\subfigure[]{\raisebox{0.025in}{\includegraphics[height=1.0in]{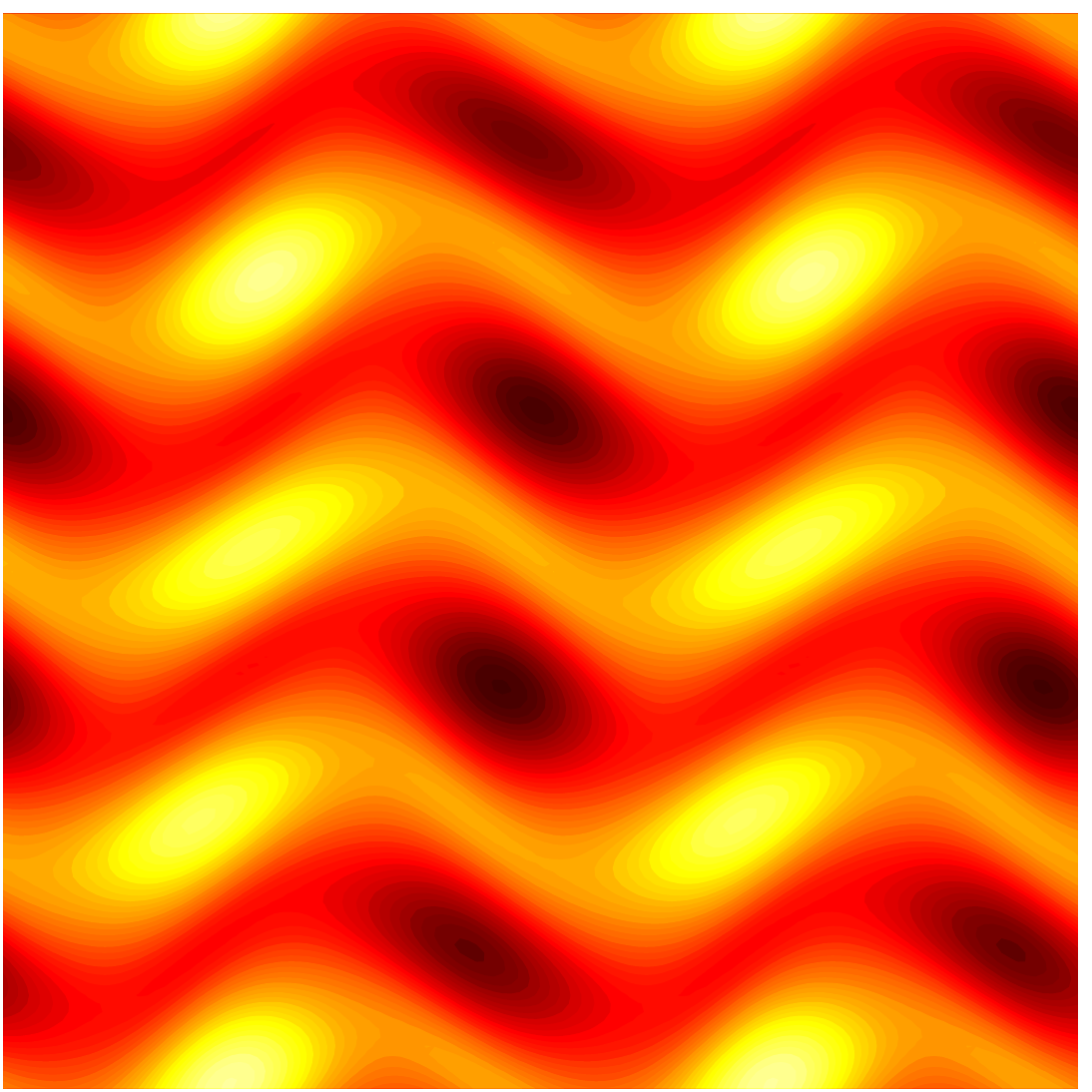}}}
\subfigure[]{\raisebox{0.025in}{\includegraphics[height=1.0in]{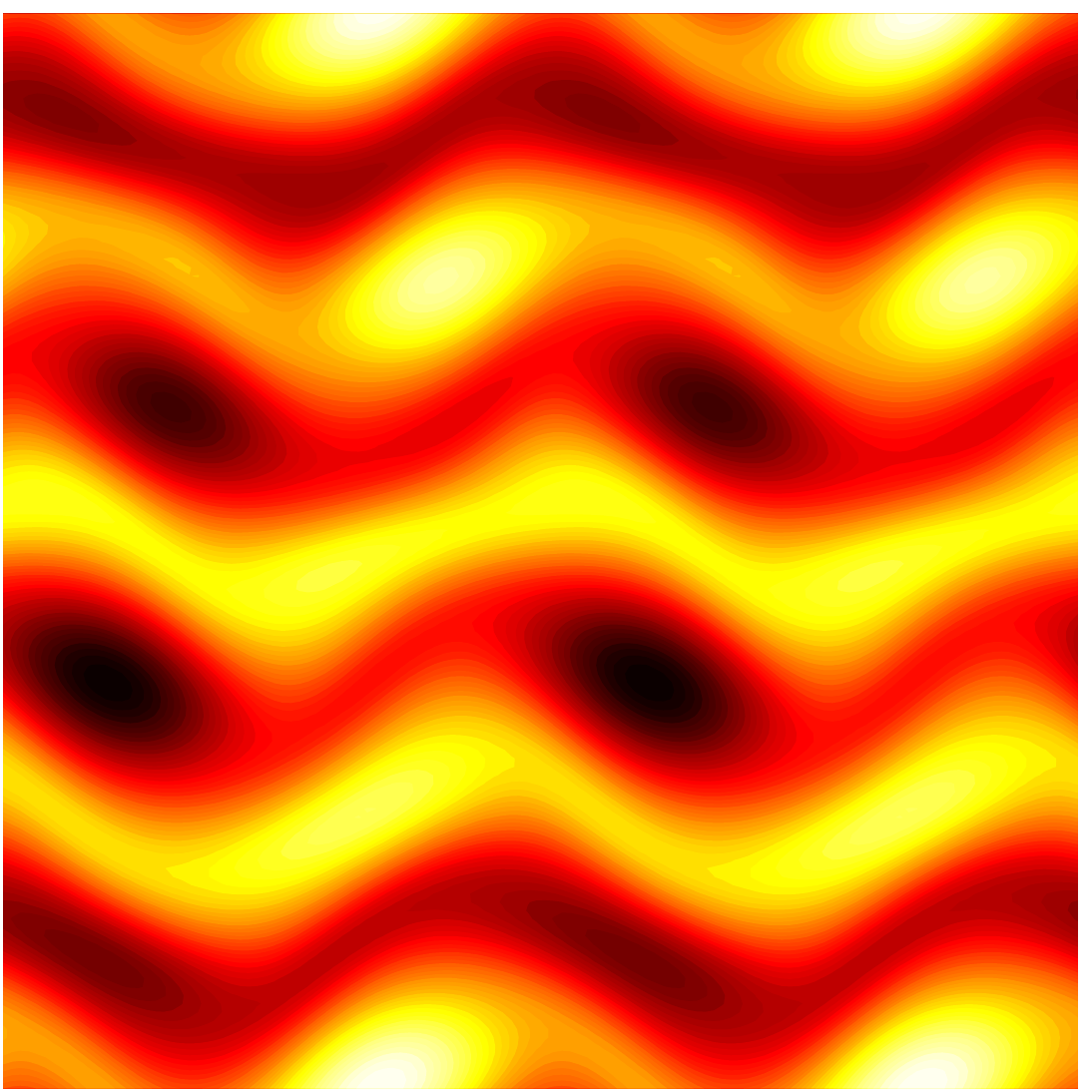}}}
\caption{Three snapshots of the vorticity fields $\omega$ from the stable relative periodic orbit of the Kolmogorov flow, found at $Re = 25.43$. The vorticity fields correspond to the (a) diamond, (b) square, and (c) circle in Figure \ref{fig:projections}(a).  The persistence diagrams for these three snapshots are generated and compared in Figure \ref{fig:3StatesDiagram}. Differences between the persistence diagrams are  qualitatively measured  by the distances shown in Table~\ref{table:ThreeNetworksFullDistance}. }
\label{fig:3States}
\end{figure}

\begin{figure}[]
\centering
\includegraphics[width=1.7in]{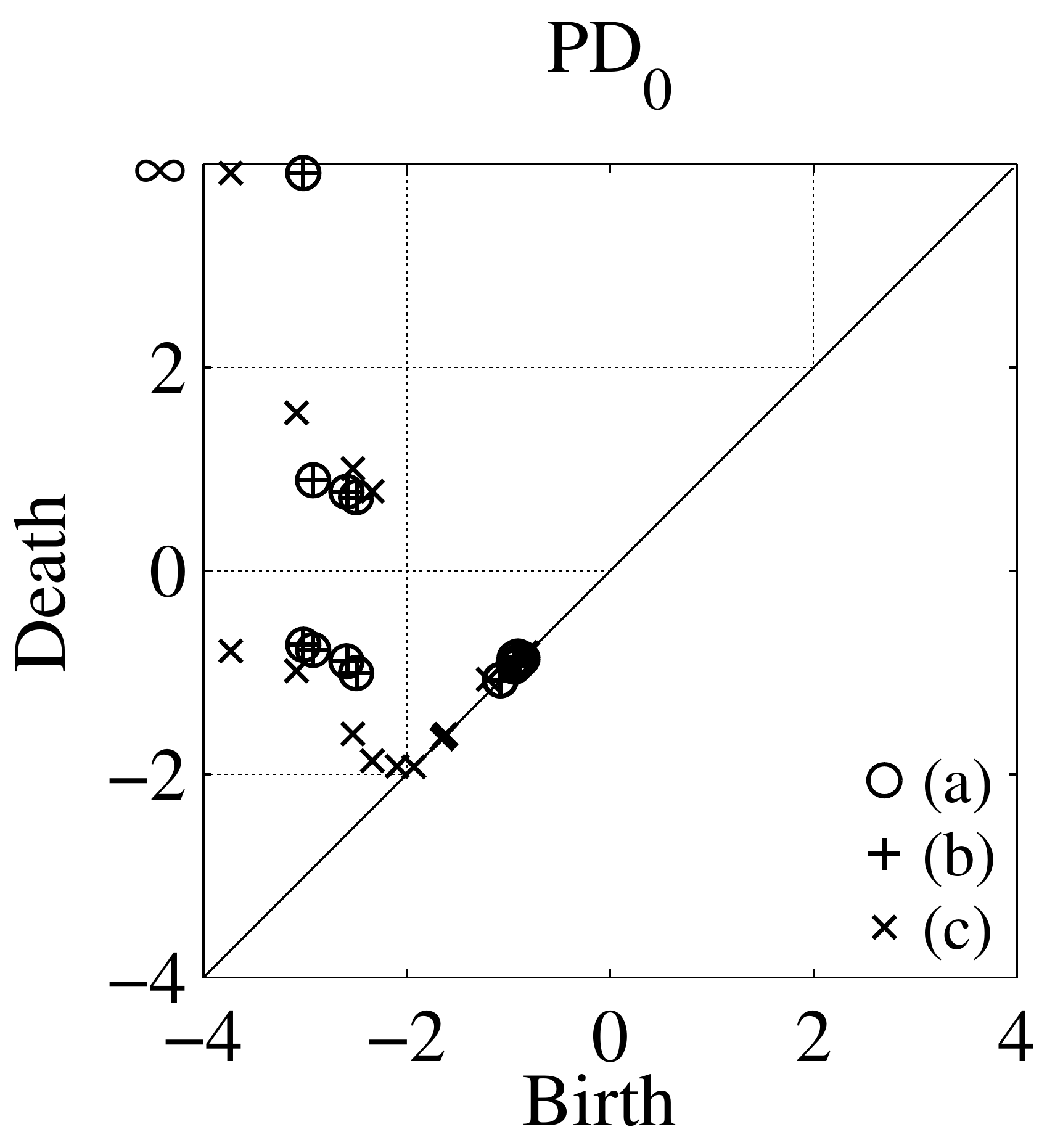} 
\includegraphics[width=1.7in]{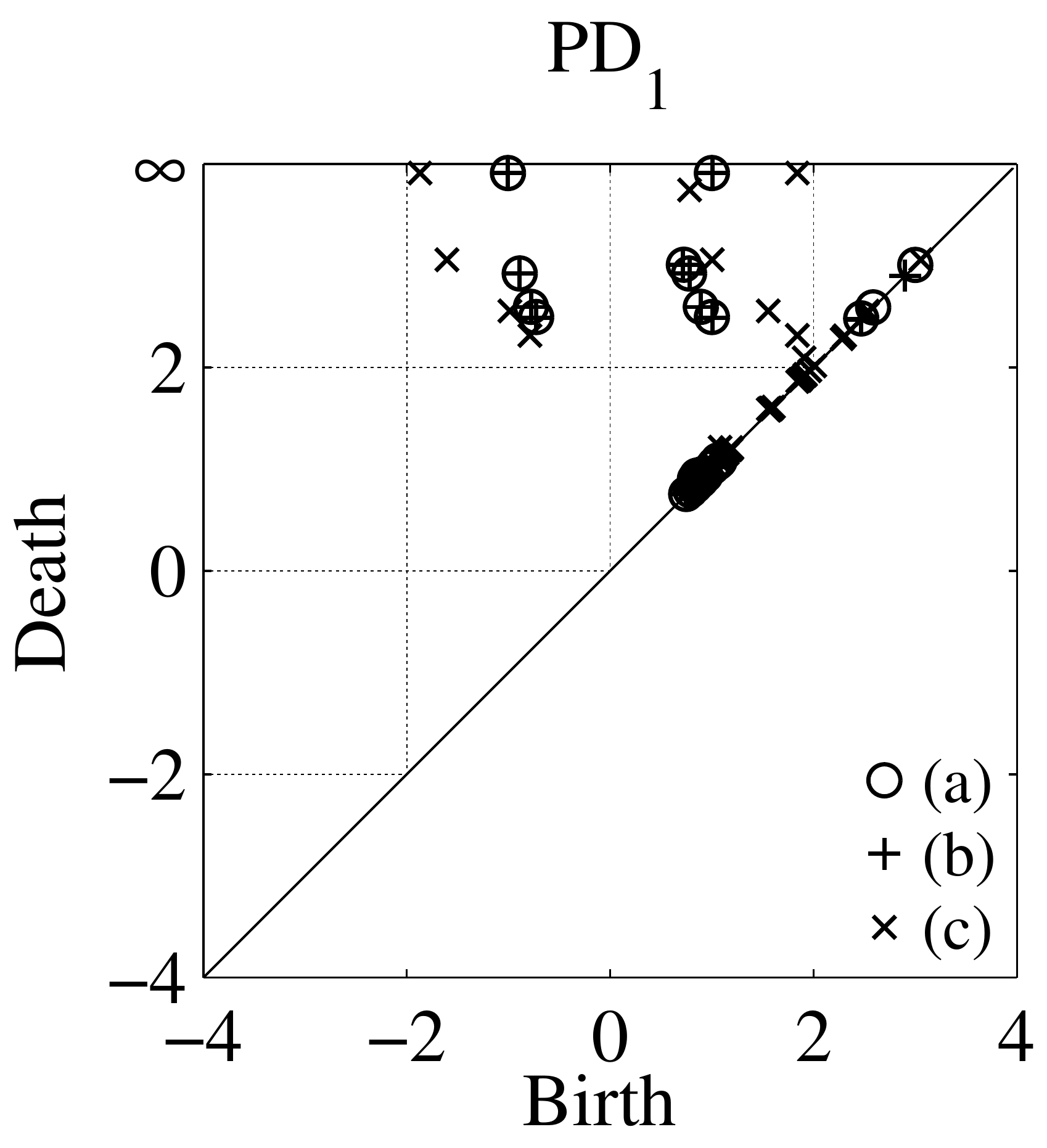} 
\caption{$\pd_0$ persistence diagrams $\pd^a, \pd^b$ and $\pd^c$ for the vorticity fields shown in  Figure~\ref{fig:3States}. 
The points in $\pd^a$ and $\pd^b$ are almost identical  because the corresponding vorticity fields are similar. The points in $\pd^c$ are more spread out and  do not shadow  the points in $\pd^a$ so well.  The same is true for the $\pd_1$ persistence diagrams which are not shown.  So  $d_\star(\pd^a,\pd^b) < d_\star(\pd^a,\pd^c)$, for $\star \in\setof{B,W^2,W^1}$, as indicated by Table~\ref{table:ThreeNetworksFullDistance}.  }
\label{fig:3StatesDiagram}
\end{figure}

\begin{table}
\centering 
\begin{tabular}{c  | c c c   }
\hline
 & $d_B$  & $d_{W^2}$ & $d_{W^1}$   \\
\hline
$(\pd^a,\pd^b)$ & $0.01$  & $0.049$ & $0.497$  \\ 
$(\pd^a,\pd^c)$ & $0.864$  & $2.648$  & $12.35$ \\ 
ratio $\frac{(\pd^a,\pd^c)}{(\pd^a,\pd^b)}$& $86.4$ & $54.05$ & $24.85$ \\
\end{tabular}
\caption{Distances between selected persistence diagrams (rounded to $3$ decimal places) shown in Figure~\ref{fig:3StatesDiagram}, corresponding to  the vorticity fields  given by  Figure~\ref{fig:3States}.}
\label{table:ThreeNetworksFullDistance}
\end{table}

The patterns shown in Figure~\ref{fig:3States}(a)-(b) differ by a symmetric transformation, and so $d_B(\pd^a,\pd^b)$ can be interpreted as  a lower bound on the numerical errors. Now observe that either $\pd^a$ or $\pd^c$ must have a persistence point with life span greater than $d_B(\pd^a,\pd^c)$ (otherwise pairing the persistence points with the diagonal will produce a smaller $d_B$ distance). Since the ratio $d_B(\pd^a,\pd^c)/d_B(\pd^a,\pd^b)$ is $86.4$, we know that there is a significant distinction between differences that should be ascribed to error and differences that can be ascribed to significant geometric features.

For some applications, there might be  only two different scales at which the geometric features evolve: one scale corresponding to the signal, and the other  representing the noise fluctuations.  If that is  the case, then  there are only  two types of changes. If we suppose that the large changes are comparable to $d_B(\pd^a,\pd^c)$ and the noise fluctuations are of the order    $d_B(\pd^a,\pd^b)$, then  we  can approximate the distances $d_{W^1}(\pd^a,\pd^c)$  and  $d_{W^2}(\pd^a,\pd^c)$ as follows:
\begin{equation}
d_{W^1}(\pd^a,\pd^c)  \approx n \cdot d_{B}(\pd^a,\pd^c) + k\cdot d_B(\pd^a,\pd^b),
\label{Eqn::dW1}
\end{equation}
\begin{equation}
\begin{split}
& d_{W^2}(\pd^a,\pd^c) \approx \\ & \text{\;\;\;\;\;\;\;\;\;\;\;\;\;} \sqrt{n (d_{B}(\pd^a,\pd^c) )^2 + k (d_{B}(\pd^a,\pd^b) )^2},
\end{split}
\end{equation}
where $n$ is the number of features that change significantly and $k$ is the number of features that change very little.
We recall that $d_{B}(\pd^a,\pd^b) \ll d_{B}(\pd^a,\pd^c)$ (Table~\ref{table:ThreeNetworksFullDistance}). Hence, the significant contributions to the $d_{W^2}$ metric are of the order of $d_B(\pd^a,\pd^c)$, and the small changes do not significantly contribute  to $d_{W^2}$. This leads to the following approximation:
\begin{equation}
d_{W^2}(\pd^a,\pd^c) \approx \sqrt{n (d_{B}(\pd^a,\pd^c) )^2}.
\label{Eqn::dW2}
\end{equation}
By solving (\ref{Eqn::dW1}) and (\ref{Eqn::dW2}), we obtain $n=9$ and $k = 383$, and so the number of features that change significantly is bounded from below by $9$. Before we  discuss the values of  $n$ and $k$, let us repeat the same computation for the Rayleigh-B\'enard convection. 

\begin{figure}[]
\centering
\subfigure[]{\includegraphics[width=1.7in]{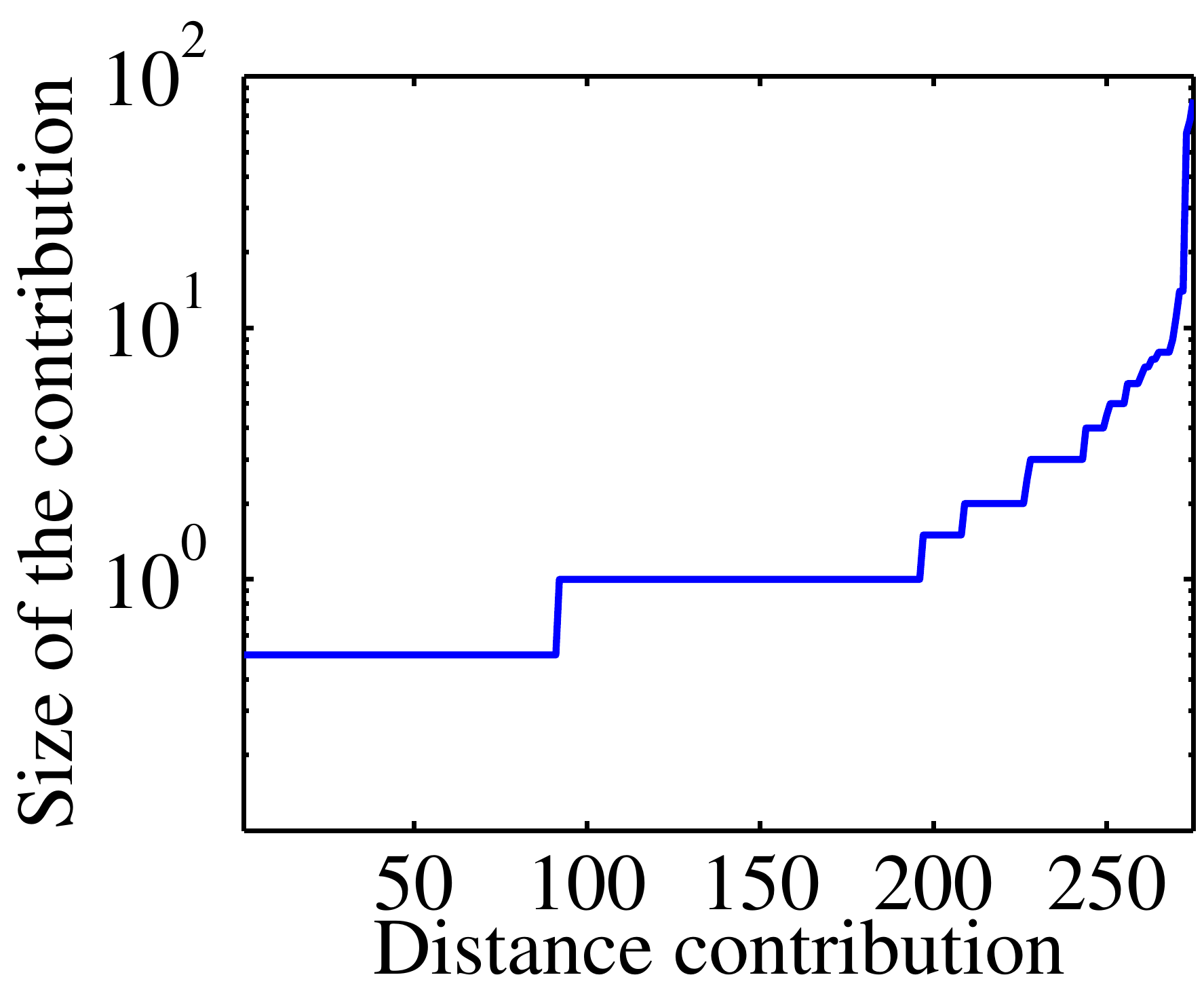}}
\subfigure[]{\includegraphics[width=1.7in]{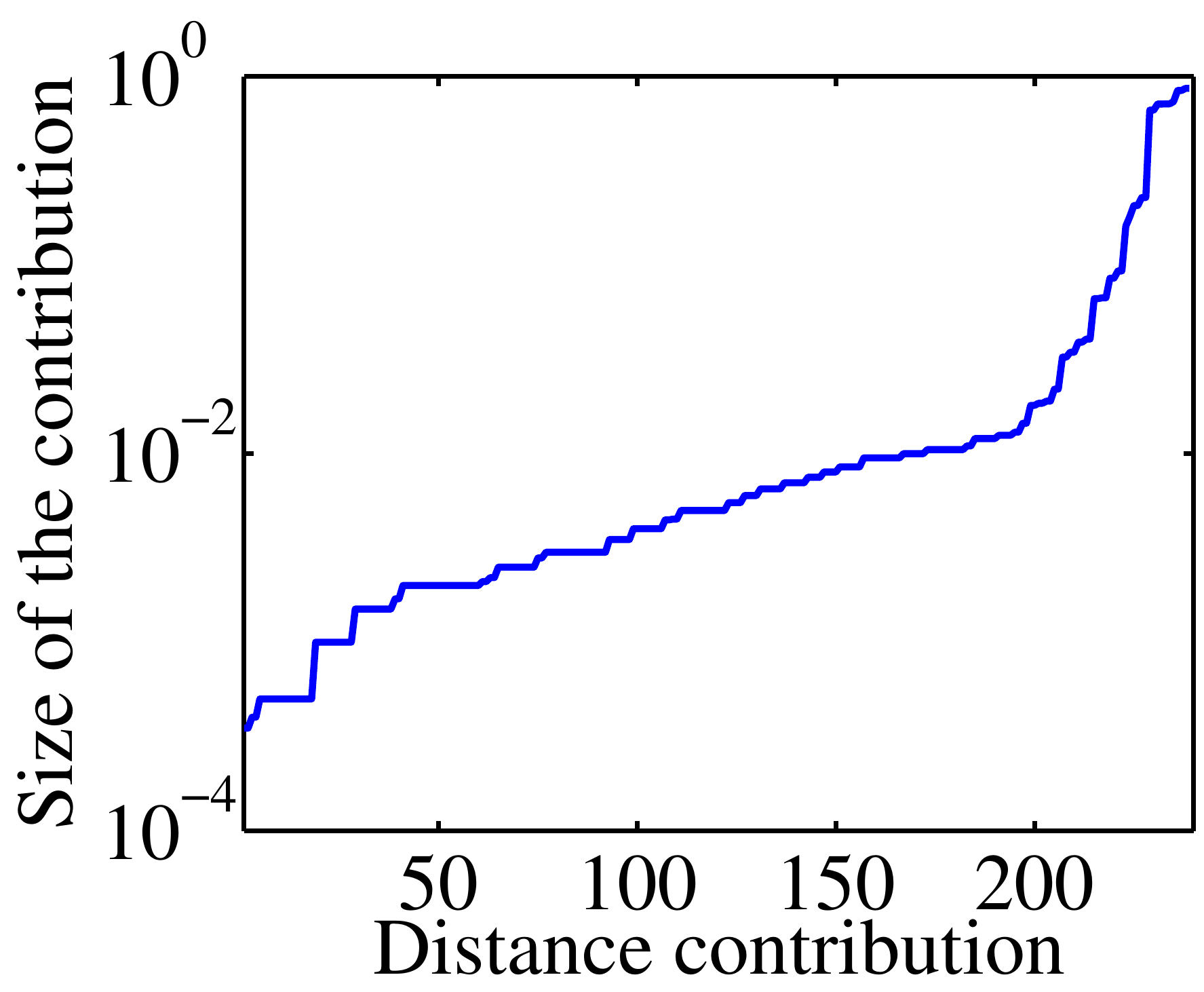}}

\caption{(a) Contributions to the $d_{W^1}(\pd^A.\pd^C)$  distance for Rayleigh-B\'enard convection.  (b) Contributions to the $d_{W^1}(\pd^a.\pd^c)$  distance for the  Kolmogorov flow.}
\label{fig:dW1_Contributions}
\end{figure}

We recall that the  $8$-bit temperature field $T^*$ is an integer-valued function with values between $0$ and $255$.  For integer-valued functions, the smallest nonzero $d_B$ distance between distinct frames is $e = 0.5$. We use this number as the lower bound on the numerical errors. The snapshots $A$ and $C$, not shown for brevity, are from a single orbit, and they realize the maximal distance between two snapshots (exact distances are given  by Table~\ref{table:ThreeNetworksRBCFullDistance}). Solving  Equations (\ref{Eqn::dW1}) and (\ref{Eqn::dW2}) yields $n = 35$ and $k = -4400$. These numbers obviously do not make sense. Therefore, the changes cannot be divided into two distinct groups, as we assumed above. Figure~\ref{fig:dW1_Contributions}(a) shows that there is  only one change on the order of $d_B(\pd^a,\pd^c)$. All the other changes are at least an order of magnitude smaller. More precisely, there are $77$ changes that are between one to two orders of magnitude smaller than $d_B(\pd^a,\pd^c)$. Moreover, these changes are at least an order of magnitude larger than our error estimate $e = 0.5$, and there are also $105$ changes of the size $2e$. Due to the significant number of contributions at all orders of magnitude between the noise estimate and $d_B(\pd^a,\pd^c)$, we cannot assume that there are predominantly two distinct types of changes: one corresponding to the noise and the other to the signal. Thus, the approximation  (\ref{Eqn::dW2}) of $d_{W^2}$ is not valid in this setting because a large part of $d_{W^2}(\pd^a,\pd^c)$ comes from contributions at the intermediate scales between $e$ and $d_B(\pd^a,\pd^b)$. 

\begin{figure*}[t]
\centering
\subfigure[]{\includegraphics[width=3.5in]{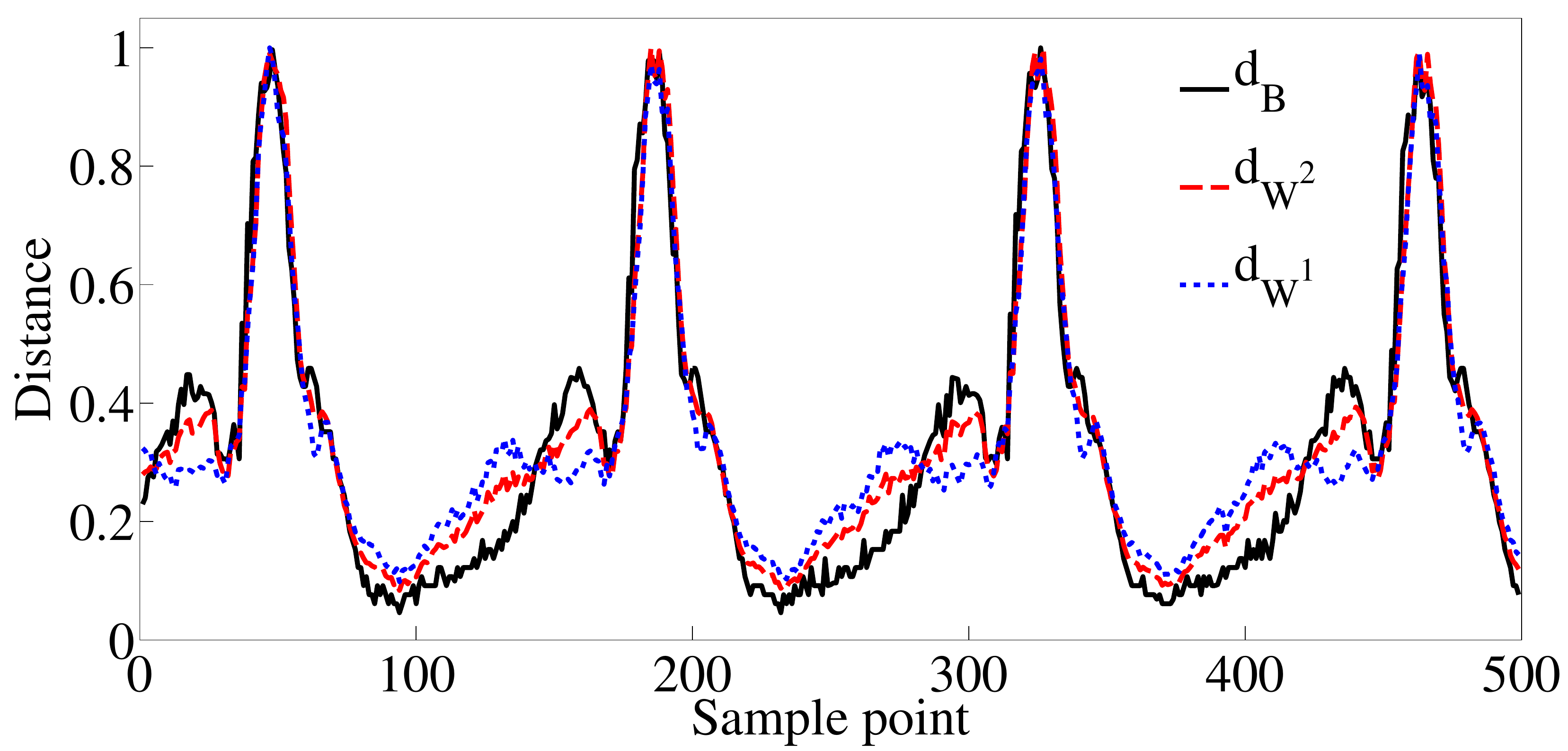}} 
\subfigure[]{\includegraphics[width=3.5in]{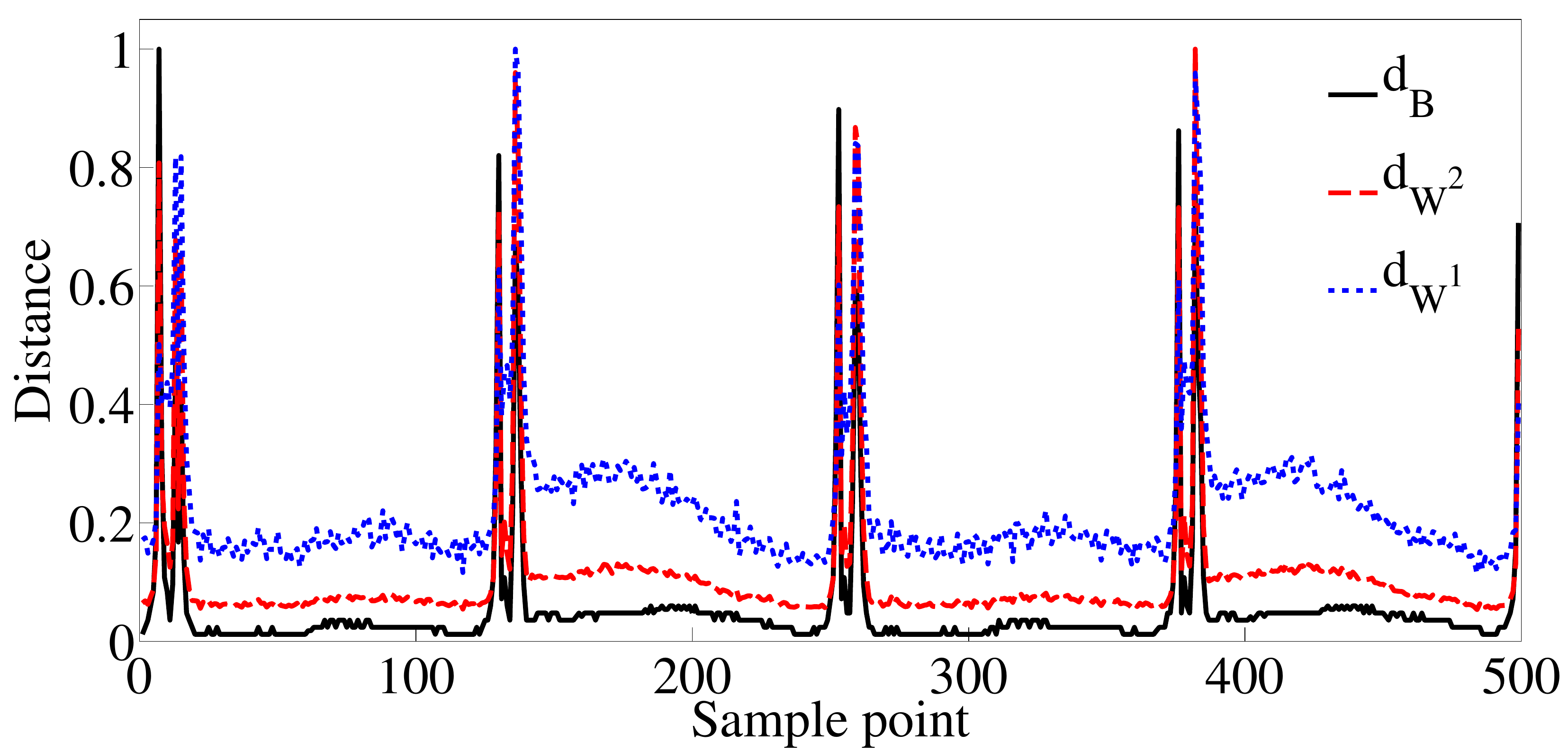} }
\caption{(a) Distances $d_\star$ between the consecutive sample points along the stable relative periodic orbit  of the Kolmogorov flow  are shown for just over three periods. Distances are normalized by their maximum value which is $0.0654$ for $d_B$,  $0.2266$ for $d_{W^2}$,  and $1.9143$ for $d_{W^1}$. Distance between the consecutive sample points  can be interpreted as the speed at which the orbit moves in the space of  persistence diagrams. Speed is not uniform along the orbit. Instead, there are parts of the orbit where the dynamics are slow, separated by relatively fast evolution. (b)  Distances $d_\star$ between the consecutive sample points along an almost-periodic orbit from Rayleigh-B\'enard convection are shown for approximately 2 periods. Distances are again normalized by their maximum value which is $83.5$ for $d_B$, $113.66$ for $d_{W^2}$, and $383$ for $d_{W^1}$. }
\label{fig:ConsecDist}
\end{figure*}

\begin{table}
\centering 
\begin{tabular}{c  | c c c   }
\hline
 & $d_B$  & $d_{W^2}$ & $d_{W^1}$   \\
\hline
$(\pd^A,\pd^C)$ & $81.5$  & $480.9$  & $650.5$ \\ 
\end{tabular}
\caption{Distances between selected persistence diagrams (rounded to $3$ decimal places)  corresponding to  two different $8$-bit temperature fields  obtained from a single orbit of  Rayleigh-B\'enard convection.}
\label{table:ThreeNetworksRBCFullDistance}
\end{table}

We now return to the values of $n$ and $k$  for the Kolmogorov flow. 
Figure~\ref{fig:dW1_Contributions}(b) shows that there are approximately $200$ changes of order smaller or equal to $d_B(\pd^a,\pd^b)$, and $11$ dominant changes of order $d_B(\pd^a,\pd^c)$. Finally, we can identify  approximately $28$ changes occurring on intermediate scales, and their sizes are at least an order of magnitude smaller than $d_B(\pd^a,\pd^c)$. In fact, most of them are not much larger than $d_B(\pd^a,\pd^b)$. Hence, the changes can be roughly divided into two classes of different order. The fact that the division is not absolutely sharp leads to $n = 9$, which is smaller than the actual number of  dominant changes on the order of $d_B(\pd^a,\pd^c)$.

We now turn to the question of understanding dynamics from the time series in $\per$. Let  $f_i$ denote the scalar field of the system at time $t_i$.  If $\Delta t = t_{i+1}-t_i$ is small and the evolution of the system is continuous, then  because $d_\star$ (for  $\star \in \setof{B,W^2,W^1}$ ) is a metric,
\begin{equation}
s_\star(t_i) = \frac{d_\star(\pd(f_i),  \pd(f_{i+1}))  }{\Delta t}
\end{equation}
can be interpreted as an average speed  in the space of  persistence diagrams over the time interval $[t_i, t_{i+1}]$.  The value of $s_\star$ depends on the choice of  metric.  For example,  $s_{d_B}$ is the rate at which the largest change between the geometric features of the scalar fields occurs. The speeds measured by  $d_{W^p}$, $p=1,2$, keep track of the rate of change between all geometric features, though to some extent, $d_{W^2}$ suppresses the effect of noise. 

Figure~\ref{fig:ConsecDist}(a) shows distances $d_*$ between consecutive sample points, normalized by the maximum distance, of samples taken along approximately three periods of the stable relative periodic orbit  of the Kolmogorov flow described in Section~\ref{Intro_Kolmogorov}.  Normalizing $s_*$ by the maximum speed along the orbit furnishes the same curves. Each of the graphs of $s_\star$ indicate that speed is not uniform along the orbit; there are parts of the orbit where the geometry is changing slowly, separated by intervals of relatively fast evolution.  The evolution is extremely slow around the states $100, 240,$ and $380$, where the speeds $s_\star$ are below the noise (fluctuation) levels  given by the first row of  Table~\ref{table:ThreeNetworksFullDistance}. This suggests that the orbit may be passing close to a fixed point. 

\begin{figure}[t]
\centering
\subfigure[]{\includegraphics[width=1.7in]{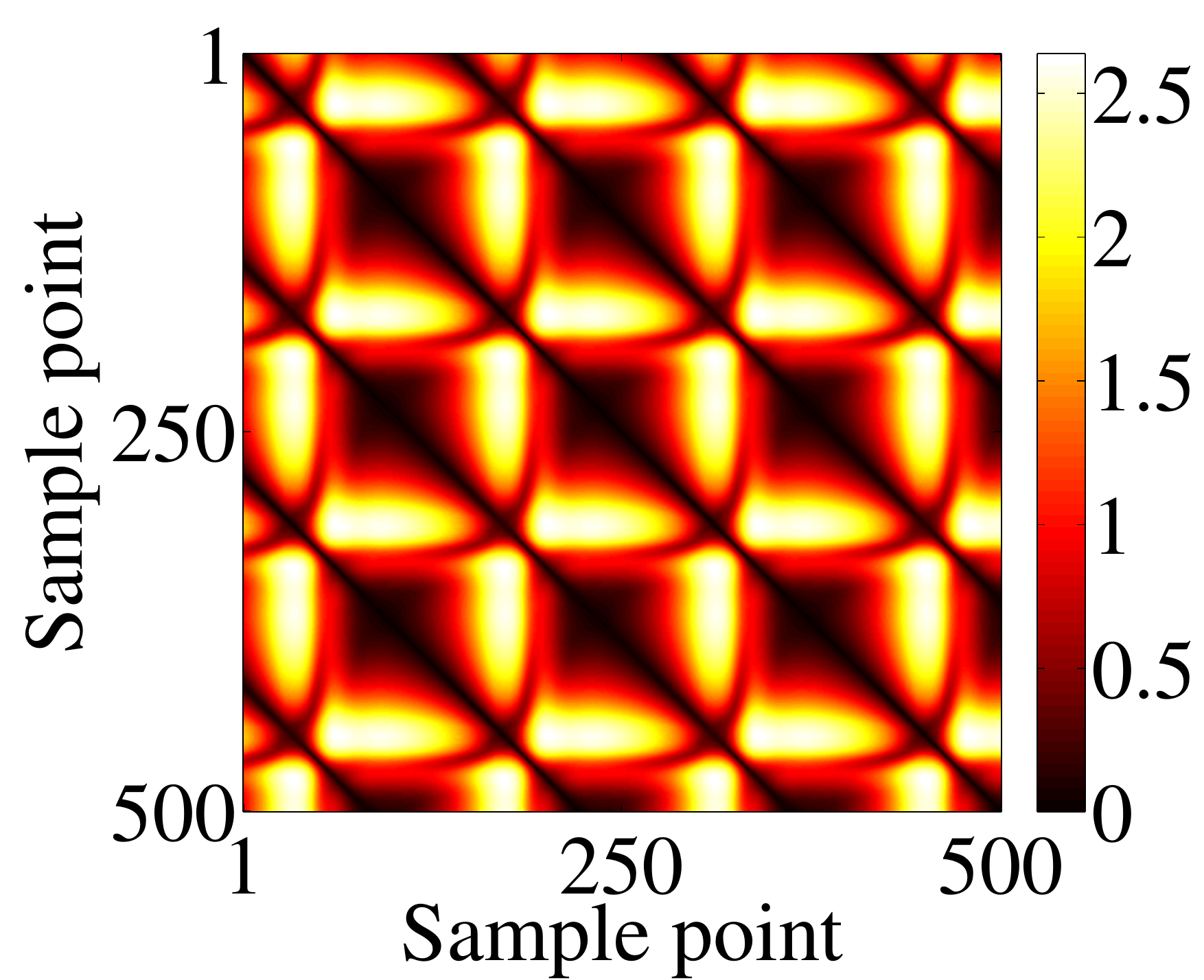}}
\subfigure[]{\includegraphics[width=1.7in]{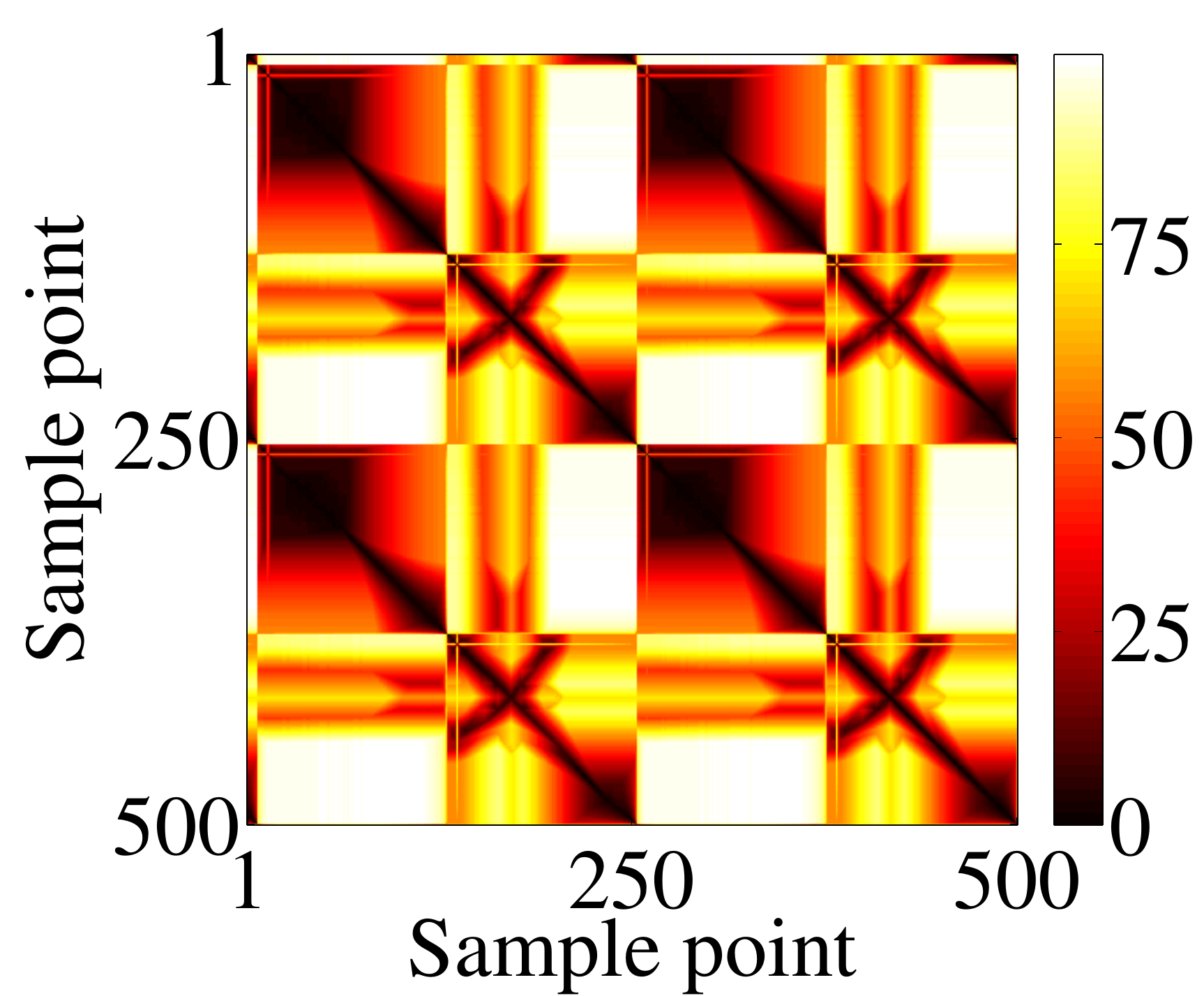}}
\caption{(a) Distance matrix $D$, generated by the $d_{W^2}$ metric,  for approximately three periods of the stable relative periodic orbit  of the Kolmogorov flow. The large black patches correspond to the parts of the orbit with slow dynamics. Equally spaced black lines parallel to the diagonal suggest periodicity of the orbit with period equal to the distance between these lines. (b) Distance matrix $D$, generated by the $d_{B}$ metric, for 2 periods of the almost-periodic orbit of Rayleigh-B\'{e}nard convection. The checkerboard pattern indicates that sampling is too sparse, and fast dynamics are not resolved with the level of sampling.}
\label{fig:DistanceMatrixPeriodicOrbit}
\end{figure}

While the general shapes of the speed profiles  for different distances are similar, there are places where the signs of their derivatives differ. As the system starts accelerating around $t = 100$, all three speeds are increasing.  Around $t = 130$, the speed $s_{W^1}$ starts decreasing while the  other two speeds are still increasing. Note that around $t = 130$, the speeds rise above the noise level (fluctuations). The fact that $s_B$ and $s_{W^2}$ are both increasing means that the changes between the  prominent geometric features are  growing in this region. The speed $s_{W^1}$ is decreasing in this region and so the noise (error) fluctuations are decreasing. At $t = 170$, the dominant geometric features start to evolve considerably. Changes of the dominant features are the most important contributions to all three metrics. Therefore, the derivatives of the speeds $s_\star$ have the same sign again (see video 3, 4, or 5 in the supplementary materials).

Figure~\ref{fig:ConsecDist}(b) shows the normalized speed profiles for the Rayleigh-B\'enard convection simulations. As in the case of the Kolmogorov flow, all three metrics indicate that there are two distinct speed scales along the orbit.  However, the speed profile for $d_B$ differs significantly from those of $d_{W^1}$ and $d_{W^2}$. In particular, it suggests that for significant time periods, the major geometric features of the temperature field vary only slightly, followed by two rapid bursts of change. This can be verified by viewing video 6, 7, or 8 of the supplementary materials. Away from the rapid bursts, $s_B$ is close to $1$. The temperature field has integer values, so the changes cannot be smaller than $0.5$. This implies that  both $s_{W^1}$ and $s_{W^2}$ are dominated by the small fluctuations which are roughly of order $1$.  Hence,  the relative speeds  $s_{W^1}$ and $s_{W^2}$ have essentially the same shape. 

\begin{figure*}[t]
\centering
\subfigure[]{\includegraphics[width=1.4in]{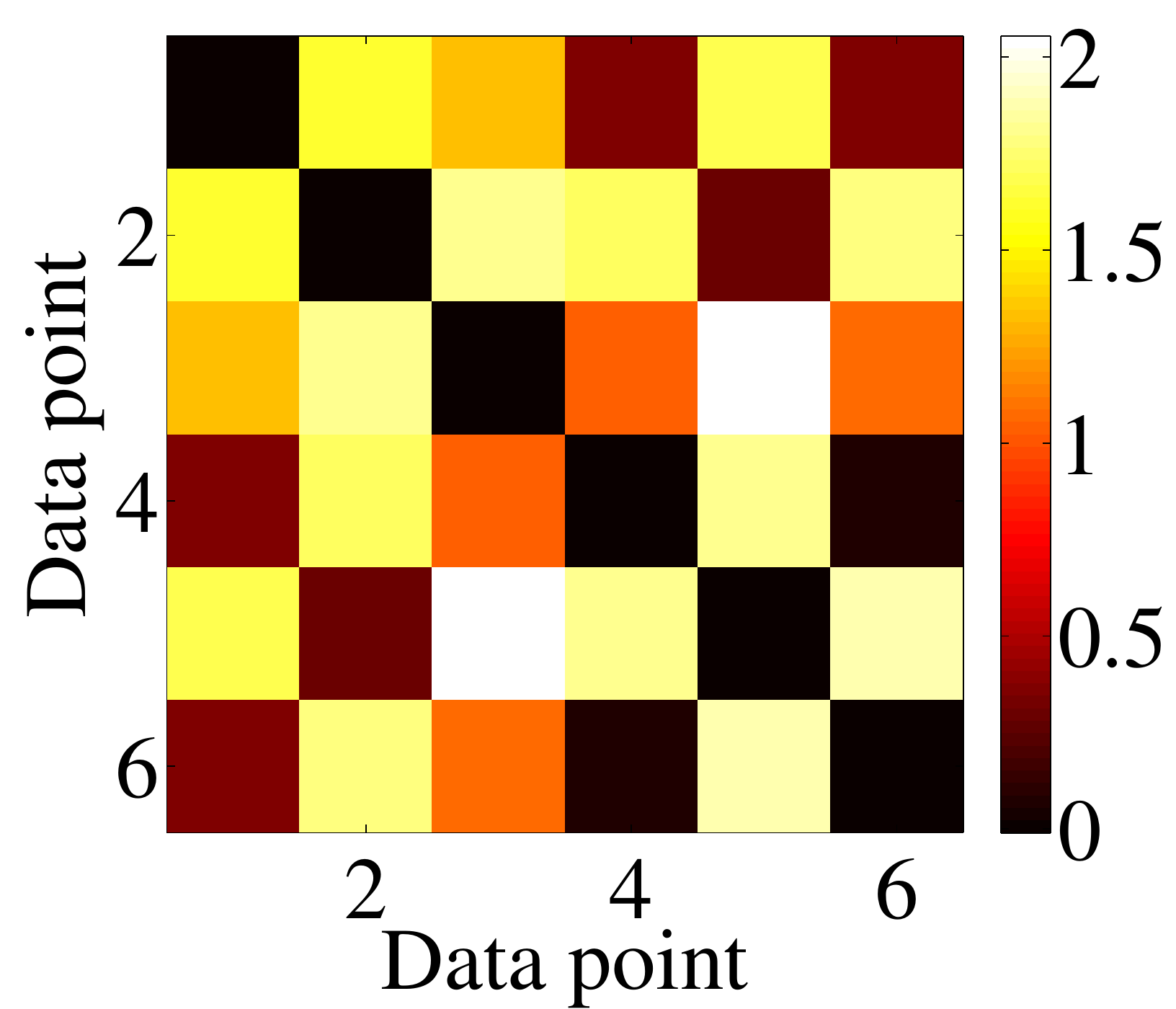}} 
\subfigure[]{\includegraphics[width=1.4in]{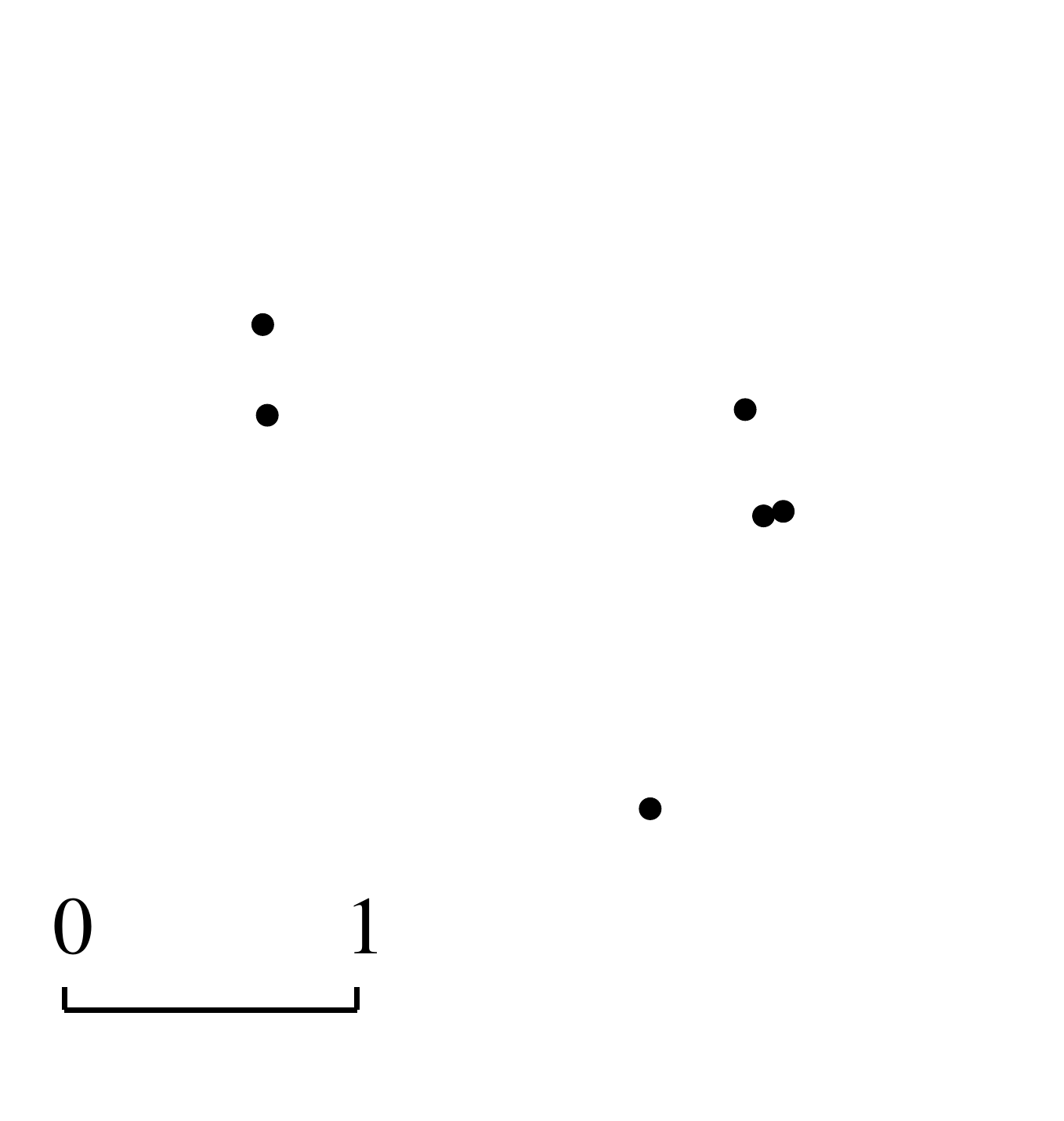}} 
\subfigure[]{\includegraphics[width=1.4in]{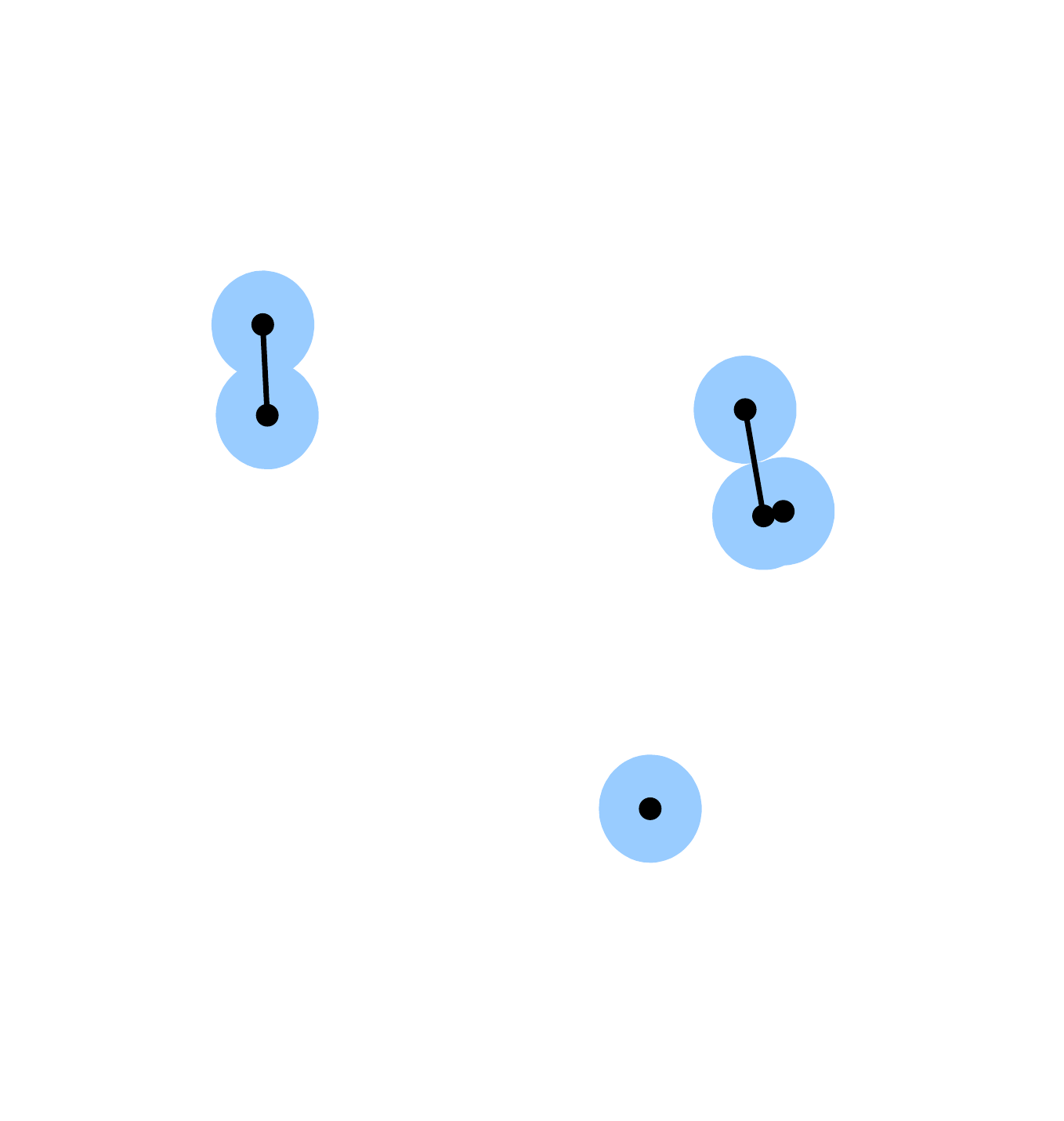}}
\subfigure[]{\includegraphics[width=1.4in]{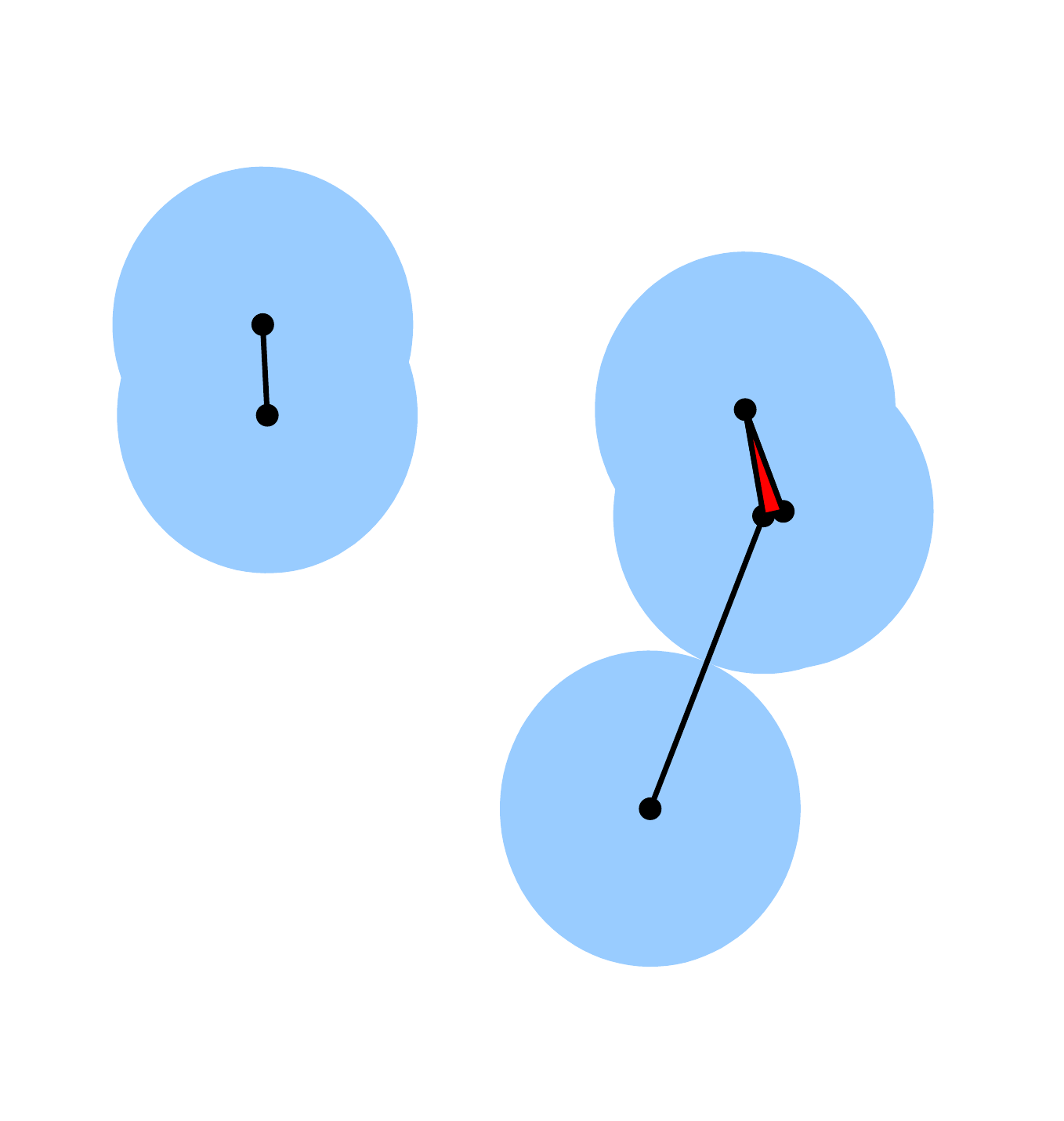}}
\subfigure[]{\includegraphics[width=1.4in]{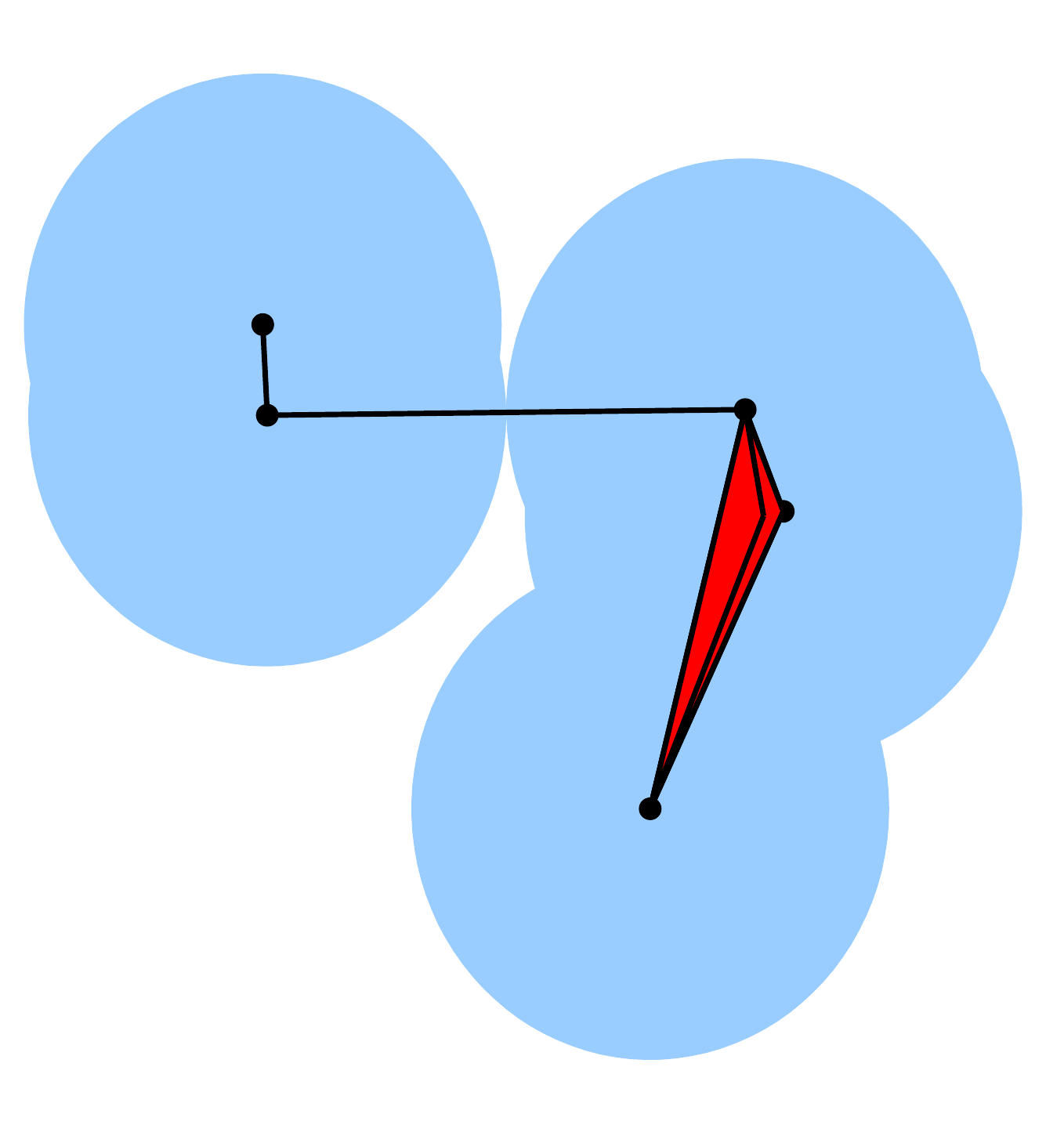}}

\caption{(a) Distance matrix representing pairwise  Euclidian distances $d_E$ between the points in (b) a point cloud $X$.  (c-e)  The blue shaded regions indicate the sub-level sets $\sC(X,\theta)$ for  $\theta = 0, 0.1755, 0.5135$, and $0.816$. The points, edges and triangles indicate the Vietoris-Rips complexes $\sR(X,\theta)$. (c) For $\theta = 0.1775$ the set $\sC(f,\theta)$ consist of three distinct connected clusters. The same is true for $\cR(X,\theta)$. The points in each connected component of $\sC(f,\theta)$ are connected by edges in  $\cR(X,\theta)$. (c) The three components remain distinct until $\theta =0.5135$, at which point two components of  $\sC(f,\theta)$ merge and an edge connecting  the points in the merged components appears in $\cR(X,\theta)$.
}
\label{fig:PointCloudOffset1}
\end{figure*}

The plots of $s_\star$ hint at the underlying dynamics being that of a periodic orbit. However, it is important to keep in mind that $\per$ is an infinite-dimensional space, and thus periodicity in the speed of a trajectory does not imply that the trajectory lies on a closed curve.  As an example, Figure~\ref{fig:ConsecDist}(b) suggests that there are just over four periods of Rayleigh-B\'enard convection shown, and that a single period is roughly 125 frames long. However, looking at video 6, 7, or 8 (in the supplementary materials), it is clear that a full period is closer to 250 frames.  Similarly, it is not obvious that extended periods of high speed imply that the pattern changes significantly over that time period (a periodic orbit of small diameter can exhibit high speed).  This requires a more global geometric analysis of the time series, which we discuss shortly.

With the same data set used to generate Figure~\ref{fig:ConsecDist}(a) and letting $\omega_j$ denote the vorticity field at time $t_j$, Figure~\ref{fig:DistanceMatrixPeriodicOrbit}(a) exhibits the $d_{W^2}$ distance matrix $D$ for Kolmogorov flow, with color-coded entries $D(i,j):= d_{W^2}(\pd(\omega_i),\pd(\omega_j))$. (The $d_B$ and $d_{W^1}$ distance matrices  look very similar and are not shown.) Observe that $D(i,i)=0$ and $D$ is symmetric since $D(i,j) = D(j,i)$. Furthermore, Figure~\ref{fig:ConsecDist}(a) is a plot of the immediate off-diagonal entries. A striking feature of the distance matrix in Figure~\ref{fig:DistanceMatrixPeriodicOrbit}(a) is the existence of dark lines parallel to the diagonal, spaced at intervals of roughly $110$ samples. This indicates that, in the space of persistence diagrams, the trajectory periodically repeats the same, or nearly the same, state. Since the diagonals are spaced at roughly $110$  samples, we can indeed say that the orbit revisits very similar states at intervals of roughly  $110$ samples. Similarly, the light regions close to the diagonal in Figure~\ref{fig:DistanceMatrixPeriodicOrbit}(a) correspond to the times in Figure~\ref{fig:ConsecDist}(a) at which the speed is large, indicating significant changes in the pattern at these times.

To obtain a more global  analysis we turn to Figure~\ref{fig:DistanceMatrixPeriodicOrbit}(b) that shows the distance matrix $D(i,j):= d_{W^2}(\pd(T^*_i),\pd(T^*_j))$ for the temperature fields $T^*_i$ corresponding to the  trajectory from Rayleigh-B\'{e}nard convection.  Distances between the consecutive temperature fields  are shown in Figure \ref{fig:ConsecDist}(b). The dark diagonal lines are spaced at intervals of roughly $250$ samples. Thus, even though the Figure~\ref{fig:ConsecDist}(b) suggests a period of approximately $125$, the orbit does not revisit the same state in the space of persistence diagrams every $125$ samples, but instead every $250$ samples.

\section{Analyzing a Point Cloud using Persistent Homology}
\label{sec:AnalyzingPointCloud}

The discussion in the previous section suggests that interesting information concerning the dynamics of the geometry of  time-evolving scalar fields can be obtained by studying the time series in the space of persistence diagrams.
Note that each scalar field is represented by a persistence diagram $\pd(f)$ and thus corresponds to a point in $\per$.
We argue that viewing the time series as a point cloud in the space of persistence diagrams and studying its geometry provides useful  information about the dynamics.

For a point   cloud  $X\subset\per$  and the scalar function $f\colon X\to [0,\infty)$ given by \eqref{eq:metricScalar} (for any of the metrics $d_B, d_{W^2}$ or $d_{W^1}$), the sub-level set $\sC(f,\theta)$ defined by \eqref{eq:sublevel} is a union of balls 
\begin{equation}
\label{eq:Cf}
\sC(f,\theta) =  \bigcup_{\pd \in X} B(\pd,\theta),
\end{equation}
where $B(\pd,\theta) = \setof{\pd' \in \per \mid d(\pd',\pd)\leq \theta}$, and $d$ is the appropriate metric.  In general, one should expect that the sets $\sC(f,\theta)$ are complicated. Therefore, computing $H_*(\sC(f,\theta))$ directly is not practical. Instead, we make use of the following complex.

\begin{defn}
\label{def:Rips}
{\em
Given a point cloud $X = \{x_0, ..., x_N\}$ in a metric space with distance function $d$,
the {\em Vietoris-Rips complex at scale $\theta$}, denoted $\sR(X,\theta)$, is the simplicial complex defined by the collection of  simplicies
\[
\setof{\ang{x_{n_0}, ..., x_{n_{k}}}\mid d(x_{n_i}, x_{n_j}) \leq 2\theta,\ \text{for all}\ i,j \in \setof{0, 1, 2, ..., k} }.
\] 
}
\end{defn}

Observe that the Vietoris-Rips complex is determined by the distance matrix associated with $X$, and hence, there is a finite set of threshold values $\Theta = \setof{\theta_i}$ at which the complex changes.  Thus, given a point cloud $X$ in a metric space with metric $d$,  the associated persistence diagrams $\pd(X,d)$ are  determined by the Vietoris-Rips complexes $\sR(X,\theta)$ for $\theta \in \Theta$.

We  emphasize that the only data used to analyze a point cloud based on the persistent homology of Vietoris-Rips complexes are the pairwise  distances between the points given by the distance matrix associated with $X$.

\subsection{Detecting Clusters}
\label{sec:Clusters}

Since $\beta_0$ counts components, it is reasonable to use persistent homology as a clustering tool. We demonstrate this idea on a point cloud with pairwise distances  given by the distance matrix shown in Figure~\ref{fig:PointCloudOffset1}(a).  A possible configuration of the six points in $\R^2$ is depicted in Figure~\ref{fig:PointCloudOffset1}(b).  Using the length scale presented in Figure~\ref{fig:PointCloudOffset1}(b) as an indicator of the order of magnitude at which we want to declare a separation length for the clusters,  there are three clusters. We now focus on the geometric information conveyed by $\pd_0(X,d_E)$, shown in Figure~\ref{fig:PointCloudData1}.

Observe that $\sC(f, 0) = \sR(X,0)$ consists of $6$ vertices. As $\theta$ increases, the distinct connected components of $\sC(f,\theta)$ (as defined in \eqref{eq:Cf})  start merging together.  In fact, when the balls $B(x_i,\theta)$ and $B(x_j,\theta)$  merge together, an edge $\ang{x_i,x_j}$ appears in $\sR(X,\theta)$.  Therefore, $H_0(\sC(f,\theta)) = H_0(\sR(X,\theta))$ for all $\theta \in \R$}, and  $\pd_0(f) = \pd_0(X,d_E)$.  Note that it is impossible for a new connected component to appear for $\theta > 0$.  Hence, all persistence points in $\pd_0(X,d_E)$ have a birth value equal to zero.  The death coordinates represent the spatial scales at which distinct connected components (clusters) merge together.  Say that we are interested in clusters where the minimal separation is on the order of  1 length scale. These clusters correspond to the  points in $\pd_0(X,d_E)$ with the death coordinate greater than approximately $0.5$, and there are three persistence points that satisfy this criterion. Thus, we declare that there are three clusters. If  the relevant scale for  separation is of an order of magnitude smaller, then  there are five clusters, since, in addition to the three points with death value greater than $0.5$, two points have death values slightly larger than $0.05$. 

\begin{figure}[]
\centering
{\includegraphics[width=1.7in]{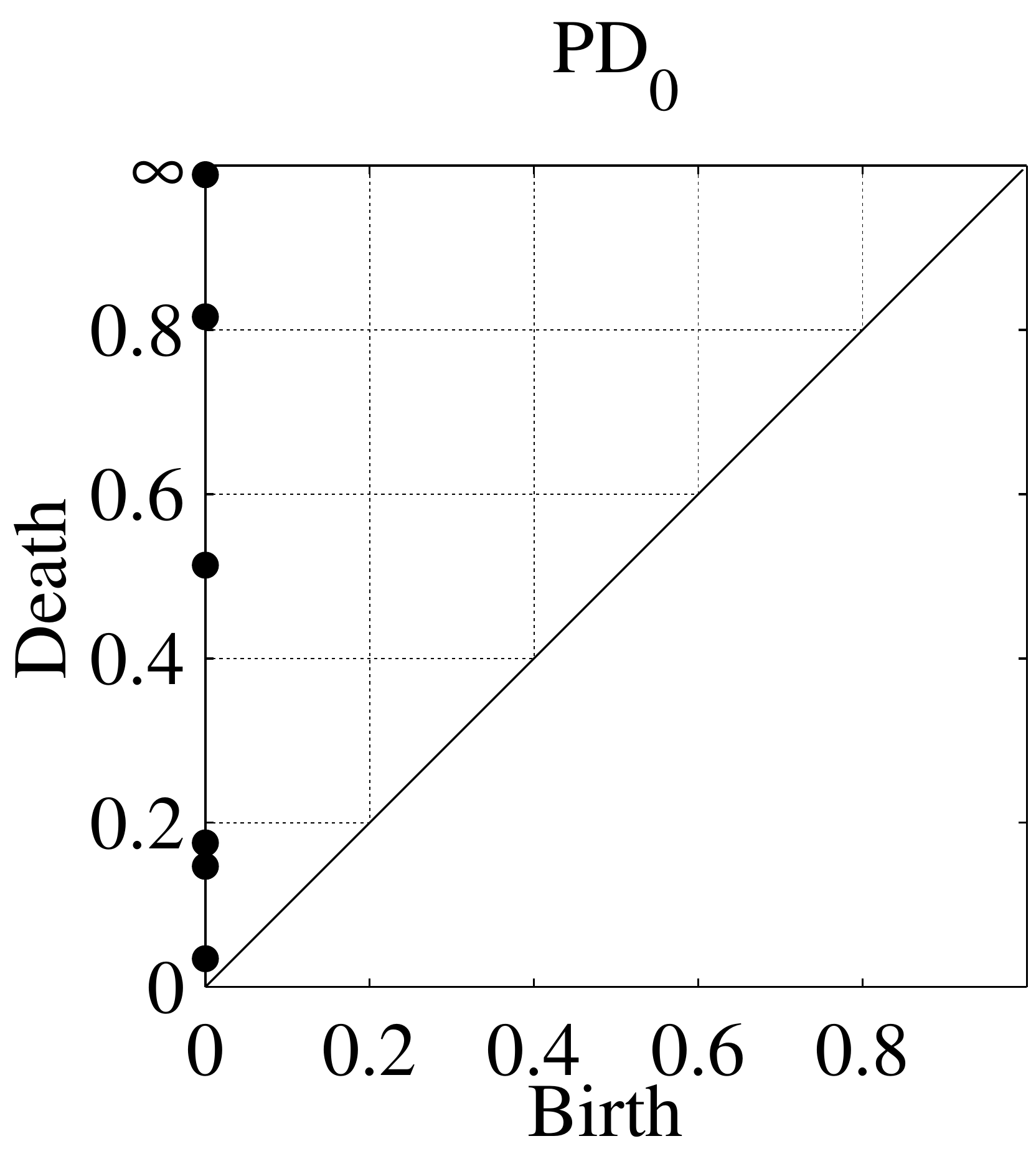}} 
\caption{Persistence diagram  $\pd_0(X,d_E)$ corresponding to  the distance matrix in Figure~\ref{fig:PointCloudOffset1}(a).}
\label{fig:PointCloudData1}
\end{figure}

Alternatively, if we are interested in dividing the data into two clusters, then $\pd_0(X,d_E)$ can  be used to determine the magnitude of the separation between the clusters. Observe that the persistence point $(0,\infty)$ corresponds to the final connected component. The persistence point $(0,0.816)$, with the largest finite death coordinate, indicates that the components merged at a distance $0.816$. Hence, the minimal distance between points from the point cloud $X$ that belong to two distinct clusters is $1.632$.

\subsection{Detecting Loops}
\label{sec:Loops}

\begin{figure*}[t]
\centering
\subfigure[]{\includegraphics[width=1.4in]{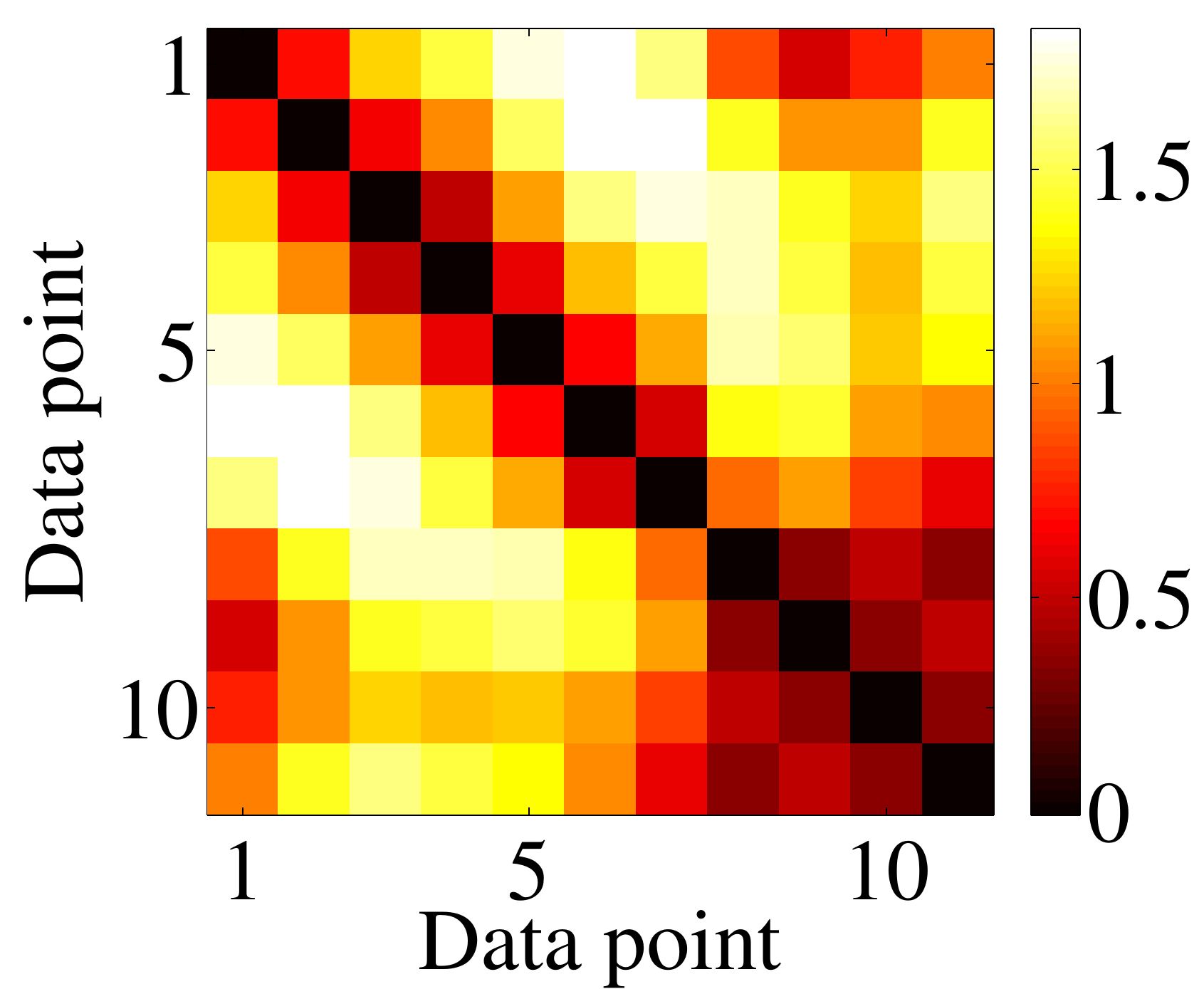}}
\subfigure[]{\includegraphics[width=1.4in]{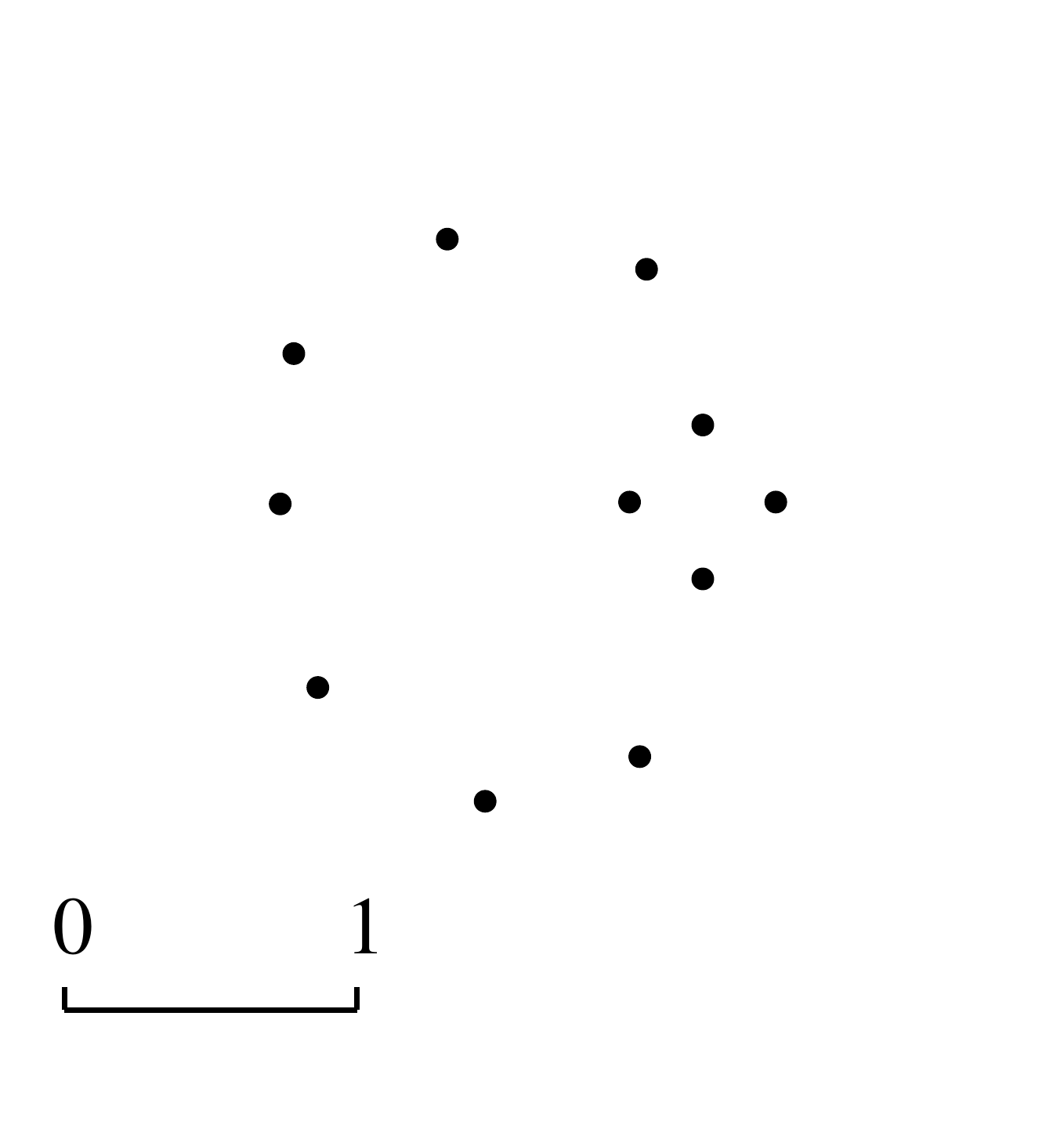}} 
\subfigure[]{\includegraphics[width=1.4in]{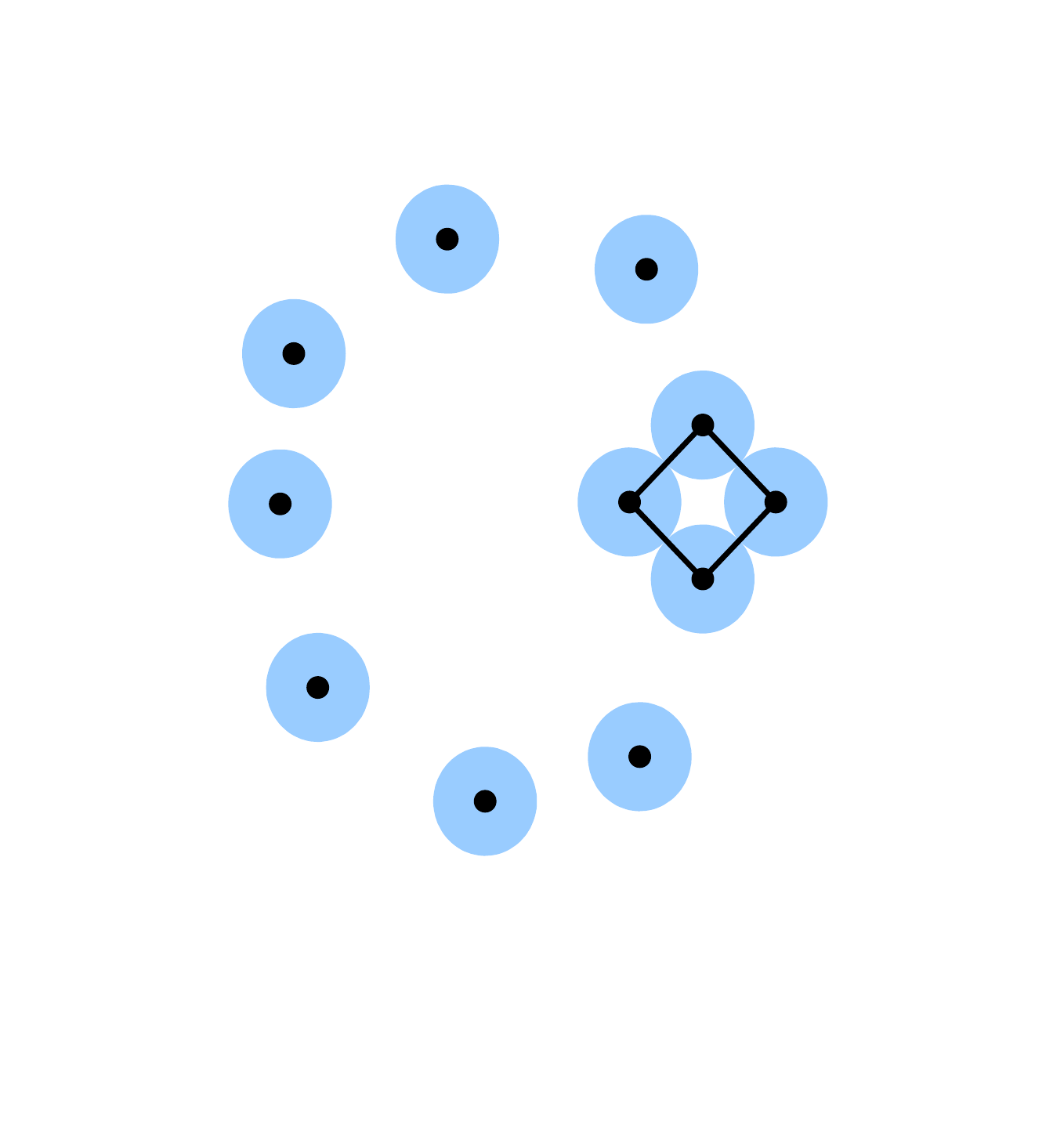}}
\subfigure[]{\includegraphics[width=1.4in]{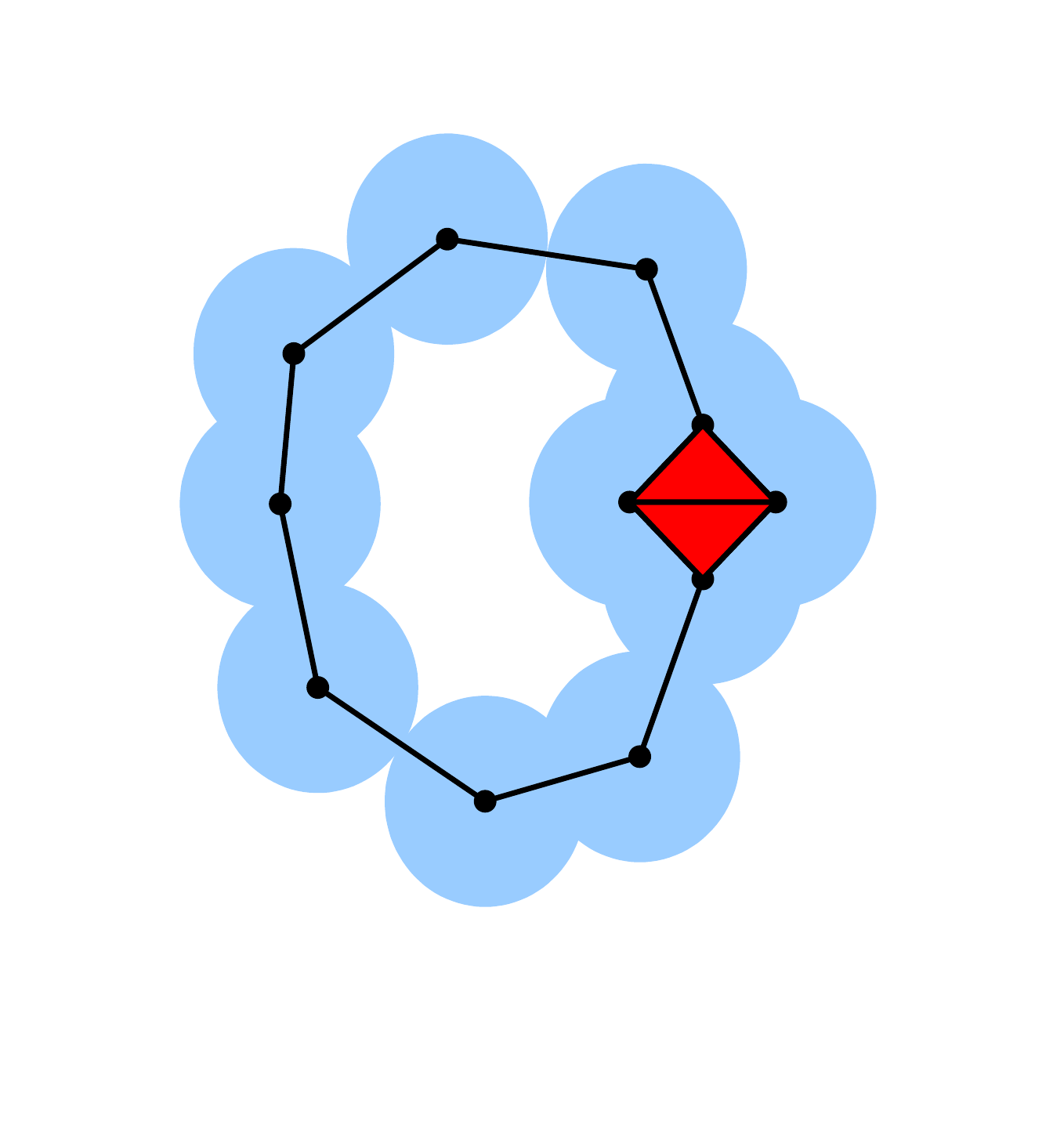}}
\subfigure[]{\includegraphics[width=1.4in]{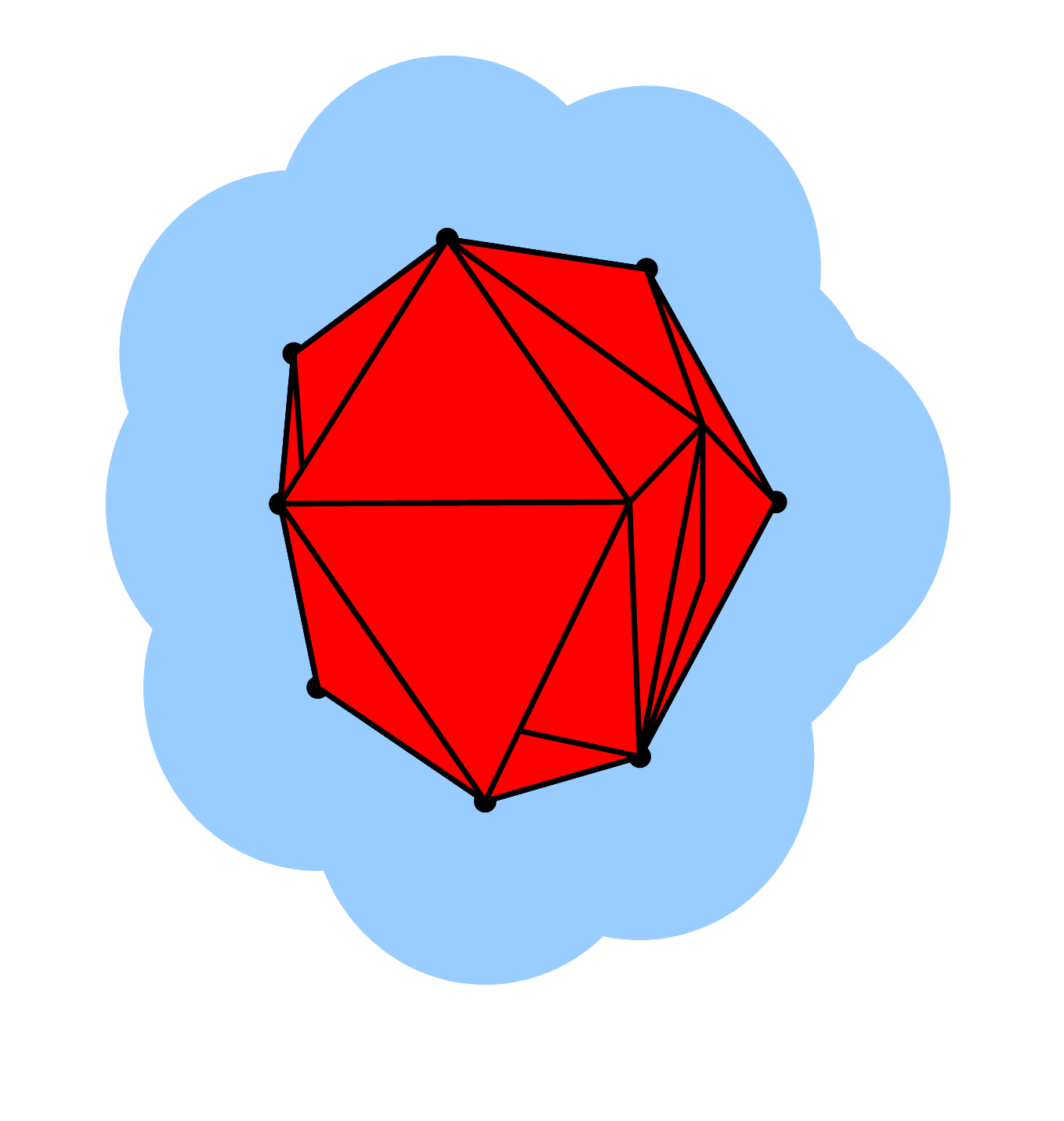}}
\caption{(a) Distance matrix representing pairwise  Euclidian distances $d_E$ between the points in (b) a point cloud $X$.  (c-e) Sets $\sC(f,\theta)$ for  $\theta = 0,  0.177, 0.343$ and  $0.596$. The homology of  $\sC(f,\theta)$ can be approximated by a Vietoris-Rips  complex $\sR(X,\theta)$ given by the vertices, edges, and  triangles shown in (b-e). The first loop in $\sC(f,\theta)$ is created at $\theta = 0.177$. This loop is due to the noisy sampling and is filled in almost immediately.  The dominant loop shown in (c) is formed at $\theta = 0.343$ and persists until $\theta = 0.596$. 
}
\label{fig:PointCloudOffset2}
\end{figure*}

\begin{figure}[]
\centering
\subfigure[]{\includegraphics[width=1.7in]{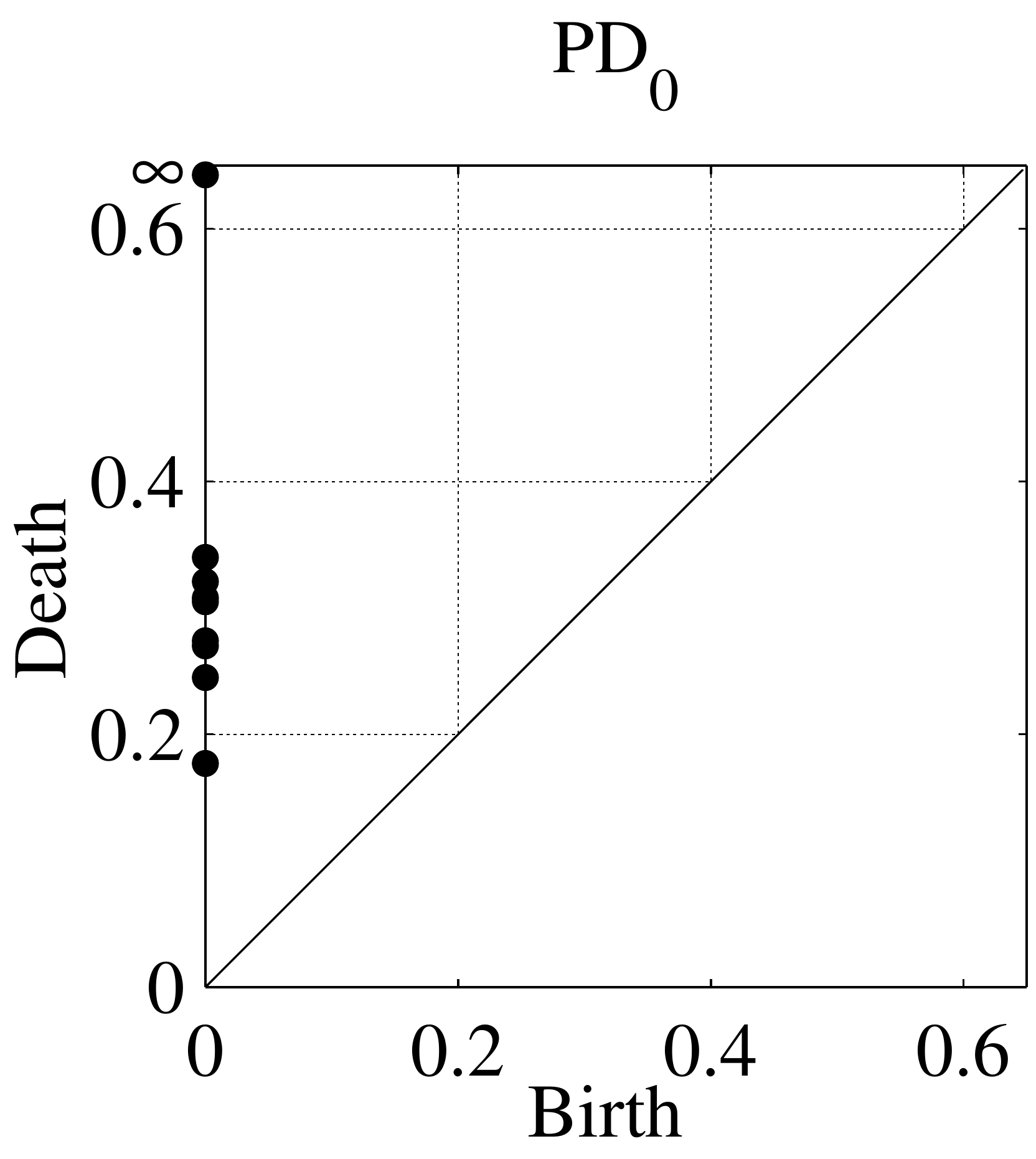}} 
\subfigure[]{\includegraphics[width=1.7in]{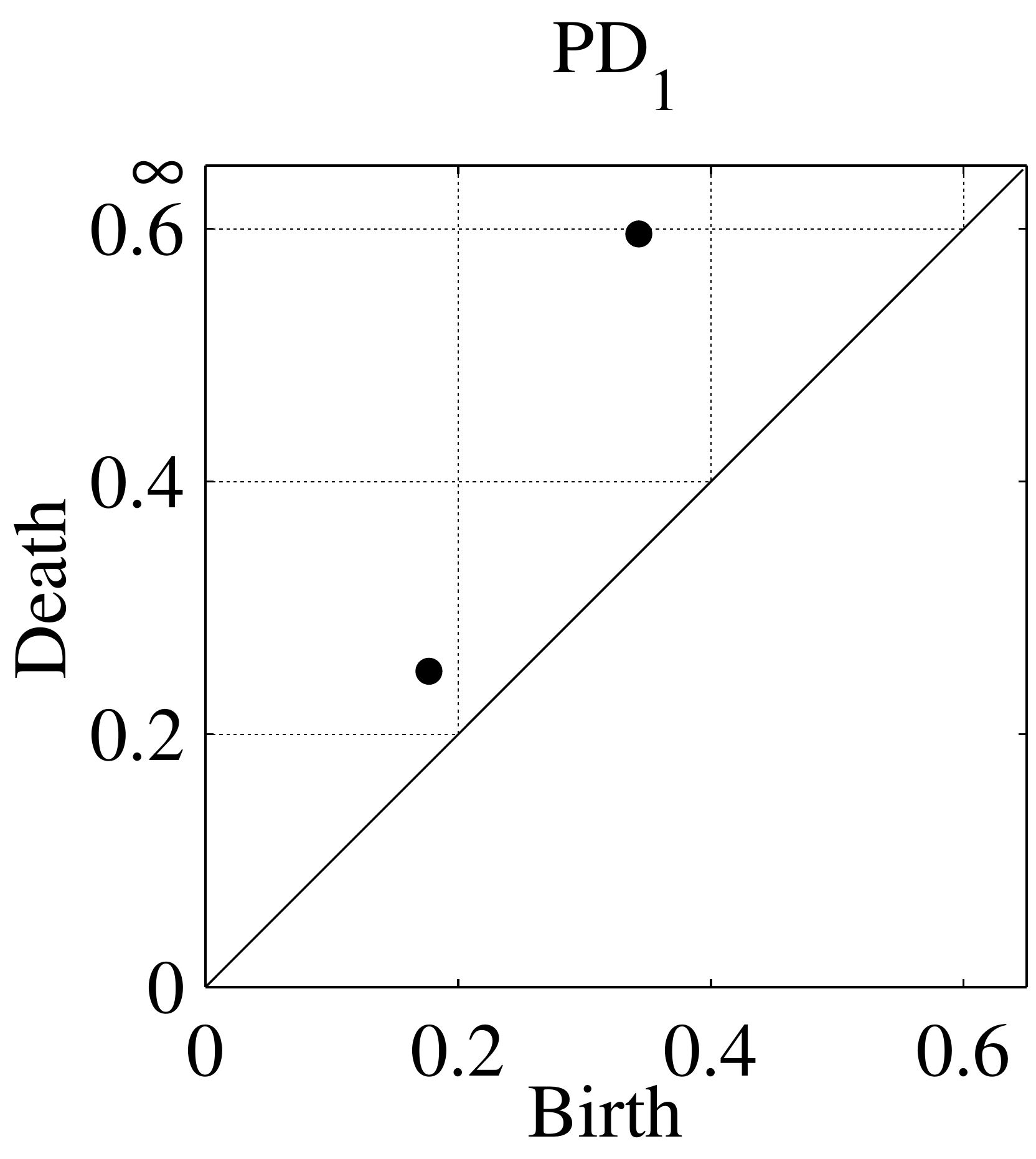}} 
\caption{Persistence diagrams (a) $\pd_0(X,d_E)$ and  (b) $\pd_1(X,d_E)$ corresponding to  the distance matrix in Figure~\ref{fig:PointCloudOffset2}(a). The persistence diagram $\pd_1(X,d_E)$ contains a dominant point $(0.343 , 0.596)$ corresponding to the robust loop  shown in Figure~\ref{fig:PointCloudOffset2}(d) while the point $(0.177 , 0.25)$ represents the small loop visible in  Figure~\ref{fig:PointCloudOffset2}(c). }
\label{fig:PointCloudData2}
\end{figure}

Since $\beta_1$ counts loops, it is reasonable to use persistent homology as a tool for identifying cycles that arise from dynamics. Consider any point cloud that generates a distance matrix as in Figure~\ref{fig:PointCloudOffset2}(a). Again, for the sake of intuition, we present in Figure~\ref{fig:PointCloudOffset2}(b) an example of a point cloud $X\subset \R^2$ with pairwise distances given by the distance matrix shown.  The persistence diagrams for the associated Vietoris-Rips complex filtrations are shown in Figure~\ref{fig:PointCloudData2}.

Applying the reasoning from the previous section, we can ask whether there is a natural or interesting clustering of the data. If, as before, we insist that we are interested in clusters where the minimal separation is on the order of  length scale 1, shown in Figure~\ref{fig:PointCloudOffset2}(b), then $(0,\infty)$ is the only persistence point with death value greater than $0.5$, i.e.\ at this scale there is only one component. Thus, we conclude that from a geometric perspective we may treat the point cloud as arising from a single dynamical structure.

We now look for cyclic structures.  Observe that $\pd_1(X,d_E)$ contains two persistence points. The life span of point $(0.177,0.250)$ is $0.06$, which is short compared to the order 1 length scale, and thus it is reasonable to think of this as a result of noise in the data.  This is substantiated by Figure~\ref{fig:PointCloudOffset2}(c), in which the loop in the Vietoris-Rips complex $\sR(X,0.177)$ consists of four edges.  An additional edge and two triangles (two 2-simplicies) appear in $\sR(X,0.250)$, see Figure~\ref{fig:PointCloudOffset2}(d). The triangles fill in the loop formed by the edges of $\sR(X,0.177)$. Two of the four data points that are involved in the construction of this loop can be viewed as arising from noise or errors associated with sampling points from a smooth cycle.

The life span of persistence point $(0.343, 0.596)$ is $0.253$ and suggests that the point cloud is generated by a loop with a minimal radius of $0.596$, which is on the order of the scale of the data.  This suggests that the associated cycle, indicated in Figure~\ref{fig:PointCloudOffset2}(d), represents an observable, robust dynamical feature.

\subsection{Application to Systems with Multiple Time Scales and Large Data Sets} 
\label{sec::LargeData}

Characterizing the geometry of a continuous orbit via an approximation by a discrete time series depends on the frequency of sampling, and thus becomes a challenge in the setting of dynamics with multiple time scales, i.e. when the rate of change of the patterns is far from constant.   If the sampling rate is too slow, then parts of the orbit will be poorly (or not at all) sampled.  Note that  the geometry of the continuous trajectory may be more complicated than that of a circle; secondary structures might occur if the orbit is twisted, pinched, or bent in $\per$. Thus, the missing parts of the orbit could distort (or entirely miss) significant features in the geometry of the sampled trajectory as compared to the geometry of the underlying (continuous) dynamics.  Thus, in order to obtain a description of the geometry on all relevant spatial scales, including information about secondary structures,  the sampling rate needs to be fast enough.

To determine if a trajectory has been sampled densely enough to resolve the geometry of the underlying dynamics, it is useful to compare the following three values related to the point cloud in $\per$: the noise threshold of the system, the maximum consecutive distance in the sampled trajectory, and the diameter of the point cloud.  Ideally, once a noise threshold has been computed, one would like  distances between consecutive points from the sampled trajectory to be on the length scale of the noise. If sampling faster than this, the features detected from the sample that are on the scale of the noise would be indistinguishable from artifacts generated from the noise in the sample. Thus, ideally, the distance profiles (e.g. Figure~\ref{fig:ConsecDist} for Kolmogorov flow and Rayleigh-B\`{e}nard convection) should have maximums no larger than the noise.  Unfortunately, this is not practical for reasons that will be explained next, and fortunately it is often not necessary.  For example, the length scale of the computational noise could be much smaller than the relevant length scale of interest for studying the geometry of the dynamics.  In this case, a comparison of the maximum consecutive distance in the sample to the diameter of the point cloud in $\per$ is often useful. For instance, if a point cloud has diameter 100 and the smallest relevant length scale for the geometry to be studied is 10, then a maximum consecutive distance of 10 is sufficient for the sampling of the time series,  even if the noise threshold is on length scale 1. Thus, it is the interplay of these three numbers that  determine if one has sampled a continuous time series densely enough.

Evaluating these three quantities from an initial time sample may indicate that an increase in the sample rate is required to resolve the dynamics at the relevant spatial scale.  In the context of a large-scale computation such as that required for the 3D simulation of Rayleigh-B\'{e}nard convection, it is easier to save the data  at a higher sampling frequency than to develop numerical methods that save data based on an adaptive time step.  In Section~\ref{sec:RBC}, we demonstrate the approach introduced here using approximately $7\times10^5$ equally-spaced snapshots of the temperature field of Rayleigh-B\'{e}nard convection.  It should be immediately apparent that the  set  $X$ is too large to compute the associated persistence diagrams $\pd(X,d_\star)$, for $\star \in B, W^2, W^1$, directly.  The first step would require computing the distance matrix for $X$, which would involve $49\times 10^{10}$ distance computations.  Note, however, that using a fast sampling rate leads to collecting unnecessarily many samples at the places where the dynamics are slow.  This suggests that an appropriate choice of down-sampling will allow us to capture the global geometry of the point cloud. 

\begin{defn} 
{\em
Let $X$ be a point cloud in a metric space $(M,d)$. Fix $\delta > 0$.  A set $Y\subseteq X$ is  a $\delta$-\emph{dense} subsample of $X$ if for every $x\in X$, there exits a $y \in Y$ such that $d(x,y) < \delta$.   
}
\end{defn}

The following theorem \cite{harker:levanger:kramar:mischaikow} guarantees that using a $\delta$-dense subsample enables us to detect geometric features with life span larger than  $\delta$.

\begin{thm}
Let $X$ be a point cloud in a metric space $(M,d)$ and $Y$ a $\delta$-dense subsample of $X$. Then $d_B(\pd(X,d), \pd(Y,d)) < \delta$.
\label{thm::subsample}
\end{thm}

\begin{rem}
According to the above theorem, there exists a bijection between the points in $\pd(Y,d)$ and  $\pd(X,d)$ such that the distance between matched points is less than $\delta$. Furthermore, it can be shown that there is a bijection $\gamma : \pd(Y,d) \to \pd(X,d)$ with the following property: if $\gamma(\theta_b, \theta_d) = (\theta'_b, \theta'_d)$ and $\theta'_b \neq \theta'_d$, then $0 \leq \theta'_b - \theta_b < \delta$ and $0 \leq \theta'_d -\theta_d < \delta$. 
\label{rem::subsample}
\end{rem}

To optimize the computational cost, we wish to choose a subsample of the point cloud $Y$ as small as possible.  A point cloud $Y$ is $\delta$-\emph{sparse}  if, for every pair of distinct points $y_1, y_2 \in Y$, the distance $d(y_1,y_2) \geq \delta$. Given a point cloud $X$ and a value $\delta \geq 0$, a $\delta$-dense, $\delta$-sparse subsample $Y$ may always be constructed \cite{dey}. Due to the size of the point cloud $X$ and the complexity of computing $d_\star$ for $\star \in B, W^2, W^1$, we use an alternate algorithm \cite{harker:levanger:kramar:mischaikow}, which takes advantages of parallel computing structures and metric trees.

\section{Distinguishing Equilibria}
\label{sec:FixedPoints}

 \begin{figure}[t!]
\centering
\subfigure[]{\includegraphics[width=1.9in]{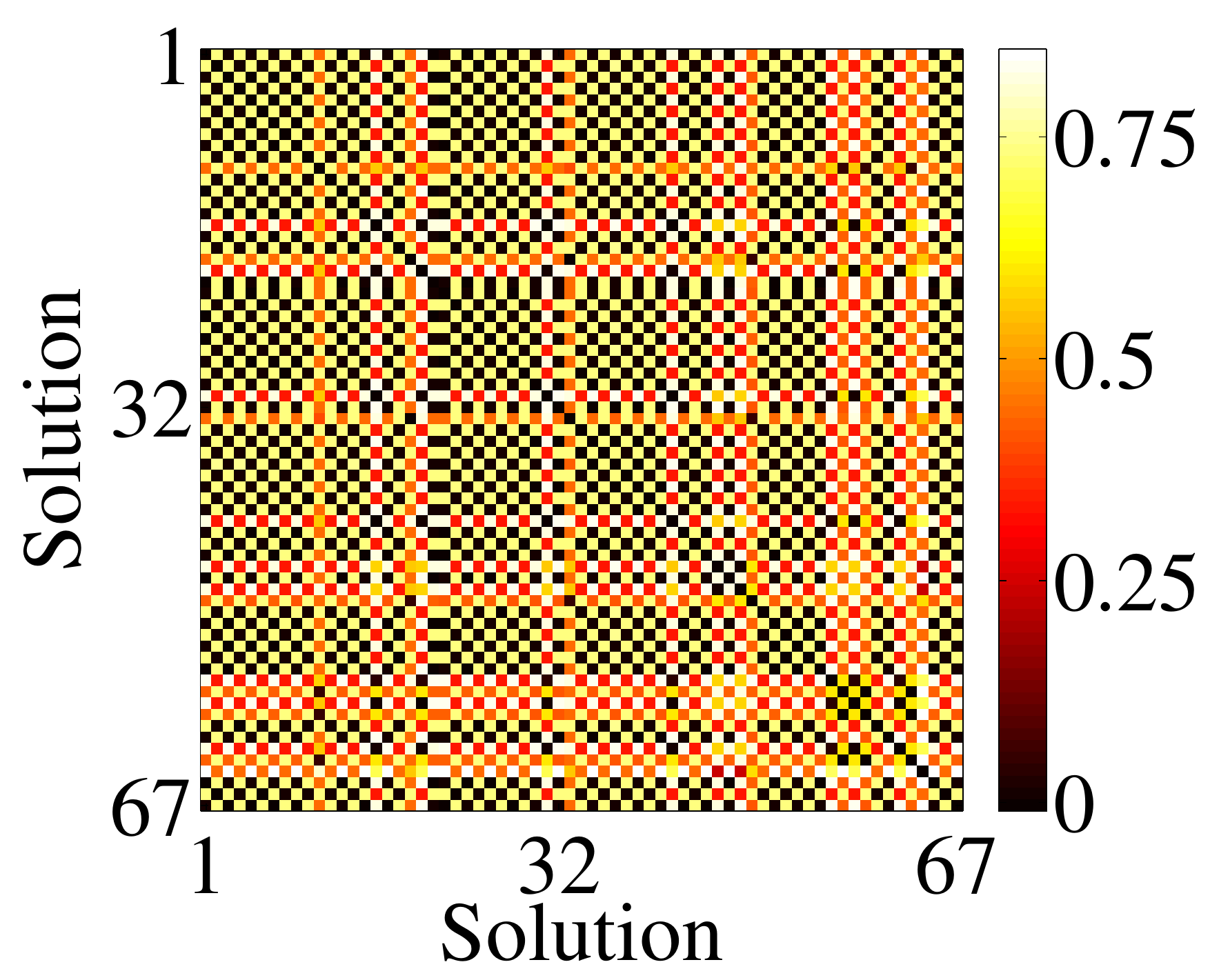}}
\subfigure[]{\includegraphics[width=1.5in]{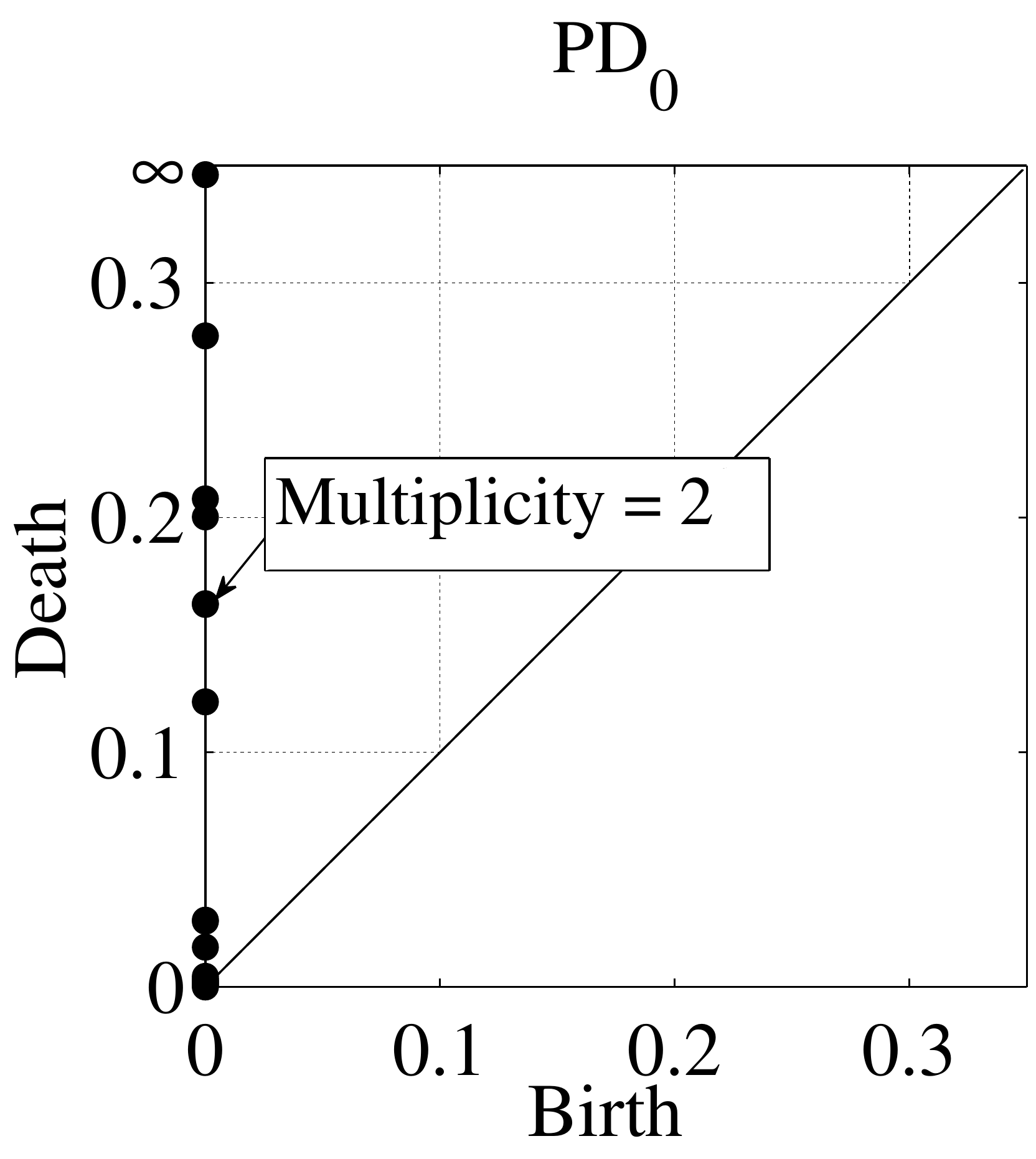}}
\subfigure[]{\includegraphics[width=1.9in]{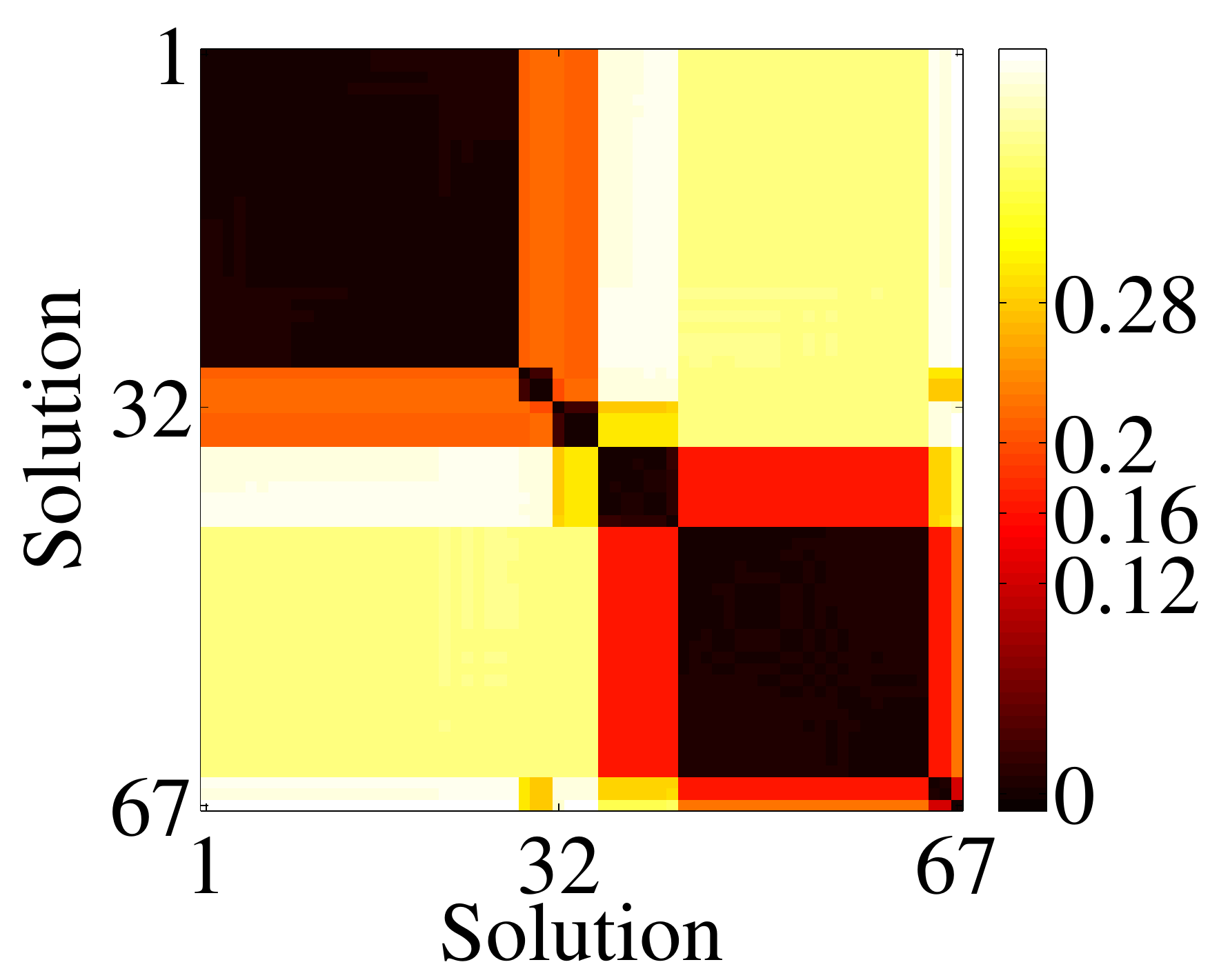}}
\caption{(a) Pairwise $d_{B}$ distances between the  EQ and REQ points in $X =  \setof{\omega_n \mid n = 1,\ldots, 67}$ of the Komogorov flow found at $Re = 26.43$ using Newton's method. (b) Corresponding persistence diagram  $\pd_0(X,d_B)$.  (c) The distance matrix (with values halved to match the values in $\pd_0(X,d_B)$) sorted so that equilibria within each of the seven clusters, detected by $\pd_0(X,d_B)$,  are grouped together. The seven black blocks on the diagonal correspond to the seven clusters, while off-diagonal blocks correspond to values of $\theta$ at which the clusters merge. 
}
\label{fig:DiffFixedPoints}
\end{figure}

We now apply the ideas presented in Section~\ref{sec:Clusters} to the problem of clustering symmetry-related equilibria of the Kolmogorov flow at $Re = 26.43$. As discussed in Section~\ref{Intro_Kolmogorov}, we sample a turbulent trajectory, shown in Figure~\ref{fig:projections}(b).   To obtain the EQ and REQ solutions, we use a Newton method.  The initial guesses for the Newton method are the vorticity fields $\omega$ that are local minima of  the $L^2$ norm of $\partial \omega/\partial t$. In this way, we find a collection  X = $\setof{\omega_n \mid n = 1,\ldots, 67}$ of vorticity fields corresponding to EQ and REQ of the Kolmogorov flow. These 67 solutions may be related to one another through any composition of the coordinate transformations listed in Section~\ref{Intro_Kolmogorov}.  Hence, it is non-trivial to determine how many unique classes of solutions there are and which solutions belong to which class. To perform this analysis, we use persistent homology.

We start by analyzing $\pd_0(X,d_B)$. The pairwise distances between the points in $X$ are shown in Figure~\ref{fig:DiffFixedPoints}(a). As is discussed in Section~\ref{sec:SpaceOfPD}, the distance between persistence diagrams of vorticity fields related by symmetry is small, while  persistence diagrams corresponding to the vorticity fields that are not symmetry related differ by a larger amount. This implies that we can reformulate the question of identifying symmetry classes of equilibria as a clustering problem.

The persistence diagram $\pd_0(X,d_B)$, depicted  in Figure~\ref{fig:DiffFixedPoints}(b), shows a clear gap between the persistence point with death value $\theta_d = 0.0285$ and the persistence point with death value $\theta_d = 0.1215$.  We interpret this gap as separation between the signal and noise (numerical errors). Indeed, $0.0285$ is just twice the estimate of the lower bound on numerical errors for the Kolmogorov flow obtained in Section~\ref{sec:Distances}. There are $7$ points in $\pd_0(X,d_B)$ with death coordinate greater than $0.12$, and so we  conclude that there are seven distinct symmetry classes of solutions.

Grouping the symmetry related equilibria corresponding to seven different clusters and reordering the distance matrix enable us to see how many solutions are in each cluster. This is done by thresholding the distance matrix so that entries greater than a certain value (in this case $\theta = 0.12$) are zeroed out, and then viewing the resulting distance matrix as an adjacency matrix, from which it is possible to determine connected components. The seven black diagonal blocks of the matrix $D(i,j) = d_B(\omega_i, \omega_j)/2$, shown in Figure~\ref{fig:DiffFixedPoints}(c), represent pairwise distances between the symmetry-related solutions in each cluster. The inter-cluster distances are given by off diagonal blocks. Reordering the distance matrix and dividing its entries by two makes it easier to tie its values to the death coordinates of the points in $\pd_0(X,d_B)$.

The values of the off-diagonal blocks between the first three blocks on the diagonal are all roughly 0.2, implying that each of these clusters will merge together at approximately $\theta \approx 0.2$. This behavior is captured by the persistence points $(0, 0.20035)$ and $(0, 0.208)$ in $\pd_0(X,d_B)$. Recall that the merging of three connected components causes the death of just two of them. The next three blocks on the diagonal have off-diagonal blocks with values at roughly $0.16$, and the deaths of two of these underlying components correspond to the persistence points $(0,  0.163)$ and $(0, 0.1633)$. The last diagonal block has a distance of roughly $0.12$ from the sixth block, and the merging of their underlying clusters corresponds to persistence point $(0, 0.1215)$. Thus, at the cutoff value $\theta = 0.21$, there are two components in the dataset, with one corresponding to the first three diagonal blocks and the other the last four. The off-diagonal block that relates the second and sixth diagonal blocks has a value of roughly $0.28$. This distance corresponds to the persistence point $(0, 0.2775)$, at which time the two large components merge together.

We validate the results of the persistence homology analysis by performing clustering using the Fourier amplitudes as follows. If $\hat{\omega}(k_x,k_y)$ is the Fourier amplitude of a mode ($k_x,k_y$), then a translation of the pattern in the $x$ or $y$ directions in real space merely adds to the phase of $\hat{\omega}(k_x,k_y)$, leaving the magnitude unchanged. Hence, by comparing the amplitudes of the Fourier modes we could group vorticity fields which are related by translations. Since the conjugate modes $\hat{\omega}(\pm k_x,\pm k_y)$ relate fields which are related by inversion, to group the vorticity fields which are related by a combination of inversion and translation, we sum the amplitudes of the conjugate modes. Adding the amplitudes of conjugate modes yields a ``reduced matrix," which is unique for all the vorticity fields related by the coordinate transformations that leave Equation \eqref{eq:q2dvor_nd} invariant. This approach also yields 7 distinct classes.

An analysis of $\pd_0(X,d_{W^p})$, $p=1,2$, yields  the same results. There are several  gaps between the death values of the points in the persistence diagrams. Again one of the gaps starts at roughly  twice the value of  the estimated lower bound of the noise.  However, the separation is less pronounced. As discussed in Section~\ref{sec:Distances},  the $d_{W^p}$ metrics capture all the differences between the persistence diagrams, and the local numerical errors are summed together. Thus,  a large number of small errors can obscure the distinction between the signal and noise.

\section{Stable Periodic Orbit of the Kolmogorov Flow}
\label{sec:PeriodicOrbit}

In the previous section, we demonstrated the practicality of using persistent homology to cluster equilibria that are symmetry-related. In this section, we extend these ideas to the setting of recurrent orbits in the context of the Kolmogorov flow.

As is indicated in Figure~\ref{fig:projections}(a), the projection of the orbit onto the real parts of the three dominant eigenvectors suggests a periodic orbit that is undergoing a slow drift in the direction of the continuous symmetry. The nature of this drift is reinforced by tracking this orbit in the space of persistence diagrams; since persistent homology is invariant under the continuous symmetry, this type of drift  is not present in $\per$.   As a result, we expect the time series to lie on a closed loop in $\per$.  This is consistent with the information provided by the distance matrix of Figure~\ref{fig:ConsecDist}, in which the dark lines parallel to the diagonal indicate that the distance between persistence diagrams becomes very small at regular time intervals.

For the remainder of this section, we use the ideas of Section~\ref{sec:Loops} to verify that a circle provides a good description of  geometry of the point cloud $X \subset \per$ generated by the time series sampled from the Kolmogorov flow.  More precisely, we show that there is a single dominant feature in $\pd_0(X,d_B)$ and a single dominant feature in $\pd_1(X,d_B)$, which agrees with the persistent diagrams for a circle.

There are two issues that need to be considered: the first is the size of the data set, and the second is the spacing between the data points. As is indicated in Section~\ref{sec:AnalyzingPointCloud}, we use the Vietoris-Rips complex to compute persistent homology of point clouds. We remark that given $N$ data points, the full Vietoris-Rips complex has $2^N$ cells. Considering this, we complete our analysis with the distance matrices corresponding to $d_B, d_{W^1},$ and $d_{W^2}$ for 500 points, or roughly three periods of the Kolmogorov flow. In the next section, we introduce techniques for computing persistence on larger point clouds, which could  arise  due to increased sampling rates, sampling more periods, or both.

\begin{figure}[t]
\centering
\subfigure[]{\includegraphics[width=1.7in]{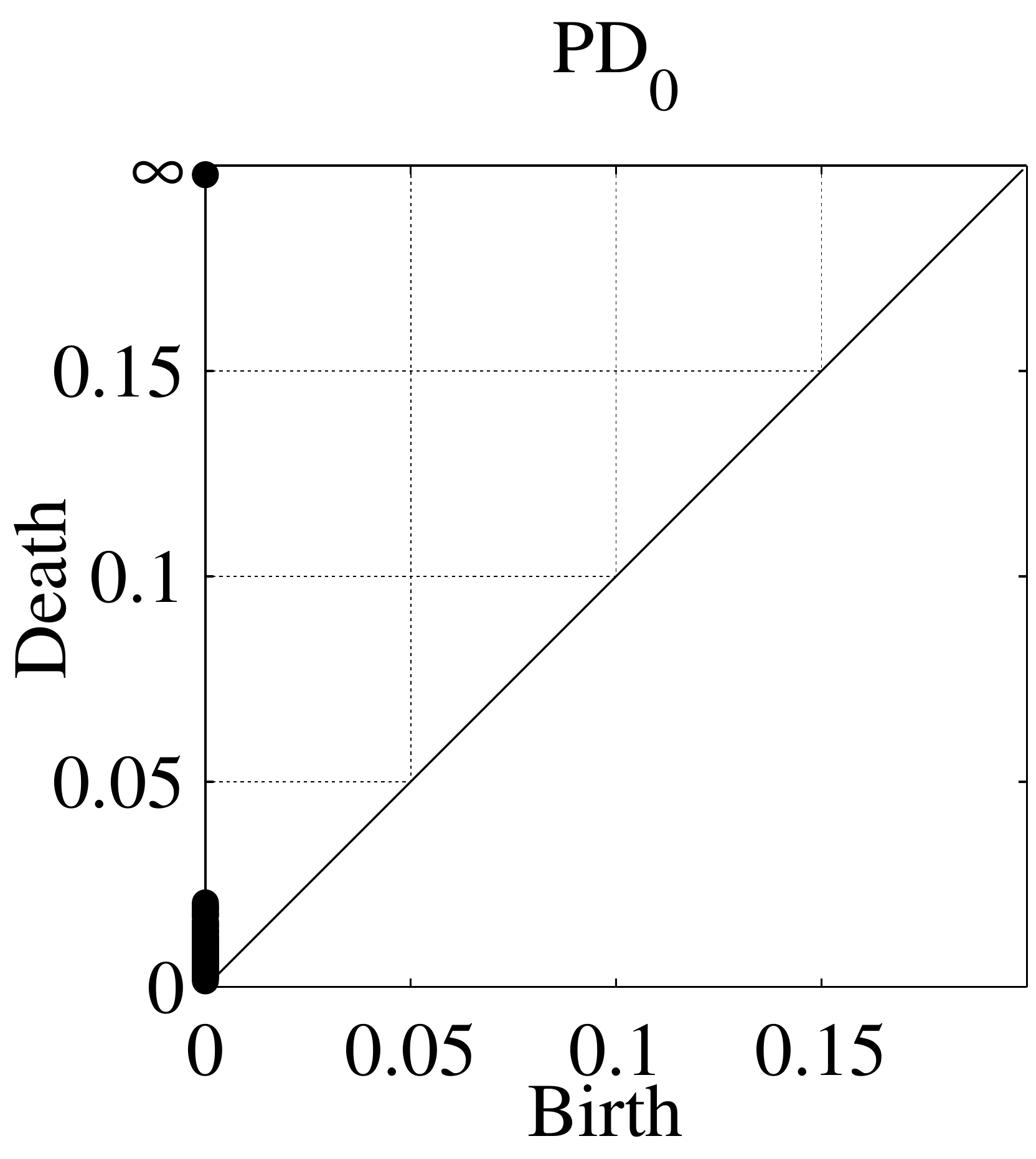}} 
\subfigure[]{\includegraphics[width=1.7in]{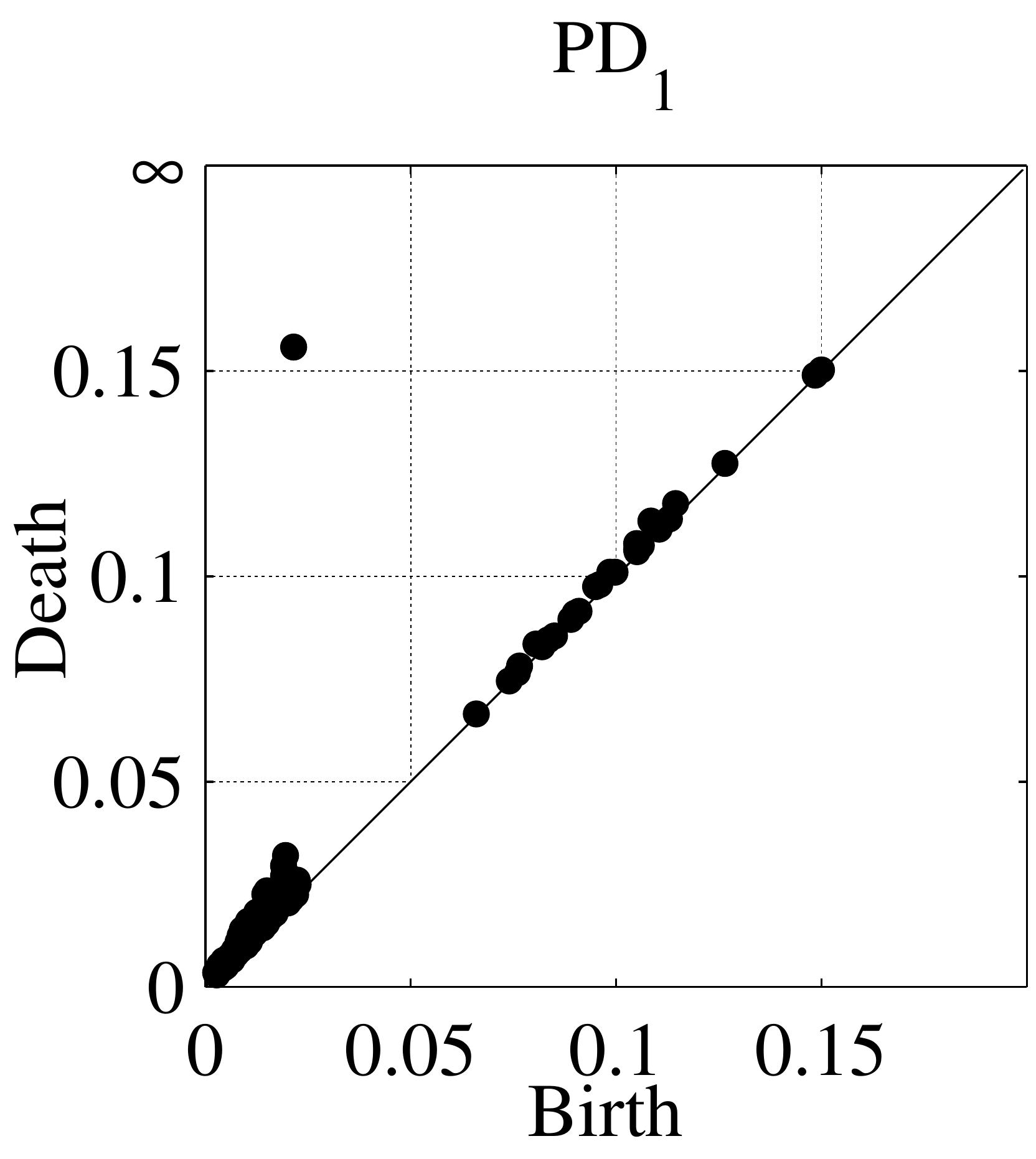}}
\caption{(a) The persistence diagram $\pd_0(X,d_B)$ for Kolmogorov flow at $Re = 26.43$. Since all points with finite death coordinates die before 0.025, there is only a single dominant point.
(b) The persistence diagram $\pd_1(X,d_B)$, showing the single dominant generator at $(0.0215, 0.1559)$.}
\label{fig:25mAOrbitPD}
\end{figure}

Since we are sampling from a single continuous trajectory, the fact that $\pd_0(X,d_B)$, as  shown in Figure~\ref{fig:25mAOrbitPD}(a), suggests the existence of a single component does not come as a surprise. The persistence diagrams for $\pd_0(X,d_{W^p})$, $p=1,2$, yield similar results and are not shown. However, it is worth noting that this is not a foregone conclusion as the location of and spacing between the points of the time series are dependent upon the speed along the periodic orbit. As is clear from Figure~\ref{fig:ConsecDist}(a), the speed of the trajectory is not constant. However, it is fairly smooth, thus we do not expect extreme differences in the spacings between points. 

As discussed in Section~\ref{sec::LargeData}, we compare the noise threshold, $\theta = 0.01$ (Table~\ref{table:ThreeNetworksFullDistance}), to the maximum consecutive sample distance, $d_B = 0.0654$ (Figure~\ref{fig:ConsecDist}(a) caption), and the diameter of the point cloud, $2.64$ (Figure~\ref{fig:DistanceMatrixPeriodicOrbit}(a)). The maximum consecutive sample distance is more than six times larger than the length scale of the noise for this system.  However, the diameter of the point cloud is more than forty times larger than the consecutive sample distance.  Thus, features on the length scale of one fortieth of the diameter of the entire point cloud will be resolved, which is sufficiently small to consider this an adequate sampling.
We will return to this issue in the next section.

As indicated in Figure~\ref{fig:25mAOrbitPD}(b), the persistence diagram $\pd_1(X, d_B)$  clearly detects a single dominant loop along which the data is organized. 
Thus, we conclude that in $\per$, equipped  with the metric $d_B$, the  point cloud $X$ generated by the time series forms a loop with a minimal radius of $0.1344$. Table~\ref{table:25mA_H1_Summary} shows the coordinates of the persistence point with the longest  life span, its life span, and the second longest life span for each of the persistence diagrams $PD_1(X,d_{\star})$, $\star \in B, W^2, W^1$. As the table indicates, the life span of the dominant point is an order of magnitude larger than the next longest life span in each case, and so there is a single dominant feature in $\pd_1(X,d_{\star})$. Additionally, note that the second longest life spans  are as small or smaller than the lower bounds on numerical errors indicated by the first row of Table~\ref{table:ThreeNetworksFullDistance}.

\begin{table}
\centering 
\begin{tabular}{c  | p{2cm} c c   }
\hline
$\pd_1$ & Dominant coordinate  & Max life span &$2^{nd}$ largest life span   \\
\hline
$d_B$ & $(0.022, 0.156)$  & $0.134$  & $0.013$  \\ 
$d_{W^2}$ & $(0.075, 0.405)$  & $1.366$  & $0.105$ \\ 
$d_{W^1}$& $(0.703, 2.069)$ & $0.330$  & $0.016$ \\
\end{tabular}
\caption{The coordinate of the dominant point in the persistence diagram  $\pd_1(X,d_{\star})$ for $\star = B, W^2, W^1$, its life span, and the second largest life span. }
\label{table:25mA_H1_Summary}
\end{table}

\section{Almost-Periodic Orbit of Rayleigh-B\'enard Convection}
\label{sec:RBC}

As mentioned in  Section~\ref{sec::LargeData}, characterizing the geometry of a continuous trajectory becomes a challenge in the setting of dynamics with multiple time scales. To demonstrate this, we consider the numerical simulation of Rayleigh-B\'{e}nard convection, where from multiple perspectives it appears that the trajectory  is close to a periodic orbit and that the rate of change in the patterns of the temperature field is far from constant. This can be clearly seen visually (see video 6, 7, or 8 in the supplementary materials). Moreover, both the speed plot, Figure~\ref{fig:ConsecDist}(b), and the distance matrix, Figure~\ref{fig:DistanceMatrixPeriodicOrbit}(b), suggest recurrent dynamics. However, we note that the rate of change, especially using the bottleneck distance, is typically small except for short periods of time at which the speed spikes. The distance matrix has a distinct checkerboard  pattern, with the edges corresponding to the spikes, again indicating a rapid and large change in location in the space of persistence diagrams.

\begin{figure}[]
\centering
\subfigure[]{\includegraphics[width=1.7in]{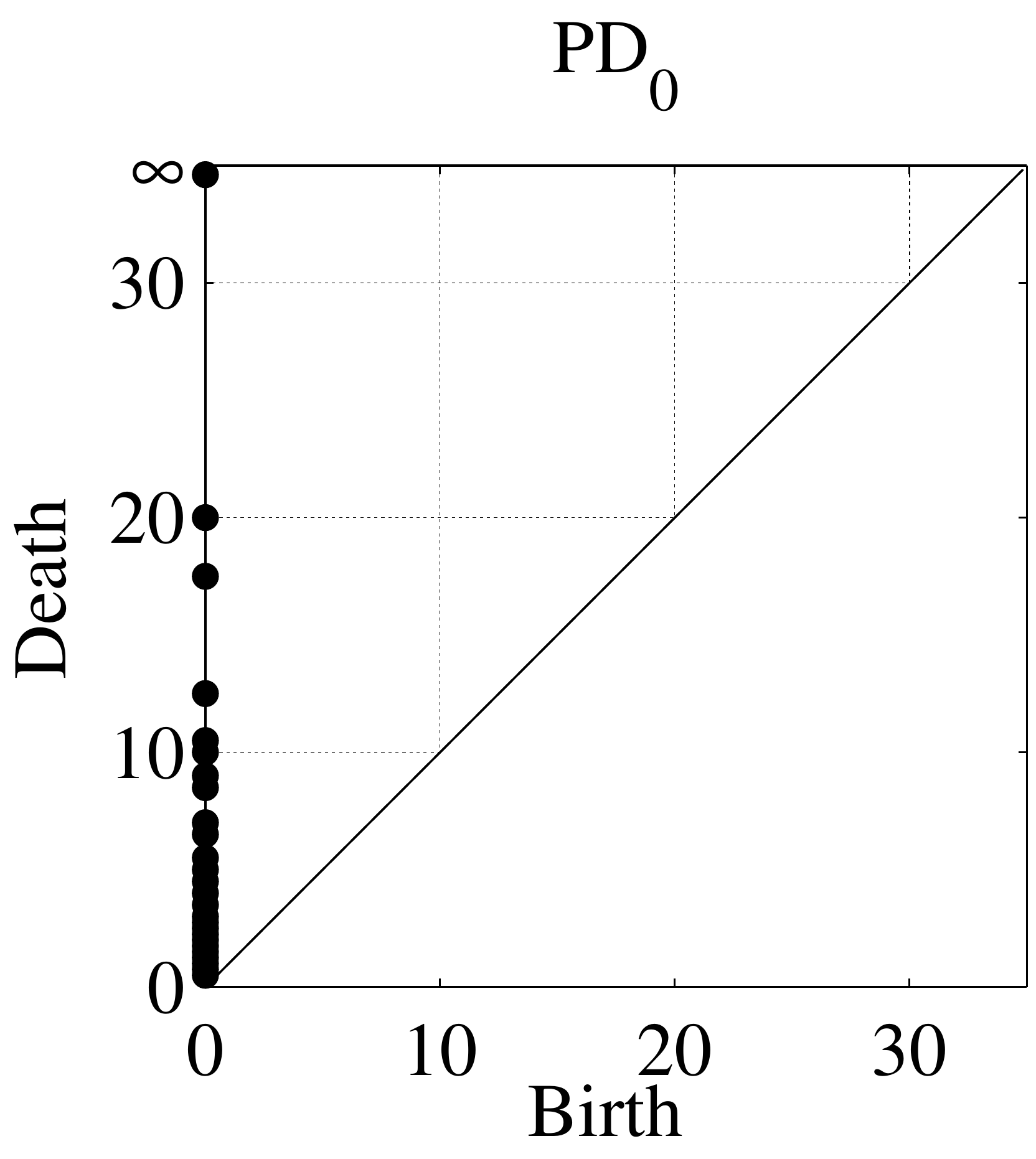}}
\subfigure[]{\includegraphics[width=1.7in]{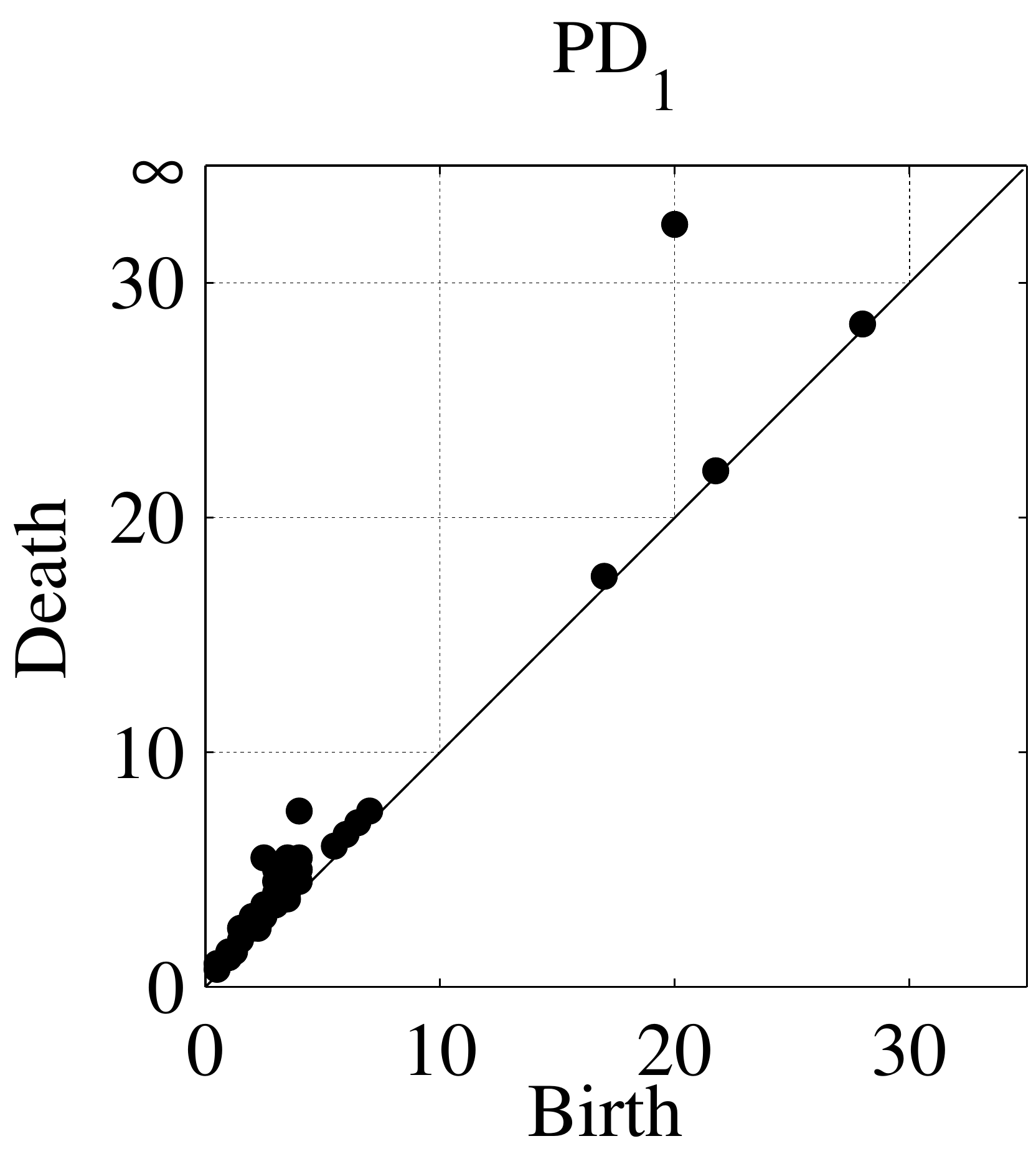}}
\caption{ Persistence diagrams for 500 points of Rayleigh-B\'{e}nard convection, generated from the distance matrix shown in Figure~\ref{fig:DistanceMatrixPeriodicOrbit}(b). (a) The diagram $\pd_0(X,d_B)$ shows the appearance of persistence points with death values significantly greater than the noise threshold, indicating that the sampled trajectory is broken into pieces and sampling is not fast enough to resolve the periodic orbit. (b) The diagram $\pd_1(X,d_B)$  shows the presence of a loop that is born when the pieces of the orbit merge together. } 
\label{fig:PersistenceDiagramsRBC500}
\end{figure}

The maximum bottleneck distance between the consecutive sampling points  is $83.5$  (Figure~\ref{fig:ConsecDist}(b) caption), while the diameter of the point cloud is only $d_B = 99.5$ (Figure~\ref{fig:DistanceMatrixPeriodicOrbit}(b)).  Therefore, we expect that significant portions of the trajectory are missing.   Indeed, Figure~\ref{fig:PersistenceDiagramsRBC500}(a) shows that there are several persistence points in $\pd_0(X,d_b)$ with a (finite) death coordinate  larger than ten.  Thus, at a length scale of $20$ (which is forty times larger than the noise threshold), the sample of the trajectory is broken into several pieces.   The largest gap between different pieces of the trajectory is $40$, as indicated by  the persistence point with coordinates $(0, 20)$. This means that the sampling rate is  far from adequate. 

The diagram $\pd_1(X,d_B)$  in Figure~\ref{fig:PersistenceDiagramsRBC500}(b)  contains a single dominant point at $(20, 32.5)$ with life span $12.5$. However, unlike in our analysis of the Kolmogorov flow in the previous section,  we cannot argue that this point corresponds to a single dominant loop along which the data is organized because of the gaps in the sampling of the orbit. As mentioned in Section~\ref{sec::LargeData}, the missing parts of the orbit could introduce loops of similar size corresponding to secondary structures.  These structures might occur due to the fact that the loop corresponding to the underlying almost-periodic dynamics might be twisted, pinched, or  bent in $\per$.  In order to obtain information about secondary structures, we require a faster sampling rate.  

We increased the sampling rate considerably and collected approximately $7\times10^5$ equally-spaced snapshots of the temperature field over four-and-a-half periods and compute the associated persistence diagrams, producing a point cloud $Y\subset \per$. The maximal distances between the consecutive frames for the increased sampling rate drop to $d_B = 4$, $d_{W^1} = 28,$ and $d_{W^2} = 6.52$.  The new value of $d_B$ is much closer to our estimate of the numerical error and it is more than 24 times smaller than the diameter of the point cloud generated from the slower sampling. Since the point cloud could only increase in diameter through increasing the sample rate, we consider this sampling rate to be satisfactory. 

\begin{figure}[t]
\centering
\subfigure[]{\includegraphics[width=1.7in]{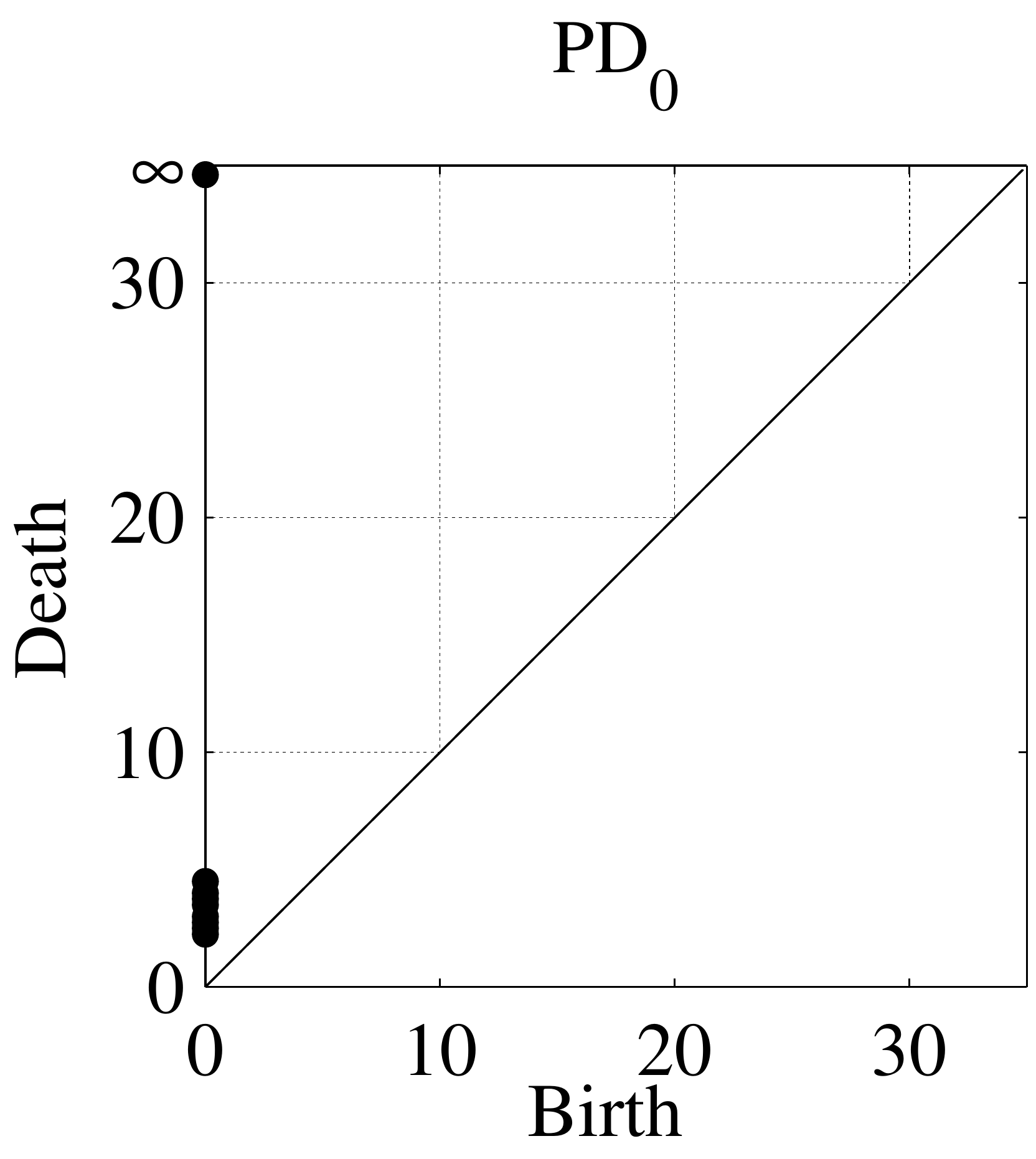}}
\subfigure[]{\includegraphics[width=1.7in]{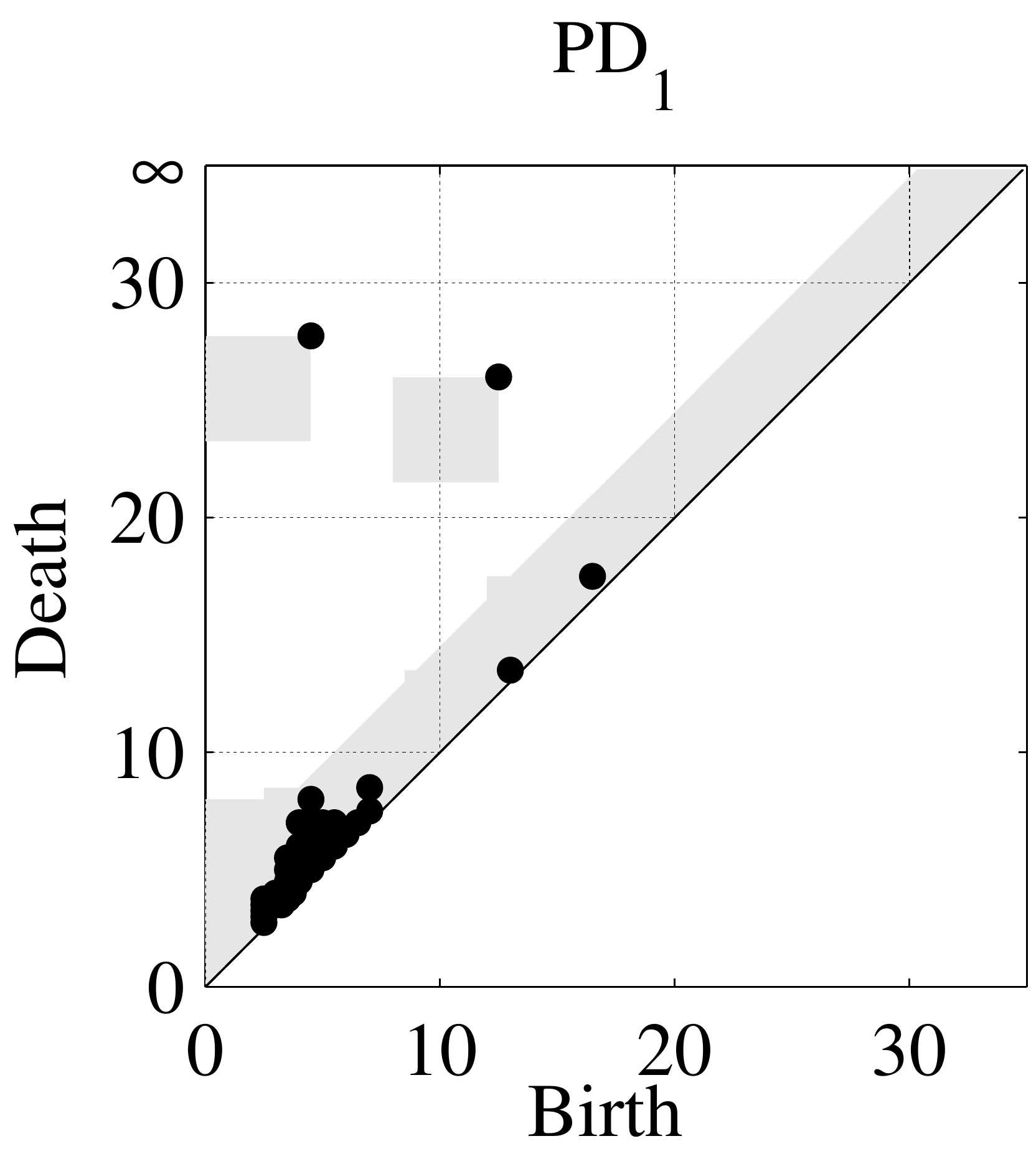}} 
\caption{Persistence diagrams for 70,000 points of Rayleigh-B\'enard convection subsampled with $\delta=4.5$, resulting in a point cloud $Y' \subset \per$ with 523 points.  (a) Persistence diagram $\pd_0(Y',d_B)$ indicating a single dominant component above the noise threshold. (b) Persistence diagram $\pd_1(Y',d_B)$ with subsampling error bounds shaded in gray.  }
\label{fig:PersistenceDiagramsRBC}
\end{figure}

Our next step is to use the ideas introduced in Section~\ref{sec::LargeData} to reduce the size of the sample and to complete our analysis.  First we construct a  $\delta$-dense, $\delta$-sparse subsample $Y' $ of the point cloud $Y$.  The smallest value of $\delta$ for which we were able to compute the persistence diagrams $\pd(Y',d_B)$, using $256$ GB of memory, is $\delta = 4.5$.  This value is only slightly larger than the largest distance between the consecutive states and,  since the diameter of the subsampled point cloud is $99.5$, the relationship between the length scale of the smallest detectable feature and the length scale of the diameter of the point cloud is still sufficient to resolve the geometry of the dynamics. The resulting persistence diagrams $\pd(Y',d_B)$ are shown in  Figure~\ref{fig:PersistenceDiagramsRBC}. 

As shown by $\pd_0(Y',d_B)$, Figure~\ref{fig:PersistenceDiagramsRBC}(a), the point cloud  merges to a single connected component at $\theta = 4.5$. This indicates that the sample of the trajectory is not broken into different pieces separated from each other. Since the maximum consecutive distance between any two points in $Y$ is $4$, the loop along which the data is organized should be present for $\theta = 2$. However, after subsampling, it is possible that the loop will not be born until $\theta = 2 + \delta$.  Looking at  the diagram $\pd_1(Y',d_B)$  in Figure~\ref{fig:PersistenceDiagramsRBC}(b), we see that it contains a dominant point at $(4.5, 27.75)$, and so the loop was indeed born before $\theta = 2 + \delta$. This is the loop along which the point cloud is organized. Now, there is another point,  $(12.5, 26) \in \pd_1(Y',d_B)$, with life span $13.5$. This point corresponds to a secondary structure of the orbit. Indeed, it can be seen from the distance matrix for the $\delta$-sparse, $\delta$-dense subsample (not shown for brevity) that the part of the orbit corresponding to the fast dynamics (missing for the slow sampling rate) revisits very similar states before continuing along the main loop. However, the development of more sensitive tools is required to fully understand these secondary features. 

We now turn our attention to the differences between the persistence diagrams of the original point cloud $Y$ and its subsample $Y'$. Theorem~\ref{thm::subsample} implies that  $d_B(\pd(Y,d_B), \pd(Y',d_B)) \leq \delta$, and so there exists a bijection between the points in $\pd(Y,d_B)$ and  $\pd(Y',d_B)$ such that the distance between matched points is less than $4.5$.  According to Remark~\ref{rem::subsample}, for the dominant point  $(4.5, 27.75) \in \pd_1(Y',d_B)$, there is exactly one corresponding point in $\pd_1(Y,d_B)$. This point is the unique point in $\pd_1(Y,d_B)$  that lies inside of the shaded box touching the point  $(4.5, 27.75)$, see Figure~\ref{fig:PersistenceDiagramsRBC}(b). The same is true for the other dominant point. Moreover, there are no points in $\pd_1(Y,d_B)$ outside of the shaded regions. Points in $\pd_1(Y,d_B)$ that do not correspond to the off-diagonal points in    $\pd_1(Y',d_B)$ can appear only $\delta/\sqrt{2} \approx 3.18$  far away from the diagonal.

\section{Conclusion}
\label{sec:Conclusion}

We have shown how persistent homology can be used to identify equilibria and study periodic dynamics, and how this method is particularly natural when solutions must be identified that lie on a group orbit.  We study two regimes in Kolmogorov flow: chaotic dynamics due to the appearance of unstable fixed points, and a periodic flow that exhibits drift in a direction of continuous symmetry.  We also study an almost-periodic orbit from Rayleigh-B\'{e}nard convection.  We solve for the unstable equilibria in the first case and sample the periodic orbits in the other two cases, and use persistent homology to project these solutions to the space of persistence diagrams.  We provide theoretical results that show this projection is stable with respect to numerical errors and discuss how the projection naturally identifies symmetry-related solutions.  We give three different metrics on the space of persistence diagrams that can be used to study pattern evolution on large versus small spatial scales, and how these metrics can be used to estimate numerical error in the space of persistence diagrams. We develop an intuition for studying dynamics in the space of persistence diagrams by looking at point clouds in two-dimensional Euclidean space, and discuss methods for determining if a continuous trajectory has been sampled densely enough to resolve the underlying dynamics, as well as mathematical methods used to address issues associated with computing on large sample sets.  We demonstrate the efficacy of these methods on Kolmogorov flow and Rayleigh-B\'{e}nard convection, comparing our methods to traditional Fourier methods where appropriate.  Our results show that the geometry of the dynamics are recovered in each case.  For Rayleigh-B\'{e}nard convection in particular, we show that the dynamics are recovered even after truncating the simulated data to an 8-bit temperature field, and so this approach is suitable for studying data collected experimentally, rather than numerically.  Also for this flow, we recover more subtle aspects of the geometry in the space of persistence diagrams.  In summary, we have shown that this method is both robust to noise and sensitive to more complicated dynamics, and that it is appropriate for studying dynamics on datasets obtained experimentally. 
Our ongoing research will further refine these tools.

\section*{Acknowledgments}

The work of MK, RL, and KM has been partially supported by NSF grants NSF-DMS-0835621, 0915019, 1125174, 1248071, and contracts from AFOSR and DARPA.  The work of JRF, BS and MFS has been partially supported by NSF grants DMS-1125302, CMMI-1234436.

\section*{Appendix A. Homology of Sets on a Torus}
\label{sec:HomologyOfTorus}

Topologically, a torus $T = S^1 \times S^1$ is a closed surface defined as the product of two circles. It can be also described as a quotient of the Cartesian plane under the identifications $(x,y) \sim (x+1,y) \sim (x,y+1)$. The homology groups of $T$ are given by
\[
H_n(T) = \begin{cases} \Z &\mbox{if } n = 0 , 2,\\ 
\Z^2 & \mbox{if } n = 1, \text{\;\;and} \\ 
0 &\mbox{otherwise}. 
\end{cases} 
\]
Intuitively, this means that $T$ has a single connected component {($n=0$)}, two independent loops {($n=1$)}, and a single cavity {($n=2$)}. For a more detailed treatment of the following material, see See Hatcher, Ch 0 for a reference to homotopies of maps, and Hatcher Ch 1 for a reference on identifying independent loops in a space, or the notion of the fundamental group.

In this section, we explain the notion of independent loops of subsets of a torus.  {By a loop, we mean a continuous path $\gamma : [0,1] \rightarrow T$ such that $\gamma(0) = \gamma(1)$. We will also be using the notion of a homotopy of loops, which can be thought of as deforming one loop continuously to another. More precisely, a homotopy of loops from a loop $\gamma_0$ to a loop $\gamma_1$ is a continuous function $F:[0,1]\times[0,1] \rightarrow T$ such that $F(s,0) = \gamma_0(s)$ and $F(s,1) = \gamma_1(s)$ for all $s \in [0,1]$, and $F(0,t) = F(1,t)$ for all $t \in [0,1]$.  Define $-\gamma(t) := \gamma(1 - t)$, which runs the loop $\gamma$ backwards in time, and define $n\gamma(t) := \gamma(\lfloor nt \rfloor)$, which traverses the loop $\gamma$ $n$ times. Finally, given two loops $\alpha_1$ and $\alpha_2$, we can form their sum $\alpha_1 + \alpha_2$ by taking a path $\delta : [0,1] \rightarrow T$ such that $\delta(0) = \alpha_1(0)$ and $\delta(1) = \alpha_2(0)$ and form the loop}
 \begin{displaymath}
   (\alpha_1 + \alpha_2)(t) = \left\{
     \begin{array}{lr}
       \alpha_1(4t) & : t \in [0,1/4]\\
       \delta(4t - 1) & : t \in [1/4,1/2]\\
       \alpha_2(4t-2) &  : t \in [1/2, 3/4]\\
       -\delta(4t-3) & : t \in [3/4, 1].
     \end{array}
   \right.
\end{displaymath} 
{Algebraically, this can be written as $\alpha_1 + \delta + \alpha_2 - \delta = \alpha_1 + \alpha_2$. We say that a loop $\gamma$ is independent of a collection of loops $\alpha_1, ..., \alpha_k$ if there does not exist a homotopy of loops from $\gamma$ to a linear combination of the loops $\alpha_1, ..., \alpha_k$. }

\begin{figure*}[]
\centering
\subfigure[$X_0$]{\includegraphics[width=1.75in]{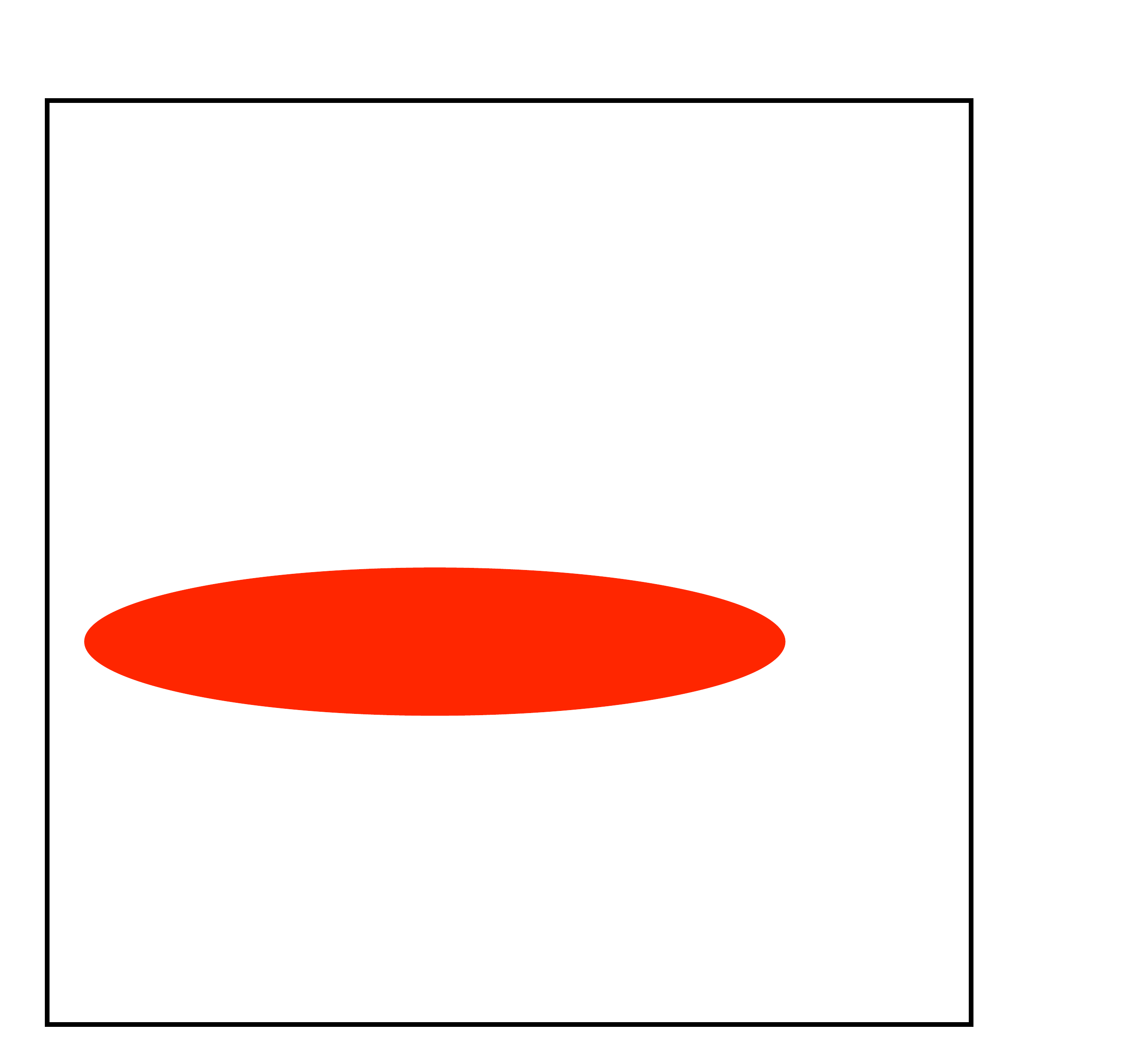}} 
\subfigure[$X_1$]{\includegraphics[width=1.75in]{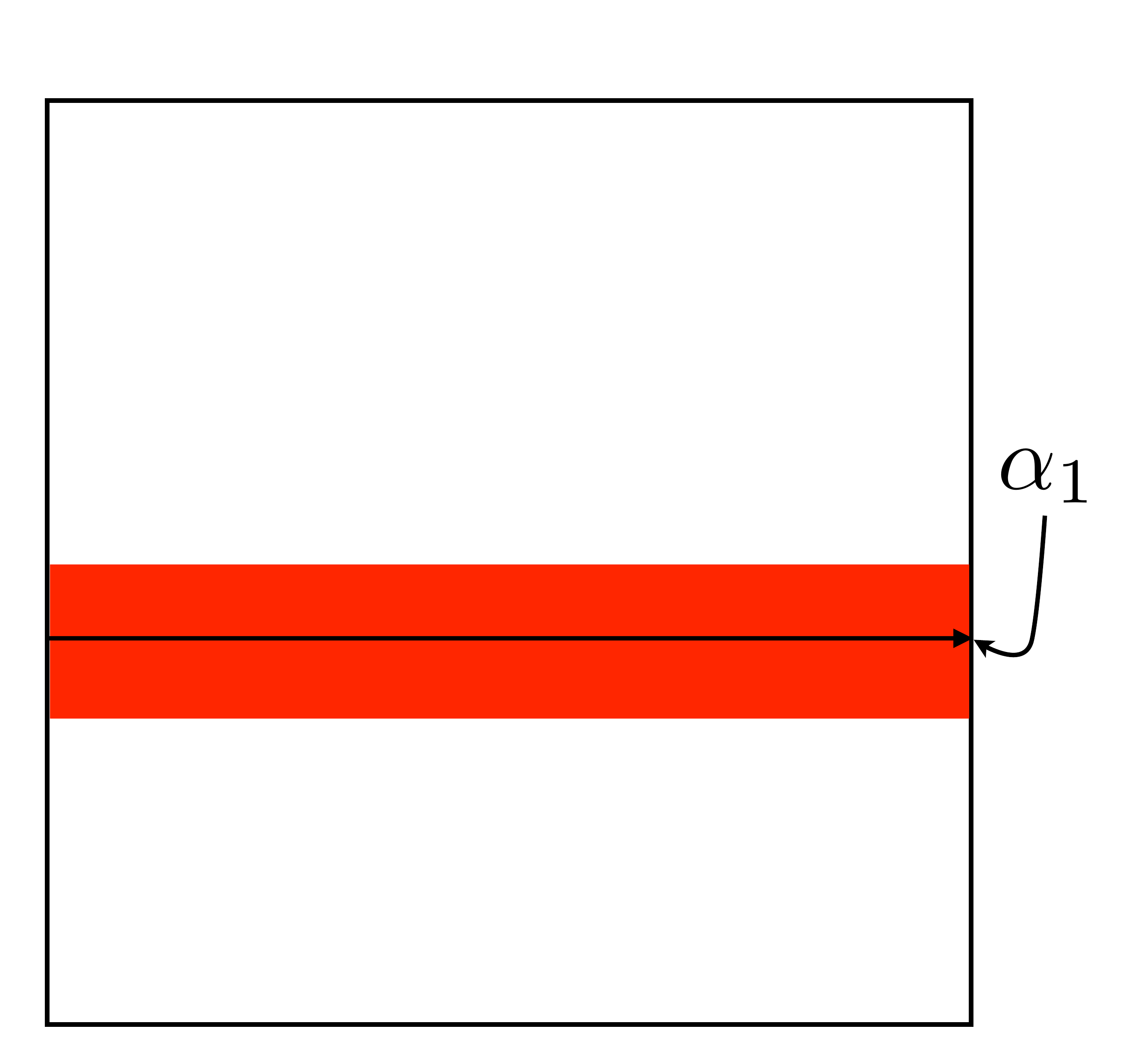}} 
\subfigure[$X_2$]{\includegraphics[width=1.75in]{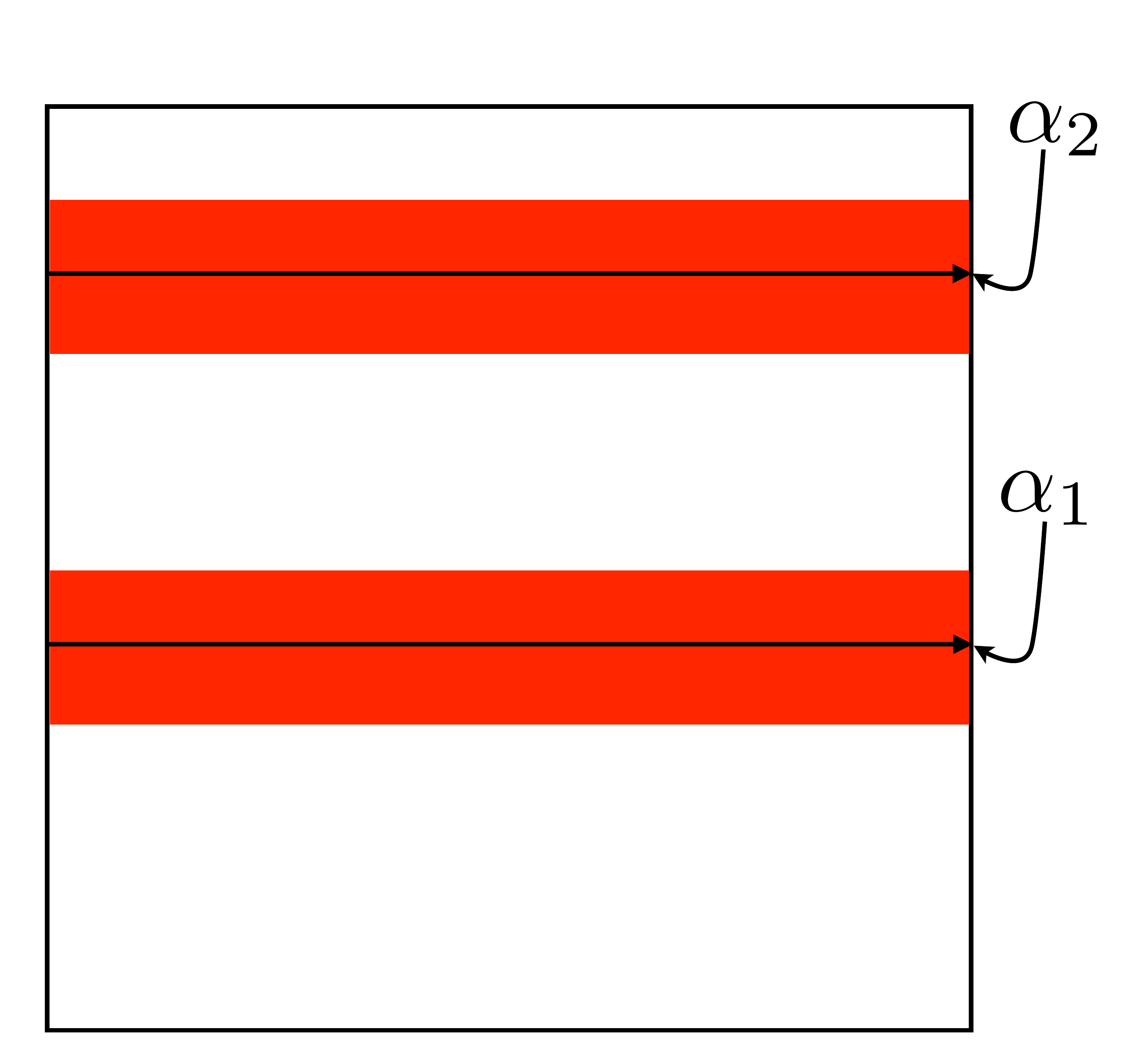}} 
\subfigure[$X_3$]{\includegraphics[width=1.75in]{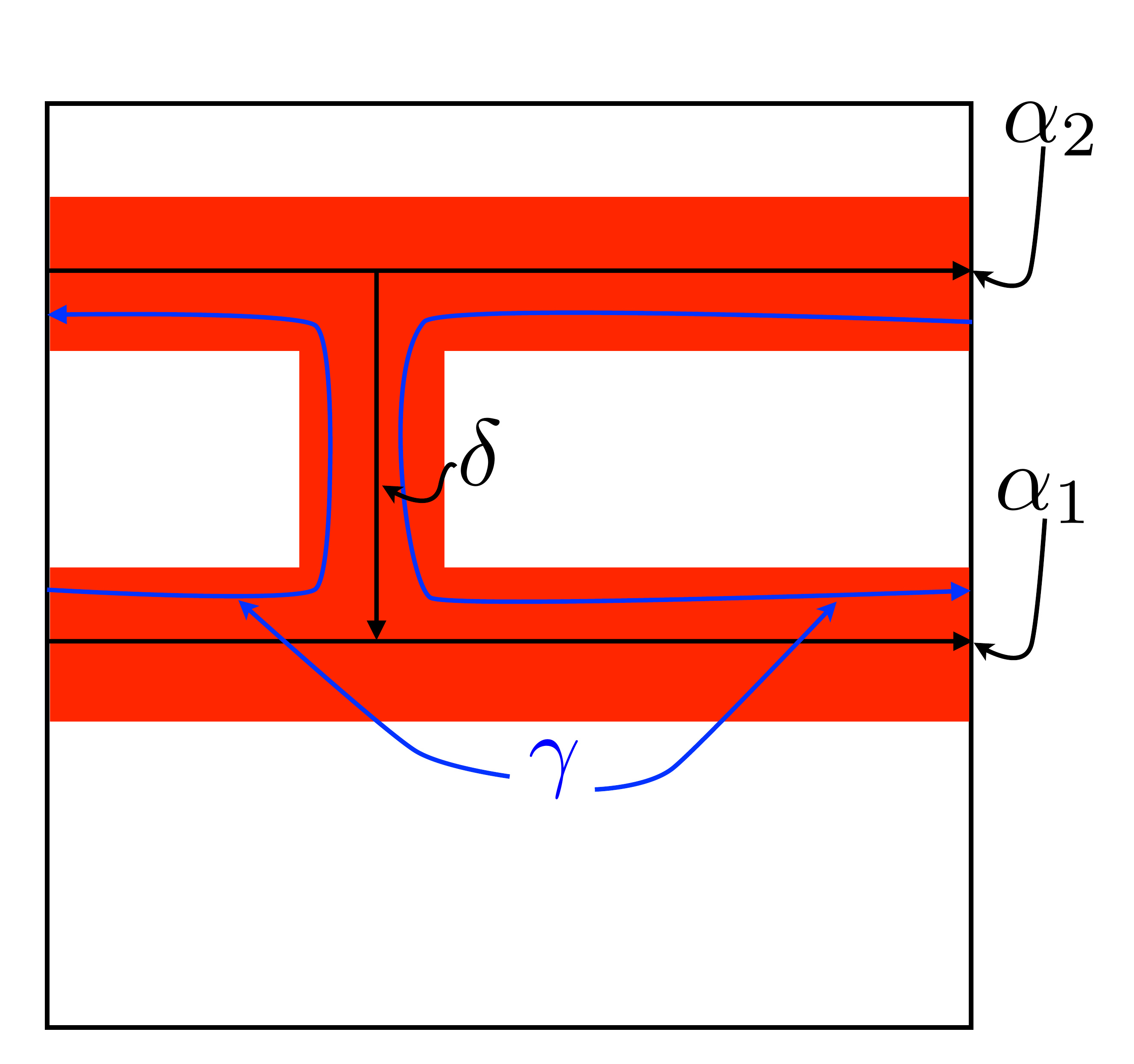}} \\
\subfigure[$X_4$]{\includegraphics[width=1.75in]{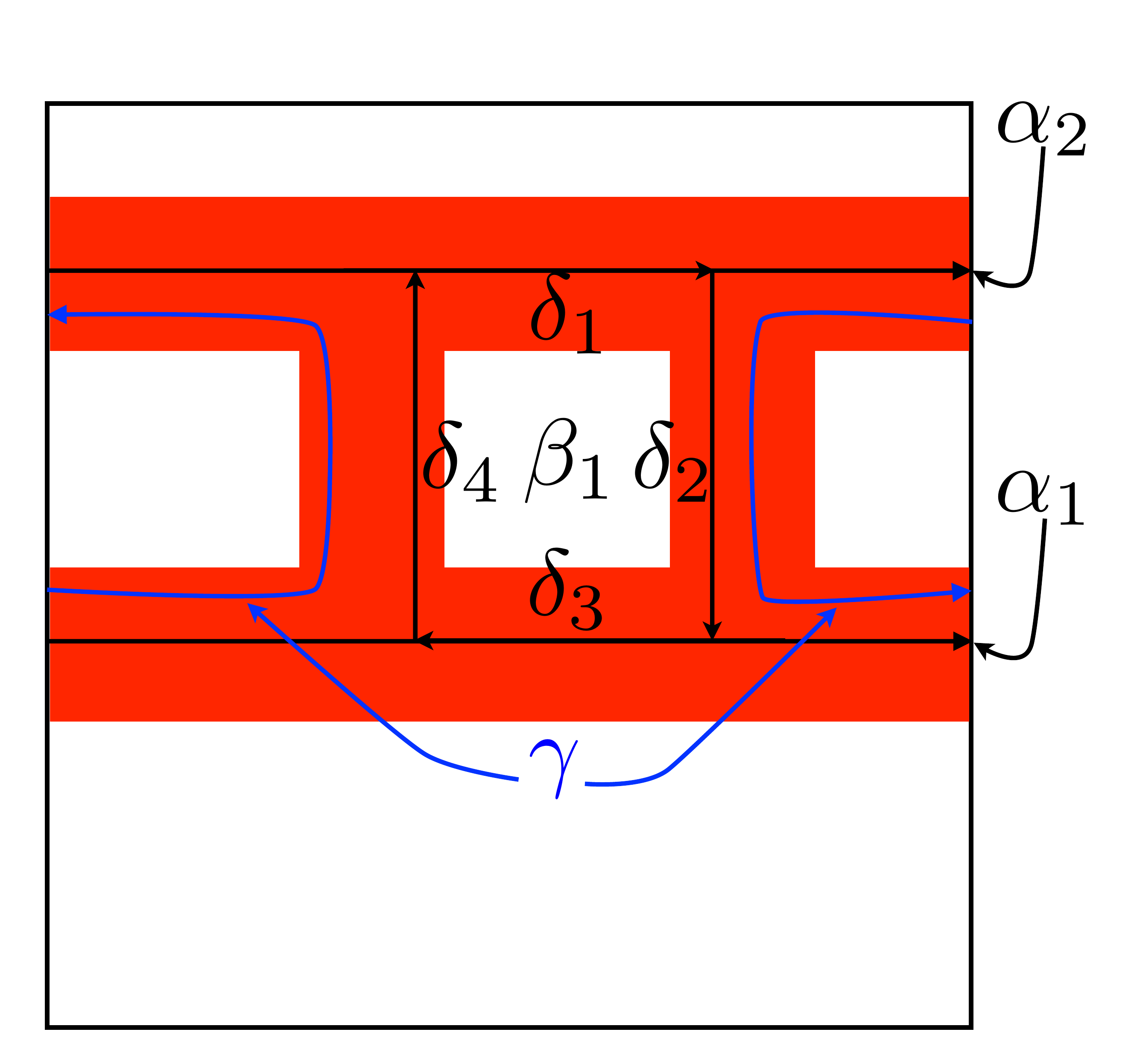}} 
\subfigure[$X_5$]{\includegraphics[width=1.75in]{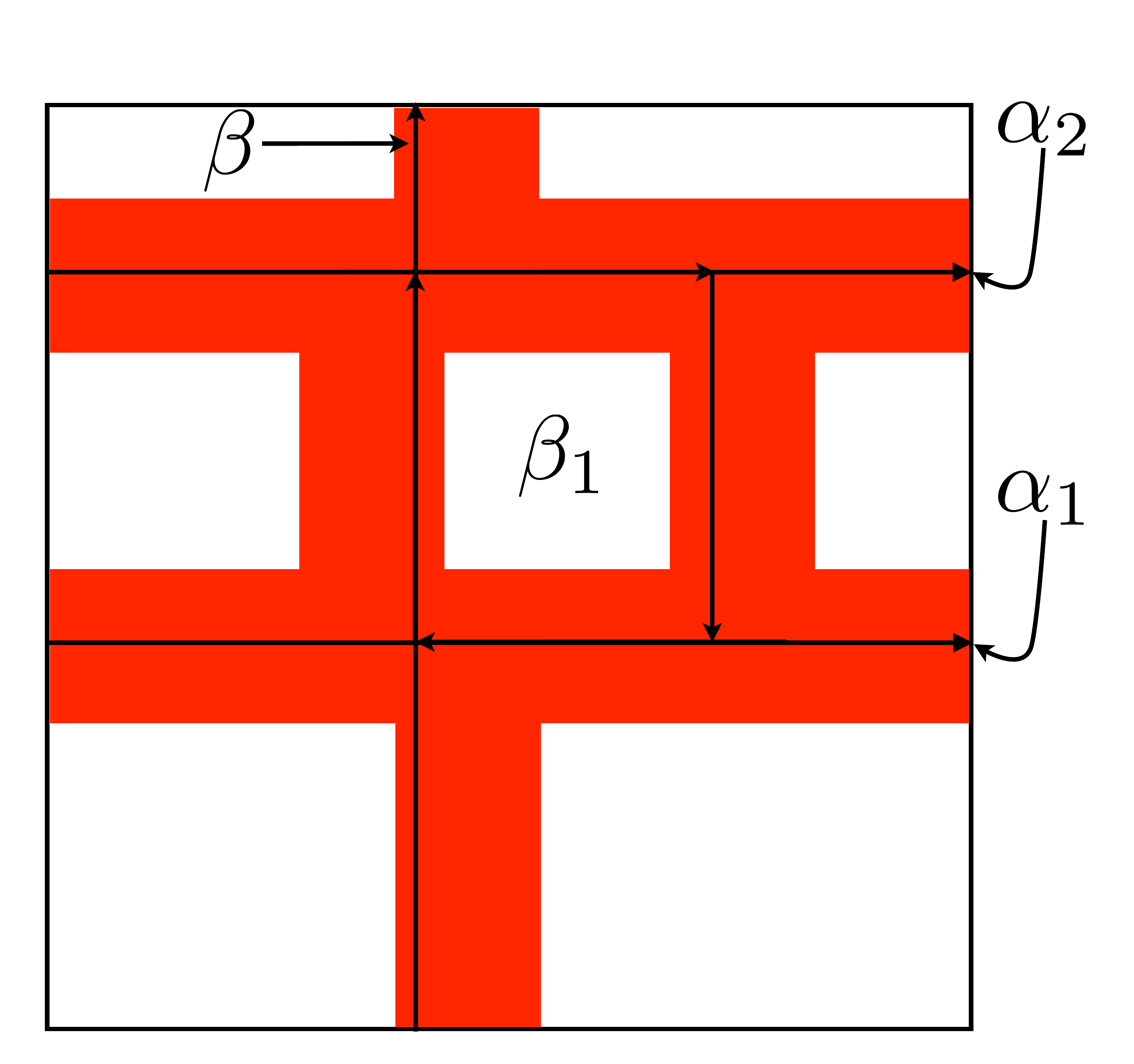}} 
\subfigure[$X_6$]{\includegraphics[width=1.75in]{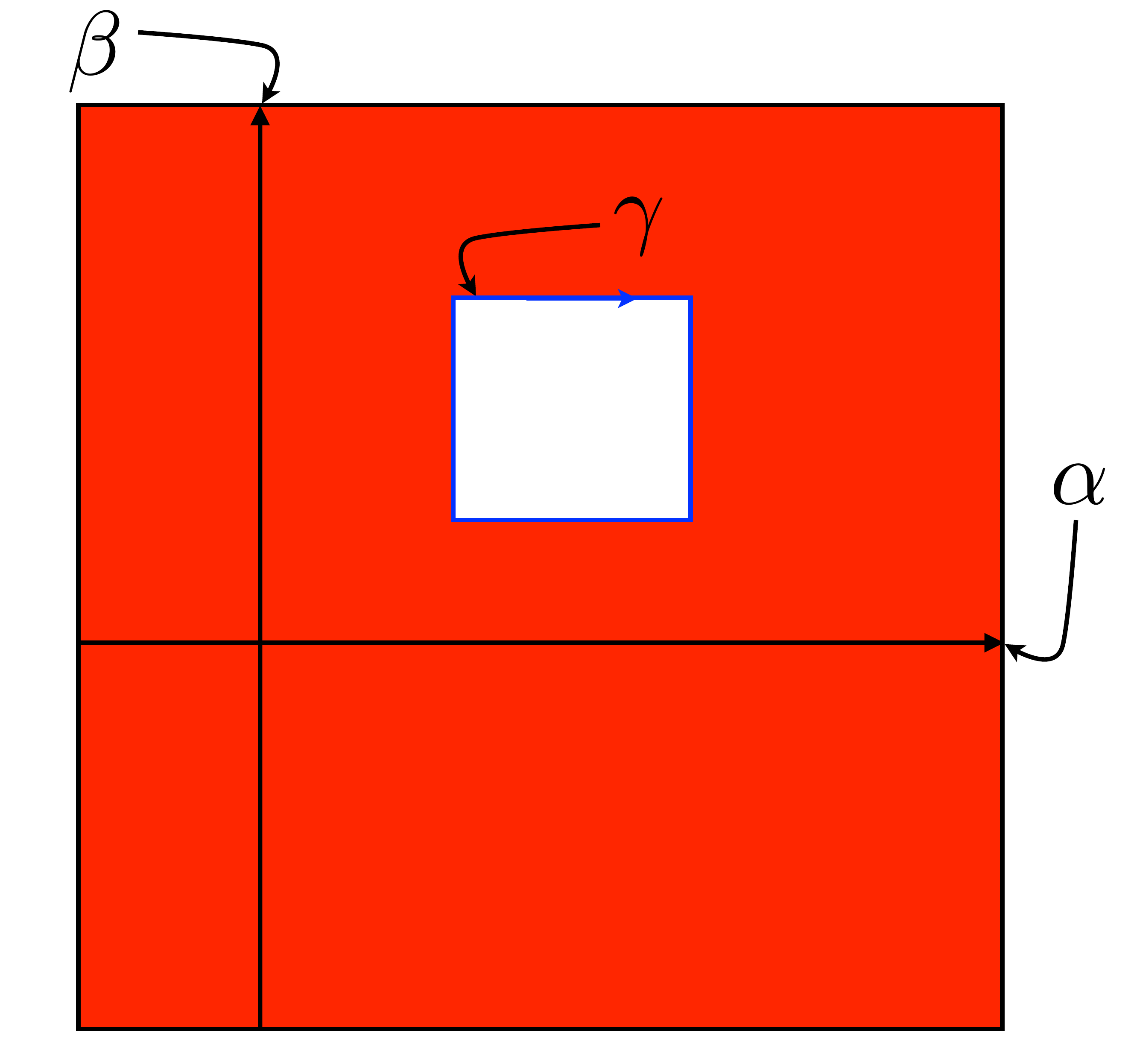}} 
\subfigure[$X_7$]{\includegraphics[width=1.75in]{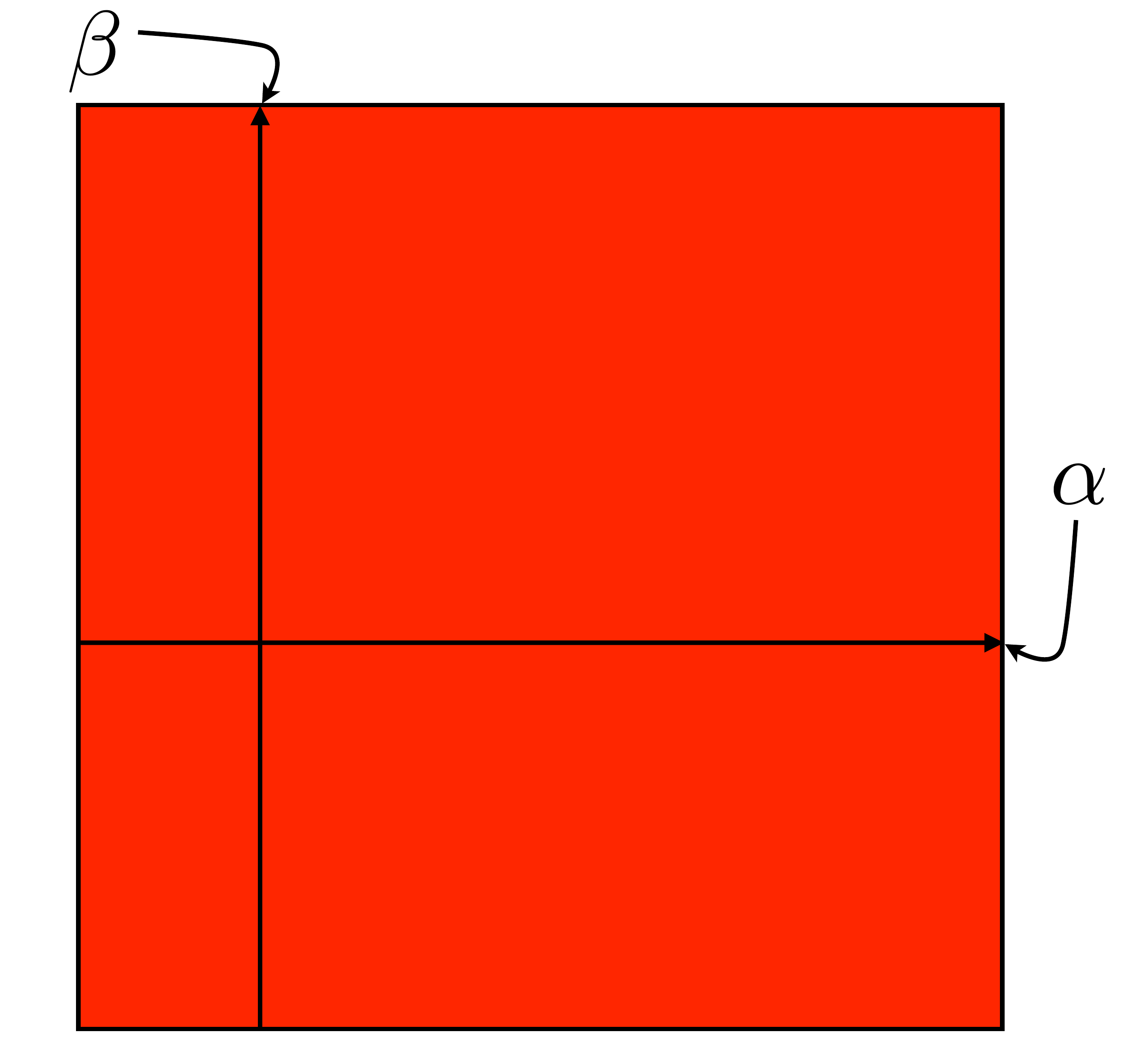}} 
\caption{Different subsets of the torus. (a) The set $X_0$ does not contain any independent loops because any loop in $X_0$ can be deformed to a point inside of $X_0$. (b) The set $X_1$ contains a single independent loop $\alpha_1$. All the other loops can be either deformed to $\alpha_1$ or to a point inside of $X_1$. (c) The set $X_2$ has two independent loops $\alpha_1$ and $\alpha_2$. (d) The loop $\gamma$ is not independent in $X_3$ because it can be deformed to the linear combination of the loops $\alpha_1$ and $\alpha_2$. (e) Adding an extra link creates a new independent loop $\beta_1$ in $X_4$. Again $\gamma$ is not independent because it can be deformed to a linear combination of the other loops. (f) The set $X_5$ is produced by adding another link which produces one new loop.    (g) The loop $\gamma$ in $X_6$ is not independent because it can be deformed to a point. (f) $X_7 = T$ contains two independent loops corresponding to two copies of $S^1$ generating the torus.}
\label{fig:HoleInTorus}
\end{figure*}

Figure~\ref{fig:HoleInTorus} shows eight subsets $\setof{X_i}_{i = 0}^7$ of a torus. Note that $X_i \subset X_j$ for $ i < j$, and the sets can be considered as sub-level sets of some scalar function $f$. We will now examine each set and identify the independent loops in each.

The set $X_0$ is contractable. Hence, every loop inside $X_0$ can be deformed  to a point inside of  $X_0$. This means that there is no independent (nontrivial) loop present in this set. 

The set $X_1$ cannot be contracted to a point. It forms a band that wraps around the torus. There are many different loops (wrapping once around the torus from left to right in the picture)  inside of this band. However, we can choose a single loop $\alpha_1$ that represents all of them; every other loop can be either continuously  deformed to a linear combination of the loop $\alpha_1$, or contracted to a point.  Similarly, the set $X_2$ contains two independent loops.   

The set $X_3$ is formed by linking the  horizontal bands present in $X_2$. The loops $\alpha_1$ and $\alpha_2$ are still independent in  $X_3$ (one cannot be deformed to the other inside $X_3$).  It might seem that there is a new independent loop, $\gamma$. However, this is not case  because $\gamma$ can be deformed (inside of $X_3$) to the union of the black lines corresponding to $\alpha_1, \alpha_2$ and $\delta$. After this deformation, the loop traverses $\delta$ twice: the right  part of the deformed loop traverses $\delta$ from the top to the bottom, and the left part in the opposite direction. Algebraically, the deformed loop can be expressed as $\alpha_1 - \delta - \alpha_2  + \delta = \alpha_1 - \alpha_2$. This shows that  $\gamma$ can be deformed to a linear combination of the loops $\alpha_1$ and $\alpha_2$. Thus, $\gamma$ is not a new independent loop.

The set $X_4$, obtained from $X_3$ by adding  another  link between the  horizontal  bands,  contains a new independent loop, $\beta_1$, consisting of the edges $\delta_1, \delta_2, \delta_3$ and $\delta_4$ ( $\beta_1 = \delta_1+ \delta_2+ \delta_3 +\delta_4$). This means that the loop $\beta_1$ cannot be deformed inside of $X_4$ to a linear combination of the loops $\alpha_1$ and $\alpha_2$ . Again, the loop $\gamma$ is not independent from $\alpha_1$, $\alpha_2$, and $\beta_1$ because it can be perturbed to $\alpha_1 + \delta_3 - \alpha_2 + \delta_4 - \alpha_2 + \delta_1 + \delta_2 = \alpha_1 - \alpha_2 + \beta_1$, which is a linear combination of the loops $\alpha_1, \alpha_2$, and $\beta_1$. Therefore, there are three independent loops in this case. Adding another link between the horizontal bands creates another independent loop. Hence, the number of independent loops for two bands with $n$ links is $n+1$.

Alternatively, we can view the set $X_3$  as a single band with one puncture, and $X_4$ as a single band with two punctures.  The number of independent loops is $n+1$, where $n$ is the number of punctures, and the extra loop is generated by the band. 

Due to the identification $(x,y) \sim (x,y+1)$, the set $X_5$ contains another link between the horizontal bands. This band creates another puncture. In Figure~\ref{fig:HoleInTorus}(f), this puncture seems to have four distinct components (white blocks in the corners). However, under the boundary identification, they correspond to  a single component. Therefore, there are four independent loops. 

The independent loops start disappearing as the punctures are filled in. The set $X_6$ contains a single puncture, and according to the previous argument, there are two independent loops, $\alpha$ and $\beta$. In this case, the loop $\gamma$ can be deformed to a point inside of the set $X_6$. Finally, the set $X_7 = T = S^1\times S^1$ contains two independent loops corresponding to the two copies of $S^1$ that generate the torus. 

\bibliographystyle{ieeetr}

\end{document}